\documentstyle[11pt]{article}

\date{}
\begin{document}
\begin{center}
{\Large \bf On Holomorphic Jet Bundles}

\vspace{0.1in}
{ by  Pit-Mann Wong and Wilhelm Stoll\\Department of
Mathematics\\ University of Notre Dame\\Notre Dame, Indianna
46556\\USA}
\end{center}

\bigskip
\noindent
{\bf Introduction}

\bigskip
In this article we provide a more detailed discussion (see
[W4]) of the jet bundles introduced by Green and Griffiths
[G-G]. In section 1 some basic facts about these jet
bundles (which are different from the usual jet bundles used
in analysis) are established with the most important one
being a Theorem of Green and Griffiths concerning the natural
filtration of the sheaf, denoted ${\cal J}_k^m
X$, of
$k$-jet differentials of weight $m$. With this reuslt many
computations (of Chern classes) and properties of ${\cal
J}_k^m X$ can be obtained or inferred from the more familiar
objects $\odot^{i_1} T^{*} X
\otimes ... \otimes \odot^{i_k} T^{*}X$ (satisfying the
condition $i_1 + 2i_2 = ... + ki_k = m$). The calculation of
Chern classes of ${\cal J}_k^m X$ are carried out in section
1 for curves and in section 2 for surfaces. These are needed
later in applying the Riemann-Roch Theorem for
${\cal J}_k^m X$ and its corresponding line sheaf, ${\cal
O}_{{\bf P}(J^kX)}(m)$, over the projectivized bundle ${\bf
P}(J^k X)$. There are some complications in working with
these sheaves due to the fact that the fibers of ${\bf
P}(J^k X)$ are weighted projective spaces and hence not
smooth and moreover, the natural sheaves ${\cal
O}_{{\bf P}(J^kX)}(m)$ are not necessarily
locally free if $m$ is not divisible by $k!$. These minor
difficulties are clarified in section 3 and is readily
seen to be rather harmless. In section 4 we consider
the case of surfaces of general type and here there is
another complication due to the fact that, as oppose to the
bundles
$\odot^{i_1} T^{*} X
\otimes ... \otimes \odot^{i_k} T^{*}X$, the sheaves
${\cal J}_k^m X$ are not semi-stable (with respect to
the canonical bundle of $X$). This difficulty, however,
can be overcome rather easily and as a result we obtain
applications in the theory of holomorphic curves in
surfaces of general type (hypersurfaces in
${\bf P}^3$ in particular). These are presented in section
4. We also include two appendices. In appendix A the lemma of
logarithmics derivatives and a version of Schwarz Lemma are
presented (see [L], [L-Y], [DSW1], [DSW2], [S-Y], [W3] and
[J]). Some combinatorics related to the symmetric groups
which we used in the computation of Chern classes (this
comes in, for example, in counting the number of positive
integer solutions of the equation $i_1 + 2i_2 + ... + ki_k =
m$) are presented in appendix B.  For higher dimensional
manifolds the approach of Nevanlinna Theory appears to work
better (see [W5]). Nevanlinna Theory for symmetric
and exterior products of the cotangent bundle can be found
in [St].

\bigskip
\noindent
{\bf \S~1 ~Holomorphic Jet Bundles}

\bigskip
We examine two concepts of "jet
bundles" of a complex manifold. The first is the
jet bundles used by analysts (PDE) and also by
Faltings in his work on rational points of an
ample subvariety of an abelian variety and
integral points of complement of an ample divsior
of an abelian variety [F]. The second is the jet
bundles introduced by Green and Griffiths [G-G].
The first notion of jet bundle shall
henceforth be referred to as the {\it full jet
bundle} while the second notion of jet bundle
shall be referred to as the {\it restricted jet
bundle}. The reason for these terminologies is that
the fiber dimension of the full jet bundle is much
larger than that of the restricted jet bundle.

\bigskip
For a complex manifold $X$ the (locally free) sheaf
of germs of holomorphic tangent vector fields
(differential operators of order 1) of
$X$ shall be denoted by $T^{1}X$ or simply $TX$.
An element of
$T^{1}X$ acts on the sheaf of germs of holomorphc
functions by differentiation:
$$(D, f) \in T^{1}X \times {\cal O}_X \mapsto Df
\in {\cal O}_X$$ and the action is linear over the
complex number field ${\bf C}$, i.e., $$D \in {\cal
H}om_{{\bf C}}({\cal O}_X, {\cal O}_X).$$ This
concept can be extended as follows:

\bigskip
\noindent
{\bf Definition 1.1}~~Let $X$ be a
complex manifold of dimension $n$ the sheaf of
germs of holomorphic $k$-jets (differential
operators of order $k$), denoted
${\cal T}^{k}X$, is the subsheaf of the sheaf of
homomorphisms
${\cal H}om_{{\bf C}}({\cal O}_X, {\cal
O}_X)$ consisting of elements (differential
operators) of the form
$$\sum_{j=1}^{k} \sum_{i_j \in {\bf N}} D_{i_1}
\circ ... \circ D_{i_{j}}$$ where
$D_{i_j} \in T^{1} X$. In terms of
holomorphic coordinates
$z_1, ..., z_n$ an element of ${\cal T}^{k}X$ is
expressed as:
$$\sum_{j=1}^{k} \sum_{1 \le i_1\le ... \le i_j
\le n} a_{i_1,..., i_j} {\partial^{j} \over
\partial z_{i_1} ... \partial z_{i_j}}$$ where
the coefficients
$a_{i_1,..., i_j}$ are holomorphic functions.
We can also drop the reqirement that the indices
be non-decreasing by requiring symmetry in the
coefficients, in other words, the elements of
${\cal T}^{k}X$ can also be expressed as:
$$\sum_{j=1}^{k} \sum_{1 \le i_1, ..., i_j
\le n} a_{i_1,..., i_j} {\partial^{j} \over
\partial z_{i_1} ... \partial z_{i_j}}$$ where
the coefficients $a_{i_1,..., i_j}$ are
symmetric in the indices $i_1,..., i_j$, i.e.,
if $\sigma$ is an element of the symmetric group
of $j$ elements then $$a_{i_{\sigma(1)},...,
i_{\sigma(j)}} = a_{i_1,..., i_j}.$$
The effect of holomorphic change of coordinates from
$z = (z_1, ..., z_n)$ to $w = (w_1, ..., w_n)$ is
given by the transistion function (for $k = 2$):
\begin{eqnarray}
 \left( \begin{array}{c}
({\partial  \over \partial z_i} )_{1 \leq i
\leq n}
 \\ ({\partial^2 \over \partial z_i \partial
z_k})_{1 \leq i \le k \leq n}
\end{array}
\right)
= \left( \begin{array}{cc}
A ~~0
 \\ B~~C
\end{array}
\right) \left( \begin{array}{c}
({\partial  \over \partial w_j} )_{1 \leq j
\leq n}
 \\ ({\partial^2 \over \partial w_j \partial
w_l})_{1 \leq j \le l \leq n}
\end{array}
\right)
\end{eqnarray} where $A$ is the $n$ by $n$
matrix:
$$A = ({\partial w_j \over \partial z_i} )_{1
\leq i,j
\leq n},$$ while $B$ is the $C^{n+1}_2$ by $n$
matrix:
$$B =
({\partial^2 w_j
\over
\partial z_i
\partial z_k})_{1 \leq i \le k \leq n}$$ $C$ is the
$C^{n+1}_2$ by $C^{n+1}_2$ matrix:
$$C = ({\partial w_j \over \partial z_i}
{\partial w_l \over \partial z_k})_{1 \leq i \le k
\leq n, 1
\leq j \le l
\leq n}
$$
and $0 = 0_{n
\times C^{n+1}_{2}}$
is the $n$ by $C^{n+1}_2$ zero-matrix
(here $C^{n+1}_{2} = (n+1)!/(n-1)! 2!$ is the usual
binomial coefficient). Note that the matrix
$A$ is the transistion function for the
tangent bundle $TX$ while the matrix
$C$ is the transistion function of $\odot^{2}
TX$, the 2-fold symmetric product of the tangent
bundle. For general $k$ the transistion function
of
$T^kX$ is of the form:
\begin{eqnarray*}
\left( \begin{array}{cccccc}
A_{1} ~~~0~~~~0~~~~0~~~~0~~~~0\\
~ *~~A_2~~~0~~~~0~~~~0~~~~0\\
~*~~~*~~~.~~~~0~~~~0~~~~0 \\
~*~~~*~~~*~~~.~~~~0~~~~0\\
~*~~~*~~~*~~*~~~.~~~~0\\
~~*~~~*~~~*~~*~~*~~A_{k}
\end{array}
\right)
\end{eqnarray*}where the $C^{n+j-1}_{j}$ by
$C^{n+j-1}_{j}$ matrix
$A_j$ is the transistion function of the bundle
$\odot^{j} TX$, the
$j$-fold symmetric product of the tangent bundle.
Here $C^{n+j-1}_{j} = (n + j - 1)!/j! (n - 1)!$
is the usual binomial coefficient.

\bigskip
The linear (and invertible) nature of the
transistion functions implies that $T^{k}X$ is
locally free. This can also be seen by observing
that
$T^{k-1}X$ injects into $T^{k} X$
and there is an exact sequence of sheaves:
\begin{eqnarray}0 \rightarrow T^{k-1}X
\rightarrow T^{k}X
\rightarrow T^{k}X /
T^{k-1}X \rightarrow 0\end{eqnarray} where
\begin{eqnarray}T^{k}X / T^{k-1}X \cong
\odot^{k} T^{1} X\end{eqnarray} is the sheaf of
germs of
$k$-fold symmetric product of $T^{1}X$, i.e.,
sheaf of germs of operators of the form:
$$\sum_{1 \le i_1\le ... \le i_j \le n} a_{i_1,..., i_j}
{\partial^{k} \over \partial x_{i_1} ... \partial
x_{i_k}}.$$ These exact sequences imply, by
induction, that ${\cal T}^{k}X$ is locally free
as the sheaves $\odot^{k} T^{1} X$, being the
symmetric product of the tangent sheaf, is
locally free. We include here the proof of the
isomorphism (3).

\bigskip
\noindent
{\bf Proposition 1.2}~~{\it With the notations
above we have:
$$ T^kX/T^{k-1}X \cong \odot^{k} TX$$
where $\odot TX$ is the $k$-fold symmetric
product of the tangent bundle.}

\medskip
\noindent
{\it Proof.}~~We shall define a surjection
from the $k$-fold tensor
product of
$TX$ onto the quotient $T^kX/T^{k-1}X$:
$$\mu :
\otimes^{k} TX \rightarrow T^kX/T^{k-1}X$$
 and then
show that the surjection factors through the
symmetric product resulting in a bijection.
The map $\mu$ is defined by:
$$\mu(D_1 \otimes ... \otimes D_k) = [D_1 \circ
...
\circ D_k]$$ where $D_i$ is (the germ of) a
vector field and $[~] : T^kX \rightarrow
T^kX/T^{k-1}X$ is the quotient map. By
definition the map $\mu$ is surjective. To see
that the map factors through to the symmetric
product it is sufficient to show that the map is
invariant by any transposition, i.e.,
$$\mu(D_1
\otimes ...
\otimes D_i \otimes D_{i+1} \otimes ... \otimes
D_k) = \mu(D_1 \otimes ... \otimes D_{i+1}
\otimes D_{i} \otimes ... \otimes D_k).$$ This
follows from the fact that the Lie bracket $D_i
\circ D_{i+1} - D_{i+1} \circ D_i$ of the vector
fields $D_i$ and $D_{i+1}$ is again a vector
field and not a 2-jet. Thus we have:
$$ D_1 \circ ... \circ (D_{i} \circ D_{i+1}
- D_{i+1} \circ D_{i}) \circ ... \circ D_k
\in T^{k-1}X$$ which implies that the map $\mu$
descends to the symmetric product $\odot^k TX$.
More precisely, if we denote the symmetrization
operator by $\sigma_k$ then
$$\overline{\mu}(D_1 \odot ... \odot D_k)
= \overline{\mu}(\sigma_k
(D_1
\otimes ...
\otimes D_k)) \stackrel{\rm def}{=} \mu (D_1
\otimes ... \otimes D_k)$$ is well-defined. It
is clear that
$\overline{\mu}$ is surjective and it remains to
show that $\overline{\mu} :  \odot^k TX
\rightarrow T^kX/T^{k-1}X$ is injective.

\bigskip
Let $(z_1, ..., z_n)$ be a local coordinate near
a point $x \in X$ then
$${\partial^k \over \partial x_{i_{1}} ...
\partial x_{i_{k}}} = {\partial \over \partial
x_{i_1}}
\odot ... \odot {\partial \over \partial
x_{i_{k}}}, ~1 \le i_1 \le ... \le i_{k} \le n$$
is a basis of $\odot^{k} TX$ at the point $x \in
X.$ If $\overline{\mu}$ is not injective then
there exists a differential operator $\Psi$ of
the form
$$\Psi = \sum_{1 \le i_1 \le ... \le i_k \le n}
a_{ i_1 ... i_k } {\partial^k \over \partial
x_{i_{1}} ... \partial x_{i_{k}}} - \Phi$$ where
$\Phi$ is a differential operator of order at
most $k - 1$ and such that $\overline{\mu}(\Psi)
= 0$.  Apply the operator $\Psi$ to the function
$f = x_{i_{1}} ... x_{i_{k}}$ shows that this is
possible only if all the coefficients $ a_{ i_1
... i_k }$ are zero. This shows that
$\overline{\mu}$ is injective and completes the
proof of the Proposition.~~QED

\bigskip
The restricted $k$-jet bundles are introduced
by Green-Griffiths in [G-G]. It is defined as
follows. Denote by ${\cal H}_{x}, x \in X$ the
sheaf of germs of holomorphic curves
$\{f : \Delta_{r} \rightarrow X, f(0) = x\}$.
Define, for
$k \in {\bf N}$, an equivalence relation
as follows. Let $z_1, ..., z_n$ be holomorphic
coordinates near $x$ and for $f \in {\cal H}_x$
let $f_i = z_i \circ f$. Two elements $f, g \in
{\cal H}_x$ are said to be $k$-equivalent,
denoted $f
\sim_k g$, if $f_{j}^{(p)}(0) =
g_{j}^{(p)}(0)$ for all $1 \leq p \le k$. The
sheaf of restricted $k$-jets is defined to be
$J^{k}X = \cup_{x \in X} {\cal H}_x /
\sim_k$. Elements of $J^{k}X$ will be
denoted by $j^{k}f(0) = (f(0), f^{'}(0), ...,
f^{(k)}(0))$. It is also clear that $J^{1}X =
T^{1}X = TX$ is the tangent bundle.

\bigskip
Note that the definition above does not depend on the choices of the
coordinates
near $x$; for if $(z_{j} \circ f)^{(p)}(0) =
(z_{j} \circ g)^{(p)}(0)$ for $0 \le
p \le k$ then $(w_{j} \circ
f)^{(p)}(0) = (w_{j} \circ g)^{(p)}(0)$ for any
other coordinates
$w_1, ..., w_n$. The effect of change
of coordinates is as follows:
\begin{eqnarray*}
&&(w_j \circ f)^{'} = \sum_{i=1}^{n} {\partial
w_j \over \partial z_{i}}(f) (z_i \circ f)^{'},\cr
&&(w_j \circ f)^{''} = \sum_{i=1}^{n} {\partial
w_j \over \partial z_{i}}(f) (z_i \circ f)^{''} +
\sum_{i,k=1}^{n} {\partial^2 w_j \over
\partial z_i \partial z_k}(f) (z_i \circ f)^{'}
(z_k \circ f)^{'}
\end{eqnarray*}
and for general $k$ we have
\begin{eqnarray}
(w_j \circ f)^{(k)} = \sum_{i=1}^{n} {\partial
w_j \over \partial z_{i}}(f) (z_i \circ f)^{(k)} +
P({\partial^l w_j \over
\partial z_{i_1} ... \partial z_{i_l}}(f), (w_j
\circ f)^{(l)})
\end{eqnarray}where $P$ is a polynomial with
integer coefficients in $\partial^l w_j /
\partial z_{i_1} ... \partial z_{i_l}, (w_j
\circ f)^{(l)}$ for $j = 1, ..., n$ and $l
= 1, ..., k$.

\bigskip
Note that the quadratic nature in
$(z_{i} \circ f)^{'}$ in the formula for
$(w_{i} \circ f)^{''}$ means that the sheaf
$J^kX$ is not locally free. It is instructive to
compare this with the matrix
$C$ in the transition formula (1) for $T^{k}X$
under the change of coordinates. In that case the
formula is quadratic in the partial derivatives
$\partial w_j/\partial z_i$ but linear in $(z_i
\circ f)^{'}$ hence the transformation can still
be represented as a linear transformation while
this is not the case for
$J^{k}X$. There is however, a natural
${\bf C}^{*}$-action on $J^{k}X$ defined via
parametrization. Namely, for $\lambda
\in {\bf C}^{*}$ and $f \in {\cal H}_{x}$
a map $f_{\lambda} \in {\cal H}_{x}$ is defined by
$ f_{\lambda}(t) = f(\lambda t)$; then
$j^{k}f_{\lambda}(0) = (f_{\lambda}(0),
f_{\lambda}^{'}(0), ..., f_{\lambda}^{(k)}(0))
= (f(0),
\lambda f^{'}(0), ..., \lambda^{k} f^{(k)}(0)).$
 In other words the $C^{*}$-action is given by
\begin{eqnarray}\lambda . j^{k}f(0) = (f(0),
\lambda f^{'}(0), ..., \lambda^{k}
f^{(k)}(0)).\end{eqnarray}

\bigskip
\noindent
{\bf Definition 1.3}~~The restricted
$k$-jet bundle is defined to be $J^{k}X$
together with the ${\bf C}^{*}$-action defined
above and shall simply be denoted by
$J^{k}X$.

\bigskip
Another difference between the full and restricted
$k$-jet bundles is that there is, in general, no
natural way of injecting
$J^{k-1}X$ into $J^{k}X$. For instance, the
coordinates transformations shows that
$$(f(0), (f_1^{'}(0), ..., f_n^{'}(0))) \mapsto
(f(0), (f_1^{'}(0), ..., f_n^{'}(0)), (0, ...,
0))$$ is not a well-defined map of $J^{1}X$ into
$J^{2}X$ as the condition $f^{''}(0) = 0$ is not
preserved by a general change of coordinates.
On the other hand, the transformation formulas
show that $\{j^{2}f = (f(0),
f^{'}(0), f^{''}(0))
\in J^{2} X~|~f^{'}(0) = 0\}$, more generally,
\begin{eqnarray}Z_{0} = \{j^{k}f = (f(0),
f^{'}(0), ..., f^{(k)}(0))
\in J^{k} X~|~f^{'}(0) = 0\}\end{eqnarray} is a
well-defined subvariety of $J^{k}X$ as the
conditon
$f^{'}(0)$ is invariant under change of
coordinates. Moreover, the transformation law
actually says that even though the condition
$f^{''}(0) = 0$ is coordinate dependent
the conditions that $f^{'}(0) = f^{''}(0) = 0$ are
independent of choices of coordinates,
in other words, the zero-section of
$J^2X$,
 more generally, the
zero-section of
$J^kX$:
\begin{eqnarray}\{j^kf(0) \in J^kX~|~f^{'}(0)
= f^{''}(0) = ... = f^{(k)}(0) = 0\}\end{eqnarray}
is well-defined.

\bigskip
\noindent
{\bf Theorem 1.4}~~{\it Let $X$ be a complex
manifold of dimension $n$ then $T^kX$ is a
holomorphic vector bundle of rank $r = n +
C^{n+1}_{2} + C^{n+2}_{3} + ... + C^{n+k
-1}_{k} = \sum_{i=1}^{k} C_{i}^{n+i-1}$ while
$J^{k}X$ is a holomorphic
${\bf C}^{*}$-bundle of rank $r = kn$ and the
zero-section of $J^{k}X$ is well-defined.}

\bigskip
As noted above there is no natural
inclusion map from $J^{k-1}X$ into $J^{k}X$ there is
however a natural projection map
$$p_{kj} : J^{k}X \rightarrow
J^{j}X$$ for any
$j \le k$ defined simply by
\begin{eqnarray} p_{kj}(j^{k}f(0)) = j^{j}f(0).
\end{eqnarray} The projection map clearly respect
the ${\bf C}^{*}$-action defined by (5) and so is
a
${\bf C}^{*}$-bundle morphism.

\bigskip
If $\Phi : X \rightarrow Y$ is a holomorphic map
between the complex manifolds $X$ and $Y$ then the
usual differentail
$\Phi_{*} : T^{1}X \rightarrow T^{1}Y$ is defined.
The same is true for the $k$-jets as the
$k$-th order differential $\Phi_{k*} : T^{k}X
\rightarrow T^{k}Y$ can be defined by
$$
\Phi_{k*} = (D_1 \circ ... \circ D_{k}) (g)
\stackrel{\rm def}{=} D_1 \circ ... \circ D_{k} (g
\circ \Phi)
$$where $g \in {\cal O}_{Y}$. The $k$-th
order differential, denoted $J^{k} \Phi : J^{k}X
\rightarrow J^{k}Y$ can also be defined:
$$
J^{k} \Phi(j^{k}f(0))
\stackrel{\rm def}{=} (\Phi \circ f)^{(k)}(0)
$$ for any $j^{k}f(0) \in J^{k}X$. For
the restricted jet bundle $J^{k}X$ there is another
notion closely related to (but not the same) the
differential: the natural lifting of a
holomorphic curve. Namely, given any holomorphic map
$f : \Delta_{r} \rightarrow X (0 < r \le \infty$),
the lifting $j^{k}f :  \Delta_{r/2} \rightarrow
J^{k}X$ is defined by:
$$
j^{k}f(\zeta) = j^{k}g(0),~\zeta \in \Delta_{r/2}
$$ where $g(\xi) = f(\zeta + \xi)$ is
holomorphic for $\xi \in \Delta_{r/2}$.

\bigskip
Consider the special case $\dim X = 1$ then
$T^{k}X$ and $J^{k}X$ have the same rank and the
underlying space of
$T^{k}X$ and
$J^{k}X$ are the same but the structures are
different. Consider the map (for
simplicity we write this out only for $k = 2$):
\begin{eqnarray}(f(0), f^{'}(0),
f^{('')}(0)) \mapsto f^{''}(0) {\partial \over
\partial z} + (f^{'}(0))^{2}
{\partial^2
\over \partial z^2}.\end{eqnarray} which is
clearly holomorphic but is not
biholomorphic. For if $$g(t) = f(0) - f^{'}(0) t
+ f^{''}(0) t^{2}/2$$ then
$j^{2}g(0) = (f(0), - f^{'}(0), f^{''}(0))$
and under the identification above
$j^{2}f(0)$ and $j^{2}g(0)$ are mapped onto
the same element. Moreover the map is a ${\bf
C}^{*}$-bundle map because $\lambda.(f(0),
f^{'}(0), f^{('')}(0)) = (f(0),
\lambda f^{'}(0), \lambda^{2} f^{('')}(0))$
is mapped onto the element
$$\lambda^{2} \{f^{''}(0) {\partial \over
\partial z} + (f^{'}(0))^{2}
{\partial^2
\over
\partial z^2}\}.$$
More generally, for $X$ of arbitrary dimension,
the map
\begin{eqnarray}(f(0), f^{'}(0),
f^{('')}(0)) \mapsto \sum_{i=1}^{n} f_{i}^{''}(0)
{\partial
\over
\partial z_{i}} + \sum_{1 \le i,
j \le n} f_{i}^{'}(0) f_{j}^{'}(0) {\partial^2
\over
\partial z_{i} \partial z_{j}}\end{eqnarray} is a holomorphic ${\bf
C}^{*}$-bundle
map from $J^{2}X$ onto a ${\bf C}^{*}$ sub-bundle of $T^{2}X$. We have already
seen the case of $n = 1$; for $n = 2$ the
second sum above has 3 terms:
$$(f_{1}^{'})^{2}, (f_{2}^{'})^{2},
f_{1}^{'} f_{2}^{'}.$$ Thus if $j^{2} f(0)$ and
$j^{2}g(0)$ have the same image then
$$(f_{1}^{'})^{2} =
(g_{1}^{'})^{2}, (f_{2}^{'})^{2} =
(g_{2}^{'})^{2}, f_{1}^{'} f_{2}^{'} = g_{1}^{'}
g_{2}^{'}$$ so that $f_{1}^{'} = \pm g_{1}^{'},
f_{2}^{'} =
\pm g_{2}^{'}$. This means that the map
is generically 2 to 1 onto its image and ramified
along the subvariety $Z_{0}$ defined by (6).

\bigskip
Returning to the case of a Riemann surface $X$
we define a map $p_{3} : J^{3}X
\rightarrow T^{3}X$ by the formula:
$$p_3(j^{3}f(0)) = f^{(3)}(0)
{\partial
\over
\partial z} + f^{''}(0) f^{'}(0)
{\partial^{2}
\over
\partial z^{2}} + (f^{'}(0))^{3}
{\partial^{3}
\over
\partial z^{3}}$$ and in general
$p_{k} : J^{k}X
\rightarrow T^{k}X$ by the formula:
$$p_k(j^{k}f(0)) =
\sum_{j=1}^{k} f^{(j)}(0)
(f^{'}(0))^{k-j} {\partial^{k-j+1}
\over
\partial z^{k-j+1}}.$$
For the higher dimensional manifold $X$ the maps
are defined by
\begin{eqnarray*}p_3(j^{3}f(0)) &=&
\sum_{i=1}^{n} f_{i}^{(3)}(0) {\partial
\over
\partial z_{i}} +
\sum_{i,j=1}^{n} f_{i}^{''}(0) f_{j}^{'}(0)
{\partial^{2}
\over
\partial z_i \partial z_j} + \cr
&&~~+ \sum_{i,j,k=1}^{n} f_i^{'}(0)
f_j^{'}(0) f_k^{'}(0) {\partial^{3}
\over
\partial z_{i} \partial z_j \partial
z_k}\end{eqnarray*} and in general
$p_{k} : J^{k}X
\rightarrow T^{k}X$ by the formula:
\begin{eqnarray}p_k(j^{k}f(0))
=
\sum_{j=1}^{k} \sum_{i=1}^{n}
\{f_i^{(j)}(0) \prod_{i_{1}, ...,
i_{k-j} \ne i} f_{i_{l}}^{'}(0)\}
{\partial^{k-j+1}
\over
\partial z_{i} \partial z_{i_1} ... \partial
z_{i_{k-j}}}.\end{eqnarray} It is clear from the
definition of the map $p_{k}$ that:

\bigskip
\noindent
{\bf Theorem 1.5}~~{\it Let $J^{k}X$ and $T^{k}X$
be, respectively, the restriced and the full
$k$-jet bundles over a complex manifold $X$. Then
the map defined by
$(10)$ is a holomorphic ${\bf C}^{*}$-bundle
map which is generically finite to $1$ onto its
image. Moreover, the map is ramified precisly
along the subvariety $Z_{0} =
\{j^{k}f(0)
\in J^{k} X~|~f^{'}(0) = 0\}.$}

\bigskip
We consider now the "dual" of the jet bundles.

\bigskip
\noindent
{\bf Definition 1.6}~~The dual of the full jet
bundles $T^{k}X$ shall be referred to as the sheaf of
germs of $k$-jet forms and shall be denoted by
$T_k^{*}X$. The global sections shall be
referred to as $k$-jet forms. For $m \in {\bf
N}$ the $m$-fold symmetric product shall be
denoted by either $\odot^m T_k^{*}X$ and its
global sections shall be referred to as $k$-jet
forms of weight $m$.

\bigskip
By definition, a $k$-jet form of weight $m$
assigns to each point $x \in X$ a homogeneous
(with respect to the standard ${\bf C}^{*}$-action
of
$T^kX$ as a vector bundle) polynomial of degree
$m$ on the fiber
$T_x^kX$ (where $T^kX$ is the $k$-jet bundle).
Let $(U, z_1, ..., z_n)$ be a local holomorphic
coordinates over $U$ then
\begin{eqnarray*}&&(e_i = {\partial \over
\partial z_{i}})_{1 \le i
\le n},\\&& (e_{i_{1} i_{2}} = {\partial^2 \over
\partial z_{i_{1}}
\partial z_{i_{2}}})_{1
\le i_{1} \le i_{2} \le n},\\&&~~~
.\\&&~~~.\\&&~~~.\\ &&(e_{i_{1} ... i_{k}} =
{\partial^k
\over
\partial z_{i_{1}} ...
\partial z_{{i}_k}})_{1
\le i_1 \le i_2 \le ... \le i_k \le
n}\end{eqnarray*} is a basis of $T^kX|_U$. The
dual basis shall be denoted, formally, by
\begin{eqnarray*}&&(e_i^{*} = dz_i)_{1 \le i
\le n},\\ &&(e_{i_{1} i_{2}}^{*} = d^2 z_{i_{1}}
z_{i_{2}})_{1
\le i_{1} \le i_{2} \le n},\\&&~~~
.\\&&~~~.\\&&~~~.\\ &&(e_{i_{1} ... i_{k}}^{*} =
d^k z_{i_{1}} ... z_{i_{k}})_{1
\le i_1 \le i_2 \le ... \le i_k \le
n}.\end{eqnarray*} An element of
$T_k^{*}X$ is then of the form:
\begin{eqnarray*}
\omega = \sum_{j=1}^{k} \sum_{1 \le i_1 \le ...
\le i_j \le n} a_{i_1, ..., i_j} e_{i_1, ...,
i_j}^{*}
\end{eqnarray*}where the coefficients $a_{i_1,
..., i_j}$ are holomorphic functions. Sometimes
it is convenient to express the sum without the
restriction as in the second sum above but
insisting on the symmetry of the coefficients
(see also definition 1.1):
\begin{eqnarray*}
\omega = \sum_{j=1}^{k} \sum_{1 \le i_1, ..., i_j
\le n} a_{i_1, ..., i_j} e_{i_1, ..., i_j}^{*}
\end{eqnarray*}where the coefficients $a_{i_1,
..., i_j}$ are symmetric in the indices.
 We can also write down a basis for the
symmetric product
$\odot^m T^kX$ and its dual basis for $\odot^m
T_k^*X$. It is convenient to use the following
notations and conventions for the index set. Let
$${\cal I}_{k} = \{I = (i_1, ...,
i_k)~|~i_{k} \in {\bf N} \cup \{0\} , 0 \le
i_{1} \le ... \le i_{k} \le n~{\rm
and~not~all}~i_j = 0\}$$ be endowed with
lexicographical order and let
$${\cal J}_m({\cal I}_{k}) = \{J = (I_1,
..., I_m)~|~I_j \in {\cal I}_{k}, I_1 \le ... \le
I_m\}.$$ With these notations, for example, the
basis for $T_k^*X$ over $U$ is simply expressed
as
$B_k^{*} = \{e_{I}^{*} = e_{i_1}^{*} \odot ...
\odot e_{i_{k}}^{*}~|~I
\in {\cal I}_{k}\}$ with the conventions that
$e_{0} = 1$. Analogously, a basis for $\odot^m
T_k^*X$ over $U$ is expressed as $B_k^{*m} =
\{e_{J}^{*} = e_{I_1}^{*} \odot ... \odot
e_{I_{k}}^{*}~|~J
\in {\cal J}_m({\cal I}_{k})\}$. Moreover, a
section $\omega
\in H^{0}(U, \odot^m
T_k^*X)$ is expressed as
\begin{eqnarray*}\omega = \sum_{J \in {\cal
J}_m({\cal I}_{k})} a_{J} e_J^{*}\end{eqnarray*}
where the coefficients are holomorphic functins
on $U$.

\bigskip
Taking the dual of the sequence (2) we get an
exact sequence:
\begin{eqnarray}0 \rightarrow \odot^{k} T_1^{*}X
\rightarrow T_{k}^{*} X
\rightarrow T_{k-1}^{*}X \rightarrow
0.\end{eqnarray}

\bigskip
For example, for $k = 3$ we have two exact
sequences:
$$0 \rightarrow \odot^{3} T_{1}^{*}X
\rightarrow T_{3}^{*} X
\rightarrow T_{2}^{*} X \rightarrow 0,$$
$$0 \rightarrow \odot^{2} T_{1}^{*}X
\rightarrow T_{2}^{*} X
\rightarrow T_{1}^{*} X \rightarrow 0.$$ In
particular, by Whitney's Formula:
$$c_1(T_{3}^{*} X) = c_1(T_{2}^{*} X) +
c_1(\odot^{3} T_{1}^{*}X) = c_1(T_{1}^{*} X)
+ c_1(\odot^{2} T_{1}^{*}X) + c_1(\odot^{3}
T_{1}^{*}X).$$ In general, we have, by induction:

\bigskip
\noindent
{\bf Theorem 1.7}~~{\it The first Chern number
of the bundle of $k$-jet forms is given by the
formula:
$$c_1(T_{k}^{*} X) = \sum_{j=1}^{k}
c_1(\odot^{j} T_{1}^{*}X).$$ In particular, if
$X$ is a Riemann surface then
$$c_1(T_{k}^{*} X) = \sum_{j=1}^{k}
jc_1(T_{1}^{*}X) = {k(k+1) \over 2}
c_{1}({\cal K}_{X})$$ where ${\cal K}_{X} =
T_{1}^{*}X$ is the canonical bundle of $X$.}

\bigskip
Note that if $X$ is a Riemann surface then the
rank of $T_{k}^{*} X$ is $k$.

\bigskip
\noindent
{\bf Corollary 1.8}~~{\it Let $X$ be a
projective manifold and suppose that the
cotangent bundle
$T_{1}^{*}X$ is ample then $T_{k}^{*} X$ is
ample for all $k$.}

\bigskip
\noindent
{\bf Definition 1.9}~~The dual of
$J^{k}X$, i.e., germs of $\omega : j^{k}X|_{U}
\rightarrow {\bf
C}$ such that $\omega(\lambda.j^{k}f) =
\lambda^{m}
\omega(j^{k}f)$ for some positive integer $m$,
shall be referred to as the sheaf of germs of
$k$-jet differentials and shall be denoted by
${\cal J}_{k}^{*}X$. A jet differential
$\omega$ satisfying the homogenity above with
integer $m$ is said to a $k$-jet differential
of weight $m$. The sheaf of $k$-jet differential
of weight $m$ shall be denoted by
${\cal J}_{k}^{m}X$.

\bigskip
It follows from the definition of the ${\bf
C}^{*}$-action on $J^{k}X$ that a
$k$-jet differential
$\omega$ of weight
$m$ is of the form:
\begin{eqnarray}
\omega(j^kf) = \sum_{|I_1| + 2|I_2| + ... +
k|I_k| = m} a_{I_1, ..., I_k}(z) (f^{'})^{I_1}
... (f^{(k)})^{I_k}
\end{eqnarray}where $a_{I_1, ... I_k}$ are
holomorphic functions, $I_{j} =
(i_{1j}, ..., i_{nj}), n = \dim X$ are the
multi-indices with ech
$i_{lj}$ being a non-negative integer and
$I_{j}| = i_{1j} + ... + i_{nj}$. In terms of a
local coordinate $(z_{1}, ..., z_{n})$ near a
point $z$,
$$(f^{'})^{I_1} ...
(f^{(k)})^{I_k} = (f_1^{'})^{i_{11}} ...
(f_n^{'})^{i_{n1}} ... (f_1^{(k)})^{i_{1k}} ...
(f_n^{(k)})^{i_{nk}}.$$ Moreover the
coefficients $a_{I_1, ... I_k}(z)$ are symmetric
with respect to the indices in each $I_j$. More
precisely,
$$a_{(i_{\sigma_{1}(1)1}, ...,
i_{\sigma_{1}(n)1}), ..., (i_{\sigma_{k}(1)k},
..., i_{\sigma_{k}(n)k})} = a_{(i_{11}, ...,
i_{n1}), ..., (i_{1k},
..., i_{nk})}$$where each $\sigma_j : \{1, ...,
n\}
\rightarrow \{1, ..., n\}, j = 1, ..., n$ is a
permutation of $n$-elements.

\bigskip
Let ${\cal L}_{k}^{m}X$ be the subsheaf of
${\cal J}_{k}^{m}X$ consisting of elements of
the form:
\begin{eqnarray}
\omega(j^{k}f) = \sum_{j=1}^{k}
\sum_{|I_{1}|+j|I_{j}|= m} a_{I_{1}}(f)
(f^{'})^{I_{1}} (f^{(j)})^{I_{j}}.
\end{eqnarray} Note that the coefficient of
$(f^{'})^{I_{1}} (f^{(j)})^{I_{j}}$ depends
only on $I_{1}$ but is independent of $I_{j}$.
This sheaf shall be referred to as jet
differentials of linear type.

\bigskip
\noindent
{\bf Lemma 1.10}~~{\it The sheaf  ${\cal
L}_{k}^{m}X$ of jet differentials of linear type
is well-defined. For $m = k = 2$ we have ${\cal
L}_{2}^{2}X = {\cal
J}_{2}^{2}X$ and if $X$
is a Riemann surface then
${\cal
L}_{3}^{3}X = {\cal
J}_{3}^{3}X$.}

\medskip
\noindent
{\it Proof.}~~The change of variable formulas (4)
shows that a jet differential of the form (14)
is invariant by change of coordinates. ~~QED

\bigskip
There is a differentiatial
operation
$d : {\cal J}_{k}^{m}X \rightarrow {\cal
J}_{k+1}^{m+1}X$ naturally defined by:
\begin{eqnarray}
d \omega (j^{k+1}f) \stackrel{\rm def}{=}
(\omega(j^{k}f))^{'}.
\end{eqnarray} It should be noted that
in contrast to exterior differentiation of forms
$d
\circ d \ne 0$ on jet differentials. In
particular, given a holomorphic function
$\phi$ defined on some open neighborhood $U$,
\begin{eqnarray}
d^{(k)} \phi(j^{k}f) = (\phi \circ f)^{(k)}
\end{eqnarray} which is a non-trivial $k$-jet
differential for general $\phi$. Another
difference between jet differentials and
exterior differential forms is that a lower
order jet differential can be naturally
identified with a jet differential of higher
order. More precisely, the natural projection
$p_{kj} : J^{k}X \rightarrow J^{j}X$ defined
for $k \ge j$ induces an injection $p_{kj}^{*} :
{\cal J}_{j}^{m}X \rightarrow {\cal J}_{k}^{m}X$
defined naturally by:
\begin{eqnarray}
p_{kj}^{*}\omega(j^{k}f) \stackrel{\rm def}{=}
\omega(p_{kj}
(j^{k}f)) = \omega(j^{j}f).
\end{eqnarray}We shall simply write
$\omega(j^{k}f) = \omega(j^{j}f)$ if no
confusion arises. Moreover, the symmetric
product of a $k$-jet differential of weight $m$
and a $k'$-jet differential of weight $m'$ is a
$(k + k')$-jet differential of weight $m + m'$.

\bigskip
Consider first the case of $k = 2$ and denote by
$$p : J^{2}X \rightarrow p(J^{2}X) \subset
T^{2}X$$ the generically 2 to 1 map onto its
image as defined in the previous section. Let
$\omega \in H^{0}(U, T_2^{*}X)$ considered as a
linear (along the fibers) functional $$\omega :
T^{2}X|_U \rightarrow {\bf C}.$$ Consider the
composite map
$$\omega \circ p : J^{2}X|_U \rightarrow {\bf
C}.$$ By the definition of $p$ we observe that
$$\omega \circ p (\lambda.j^{2}f) = \omega
(\lambda^{2}. p(j^{2}f)) = \lambda^{2}
\omega\circ p (j^{2}f)$$ is homogeneous of
degree 2. In other words, the composite $\omega
\circ p$ is a section of ${\cal J}_{2}^{2}X$
over $U$. Thus we have a well-defined ${\bf
C}^{*}$-bundle map
$$q = p^{*} : T_{2}^{*} X \rightarrow {\cal
J}_{2}^{2}X.$$

\bigskip
Consider again the special case of a Riemann
surface $X$ and $k = 2$. Let $\omega \in H^{0}(U,
T_{2}^{*}X)$ be a 2-jet form then locally
$\omega$ is simply of the form
$$\omega = a d z + b d^2z^{(2)}$$
where $dz$ is the dual of $\partial / \partial
z$ and $d^2 z^{(2)}$ (the notation is formal and
should not be confused with differentiating
$z^2$ twice) with $a$ and
$b$ being local holomorphic functions. By the
definition of $p$ we have:
\begin{eqnarray}p^{*} \omega (j^2f) = a(f) f^{''}
+ b(f) (f^{'})^{2}.
\end{eqnarray} If $p^{*} \omega (j^2f) = a(f)
f^{''} + b(f) (f^{'})^{2} = 0$ for all $j^{2}f$
then, by taking $j^{2}f = (f(0), f^{'}(0),
f^{''}(0))$ such that $f^{'}(0) = 0$ and
$f^{''}(0) \ne 0$, we see that $a(f(0)) = 0$
(as $x = f(0)$ is an arbitrary point of $U$,
we have $a \equiv 0$) and so
$p^{*}
\omega(j^{2}f) = b(f) (f^{'})^{2}$ for all
$j^{2}f$. Now choosing $j^{2}f(0)$ so that
$f^{'}(0) \ne 0$ this time shows that $b(f(0)) =
0$. In other words, the map $p^{*} : T_2^{*}X
\rightarrow
{\cal J}_{2}^{2} X$ is injective.  On the
other hand, any homogeneous polynomial of degree
$2$ of $J_{x}^{2}X$ is  the germ of a section of
the form as in (18) where $a$ and $b$ are
holomorphic functions defined on some open
neighborhood of
$x$ and that $f(0) = x$. This shows that $p^{*}$
is actually an isomorphism between $T_2^{*}X$ and
${\cal J}_2^{2}X$. Analogously, a section of
$\odot^2 T^{*}X$ is of the form
$$\omega = b (dz \odot dz)$$ where $b$ is a
holomorphic function on $U$. Then the pull-back
$$p^{*}(\omega)(j^2f) = b(f) (f^{'})^{2}.$$ In
other words $$p^{*}(\odot^{2} T^{*}
X) \cong \odot^{2} T^{*} X$$ and hence we
conclude that the pull-back of the sequence:
$$0 \rightarrow \odot^{2} T^{*}X \rightarrow
T_2^{*}X \rightarrow T^{*}X \rightarrow
0$$ yields an exact sequence:
$$0 \rightarrow p^{*}(\odot^2 T^{*}X) \cong
\odot^2 T^{*}X \rightarrow {\cal J}_2^{2} X
\rightarrow {\cal J}_2^{2} X / p^{*}(\odot^2
T^{*}X) \cong T^{*} X
\rightarrow 0.$$ The same argument works also
for $k = 3$; for general
$k$ an analogous argument shows that $p_{k}^{*}
(T_{k}^{*}X)$ is isomorphic to the sheaf
of jet differentials of linear type ${\cal
L}_{k}^{k}$:

\bigskip
\noindent
{\bf Theorem 1.11}~~{\it For a complex manifold
$X$ the pull-back $p_{k}^{*} (T_{k}^{*}X)$ is
${\bf C}^{*}$-isomorphic to ${\cal
L}_{k}^{k}X$ where $p_{k} : J^{k}X \rightarrow
T^{k}X$ is the map defined by $(11)$. Then
$p_{k}^{*} (T_{k}^{*}X)$ is
${\bf C}^{*}$-isomorphic to ${\cal
L}_{k}^{k}X$.}

\medskip
\noindent
{\it Proof.}~~We have already seen the case $k =
2$ and suppose now that $\omega \in H^{0}(U,
T_{k}^{*}X)$ is a k-jet form then locally
$\omega$ is of the form
$$\omega = a_1 d z + a_2 d^2z^{(2)} + ... + a_k
d^kz^{(k)}$$ where $d^j z^{(j)}$ is the dual of
the differential operator $\partial^{j}/\partial
z^{j}$ and each $a_j$ is a holomorphic
function on $U$. Pulling back we get:
$$p_{k}^{*} \omega(j^{k}f) = \sum_{j=1}^{k}
a_{j} f^{(j)} (f^{'})^{k - j}$$
and suppose that
$p_{k}^{*} \omega(j^{k}f) \equiv 0$. Consider
first the case
$k = 3$
then for any
$x \in U$ choosing $j^{3}f$ so that $f(0)
= x, f^{'}(0) = 0$ shows that $a_{3} f^{(3)}(0) =
0$ so $a_{3}(x) = 0$. Since $x$ is arbitray the
function $a_3 \equiv 0$. Thus
$$0 \equiv p_{3}^{*} \omega(j^{3}f) = f^{'}
\{a_{2} f^{''} + a_{1}
(f^{'})^{2}\}$$ and so we have
$$a_{2} f^{''} +
a_{1} (f^{'})^{2} \equiv 0$$ on
$J^{2}X|_{U}
\setminus \{f^{'} \ne 0\}$. Let $\phi :
\Delta_{\epsilon} \rightarrow \Delta_{\epsilon}$
be a holomorphic function such that $\phi(0) = 0,
\phi^{'}(0) = 1$ then $(f \circ \phi)^{'} =
f^{'}(\phi) \phi^{'}, (f \circ \phi)^{''} =
f^{''}(\phi) (\phi^{'})^{2} + f^{'}(\phi)
\phi^{''}$ and the condition that $f^{'}(0) \ne
0$ implies that we may choose $\phi$ so that the
condition that the first jet is non-zero,
(i.e., $(f
\circ
\phi)^{'}(0)
\ne 0$) is preserved but
$$(f \circ \phi)^{''} =
f^{''}(\phi) (\phi^{'})^{2} + f^{'}(\phi)
\phi^{''} = 0$$ i.e., choose $\phi$ so that
$\phi^{''}(0) = - f^{''}(0)/f^{'}(0)$. This
yields:
$$a_{2}(\phi) (f \circ \phi)^{''}(0) =
a_{2}(\phi) (f
\circ \phi)^{''}(0) + a_{1}(\phi) ((f
\circ
\phi)^{'})^{2}(0) = 0$$so $a_{2} \equiv
0$ (because $x = f(0)$ is an arbitrary point) and
the original equation is reduced to the equation
$a_{1}(f) (f^{'})^{3} \equiv 0$. Thus by
choosing $f^{'}(0) \ne 0$ we conclude
that $a_{1}(f(0)) = 0$; this implies that $a_1
\equiv 0$ as well. This establishes
injectivity of the map $p_{3}^{*}$;
surjectivity follows from the fact that an
element of ${\cal J}_{3}^{3}X$ is of the form
$$a_{3} f^{(3)} + a_{2} f^{'} f^{''}
+ a_{1} (f^{'})^{3}.$$ In general we
have, by setting
$f^{'} = 0$, that $a_k \equiv 0$ and then:
$$0 \equiv p_{k}^{*} \omega = f^{'}
\sum_{j=1}^{k-1} a_{j} f^{(j)}
(f^{'})^{k-1-j}$$ and so $$0 \equiv
\sum_{j=1}^{k-1} a_{j} f^{(j)}
(f^{'})^{k-1-j}$$  on $J^{k}X \setminus \{f^{'}
\ne 0\}$. This shows injectivity;
surjectivity now follows from the definition
of ${\cal L}_{k}^{k}X$. The proof is then
completed by induction and by
reparametrization.~~QED

\bigskip
The following Theorem can be found (without
proof) in Green-Griffiths [G-G], we include a proof
here for the sake of completeness:

\bigskip
\noindent
{\bf Theorem 1.12}~~{\it There exists a
filtration of ${\cal J}_{k}^{m} X$:
$${ \cal J}_{k-1}^{m} X = {\cal F}_{k}^{0}
\subset {\cal F}_{k}^{1} \subset ... \subset
{\cal F}_{k}^{[m/k]} = {\cal J}_{k}^{m} X$$
$($where $[m/k]$ is the greatest integer smaller
than or equal to $m/k)$ such that
$${\cal F}_{k}^{i}/{\cal F}_{k}^{i-1} \cong
{\cal J}_{k-1}^{m-ki} X  \otimes (\odot^{i} T^{*}
X).$$}

\medskip
\noindent
{\it Proof.}~~
The filtrations are defined as follows. Since a
$(k-1)$-jet differential of weight $m$ is also a
$k$-jet differential of weight $m$ thus
$$F_k^0 = {\cal J}_{k-1}^mX \subset {\cal
J}_{k}^mX$$which in terms of the expression (13)
for jet differentials consists of elements of
which does not contain any terms involving
$f^{(k)}$; put it another way the exponent $I_k$
for
$f^{(k)}$ satisfies the condition $|I_k| = 0$:
\begin{eqnarray*}\omega(j^kf)  &=& \sum_{|I_1| +
2|I_2| + ... + (k-1)|I_{k-1}| = m} a_{I_1, ...,
I_{k-1}}(z) (f^{'})^{I_1} ...
(f^{(k-1)})^{I_{k-1}} \\&=& \sum_{|I_1| + 2|I_2|
+ ... + k|I_k| = m, |I_k| = 0} a_{I_1, ...,
I_k}(z) (f^{'})^{I_1} ...
(f^{(k)})^{I_k}.\end{eqnarray*}
For any
$0
\le j
\le [m/k]$ we define $F_{j} \subset {\cal
J}_k^mX$ to be the sheaf of germs consisting
elements so that $|I_k| \le j$:
\begin{eqnarray}F_k^j = \{\omega|\omega(j^kf) =
\sum_{|I_1| + 2|I_2| + ... + k|I_k| = m, |I_k|
\le j} a_{I_1, ..., I_k}(z) (f^{'})^{I_1} ...
(f^{(k)})^{I_k}\}.\end{eqnarray}By definition,
we have, for $1 \le j \le [m/k]$:
\begin{eqnarray*}F_k^j/F_k^{j-1} &=&
\{\omega|\omega(j^kf) \\&=&
\sum_{|I_1| + 2|I_2| + ... + k|I_k| = m, |I_k|
= j} a_{I_1, ..., I_k}(z) (f^{'})^{I_1} ...
(f^{(k)})^{I_k}\}\end{eqnarray*}and the claim is
that
$$F_k^j/F_k^{j-1} \cong {\cal
J}_{k-1}^{m-kj} X
\otimes (\odot^{j} T^{*} X).$$ We first establish
the special case of a Riemann surface. In this
case a
$k$-jet differential of weight $m$ is of the form
$$\omega(j^kf) = \sum_{i_1 + 2i_2 + ... +
ki_k = m} a_{i_1, ..., i_k}(z) ((z \circ
f)^{'})^{i_1} ... ((z \circ
f)^{(k)})^{i_k}$$where we identify $f$ with
$z$ being a local coordinate on an
open coordinate neighborhood $U \subset X$ and
$i_j$ are non-negative integers; the subsheaves
$F_k^j$ is of the form:
$$F_k^j =
\{\omega|\omega(j^kf) =
\sum_{i_1 + 2i_2 + ... +
ki_k = m, i_k
\le j} a_{i_1, ..., i_k}(z) (f^{'})^{i_1} ...
(f^{(k)})^{i_k}\}$$ for $0 \le j \le [m/k]$ and
$$F_k^j/F_k^{j-1} = \{\omega|\omega(j^kf) =
\sum_{i_1 + 2i_2 + ... +
ki_k = m, i_k
= j} a_{i_1, ..., i_k}(z) (f^{'})^{i_1} ...
(f^{(k)})^{i_k}\}$$ for $1 \le j \le [m/k]$.
We first define a map
$$L_U : F_k^j/F_k^{j-1}|_{U} \rightarrow
{\cal J}_{k-1}^{m-kj} X \otimes (\odot^{j} T^{*}
X)|_{U}$$where
\begin{eqnarray*}&&L_U(\sum_{i_1 + 2i_2 + ... +
ki_k = m, i_k
= j} a_{i_1, ..., i_k}(z) (f^{'})^{i_1} ...
(f^{(k)})^{i_k})\\
&&= (f^{(k)})^{j} \sum_{i_1 + 2i_2 + ... +
(k-1)i_{k-1} = m-kj} a_{i_1, ..., j}(z)
(f^{'})^{i_1} ...
(f^{(k-1)})^{i_{k-1}})\end{eqnarray*} The fact that $L_U$ is an
isomorphism is clear and the fact that $L = L_U$
(where ${\cal U} =
\{U\}$ is an open cover of $X$ by coordinate
neighborhoods) follows from the following
observation that (see (4)) if $(V, w)$ is another
coordinate neighborhood then
\begin{eqnarray*}
((w \circ f)^{(k)})^j = ((\partial
w/\partial z) (z \circ f)^{(k)} + P)^{j}
= ((\partial
w/\partial z) (z
\circ f)^{(k)})^j + Q
\end{eqnarray*}where $P$ and $Q$ are polynomials
in the variables $\partial^{s} w_i/\partial
z_l^{s}, 1
\le i, l \le n, 1 \le s \le k$ and in $(z \circ
f)^{(r)}, 1 \le r \le k - 1$. In particular, $Q$
is a $(k-1)$-jet differential of total weight $m
- kj$. In orther words,
$$((w \circ f)^{(k)})^j = (\partial
w/\partial z)^j ((z
\circ f)^{(k)})^j~~~~{\rm mod}~F_k^{j-1}$$
and the transisistion function $(\partial
w/\partial z)^j$ is the same as the transistioon
function for $\odot^j T^{*}X$.

\bigskip
The higher dimensional case is notationally more
complicated but the proof is essentially
the same.~~QED

\bigskip
As an immediate consequence (see Green-Griffiths
[G-G]) we have:

\bigskip
\noindent
{\bf Corollary 1.13}~~{\it Let $X$ be a
smooth projective variety then ${\cal J}_{k}^{m}
X$ admits a composition series whose factors
contain all bundles of the form:
$$(\odot^{i_{1}} T^{*}X) \otimes ... \otimes
(\odot^{i_{k}} T^{*}X)$$
where $i_{j}$ ranges over all non-negative
integers satisfying
$$i_{1} + 2i_{2} + ... + ki_{k} = m.$$
The first
Chern number of $c_{1}({\cal J}_{k}^{m} X)$
is given by:
$$c_{1}({\cal J}_{k}^{m} X) = \sum_{i_{1}
+ 2i_{2} + ... + ki_{k} = m, i_{j} \in {\bf
Z}_{\ge 0}} ~c_{1}((\odot^{i_{1}} T^{*}X)
\otimes ...
\otimes (\odot^{i_{k}} T^{*}X)).$$
In particular, for a curve $X = C$,
\begin{eqnarray*}c_{1}({\cal J}_{k}^{m} C)
&=&
\sum_{i_{1} + 2i_{2} + ... + ki_{k} = m, i_{j}
\in {\bf Z}_{\ge 0}} ~(i_{1} + i_{2} + ... +
i_{k}) c_{1}(T^{*}C).\end{eqnarray*}}

\bigskip
The preceding Theorem can be used in calculating
the Chern classes of ${\cal J}_k^m X$.

\bigskip
\noindent
{\bf Example 1.14}~~
For example, for $m = k = 2$, the filtration is
given by: $$\odot^{2} T^{*} X = {\cal J}_{1}^{2}
X = {\cal S}_{2}^{0} \subset {\cal S}_{2}^{1} =
{\cal J}_{2}^{2} X,~~{\cal
S}_{2}^{1}/{\cal S}_{2}^{0} \cong T^{*}X$$ we
have the following exact sequence:
$$0 \rightarrow \odot^{2}
T^{*}X
\rightarrow {\cal J}_{2}^{2} X \rightarrow
T^{*}X \rightarrow 0.$$ Thus
the first Chern numbers are related by the
formula:
$$c_{1}({\cal J}_{2}^{2} X) = c_{1}(\odot^{2}
T^{*}X) + c_{1}(
T^{*}X).$$

\bigskip
The filtration of ${\cal J}_{3}^{3} X$ is as
follows:
$${\cal J}_{3}^{3} X = {\cal S}_{3}^{1} \supset
S_{3}^{0} = {\cal J}_{2}^{3} X,~
{\cal J}_{3}^{3} X/{\cal J}_{2}^{3} X = {\cal
S}_{3}^{1} / {\cal S}_{3}^{0} \cong
T^{*}X.$$ Hence we have an exact sequence:
\begin{eqnarray*}0 \rightarrow {\cal J}_{2}^{3} X
\rightarrow {\cal J}_{3}^{3} X
\rightarrow T^{*}X \rightarrow 0.\end{eqnarray*}
Now the filtration of ${\cal J}_{2}^{3} X$ is
$${\cal J}_{2}^{3} X = {\cal S}_{2}^{1} \supset
{\cal S}_{2}^{0} = {\cal J}_{1}^{3} X,
~~{\cal J}_{2}^{3} X/{\cal J}_{1}^{3} X
\cong T^{*}X \otimes T^{*}X$$ and, since
${\cal J}_{1}^{3} X = \odot^{3}
T^{*} X$, we have an exact sequence:
\begin{eqnarray*}0 \rightarrow \odot^{3}
T^{*} X \rightarrow {\cal J}_{2}^{3} X
\rightarrow T^{*}X \otimes T^{*}X \rightarrow
0.\end{eqnarray*} From these 2 exact sequences we
get
\begin{eqnarray*}c_{1}({\cal J}_{3}^{3} X) =
c_{1}(T^{*}X) + c_{1}(T^{*}X \otimes T^{*}X) +
c_{1}(\odot^{3} T^{*}X).\end{eqnarray*} From basic
representation Theory (or just simple liner
algebra in this special case) we know that
$T^{*}X \otimes T^{*}X = \odot^{2} T^{*}X
\oplus \wedge^{2} T^{*}X$ hence,
\begin{eqnarray*}c_{1}({\cal J}_{3}^{3} X) =
c_{1}(T^{*}X) + c_{1}(\odot^{2} T^{*}X) +
c_{1}(\odot^{3} T^{*}X) +
c_{1}(\wedge^{2} T^{*}X).\end{eqnarray*} In
representation theory $\wedge^{2} T^{*}X$ is the
Weyl module $W_{1,1}^{*} X$ associate to the
partition $\{1, 1\}$ (see [F-H]). Thus we have:
\begin{eqnarray}c_{1}({\cal J}_{3}^{3} X) =
\sum_{j=1}^{3} c_{1}(\odot^{j} T^{*}X) +
c_{1}(W_{1,1}^{*}X).\end{eqnarray} In the
special case of a Riemann surface $\wedge^{2}
T^{*}X$ is the zero-sheaf. Thus for
a curve we have
$$c_{1}({\cal J}_{3}^{3} X) = (1 + 2 +
3) c_{1}( T^{*}X) = 6 c_{1}(
T^{*}X).$$

\bigskip
For $m = k = 4$, we have the following
filtrations:
$${\cal J}_{4}^{4} X = S_{4}^{1}
\supset S_{4}^{0} = {\cal J}_{3}^{4}
X,~~{\cal J}_{4}^{4} X/{\cal J}_{3}^{4} X =
S_{4}^{1} / S_{4}^{0}
\cong  T^{*}X,$$
$${\cal J}_{3}^{4} X = S_{3}^{1} \supset
S_{3}^{0} = {\cal J}_{2}^{4} X,~~{\cal
J}_{3}^{4} X/{\cal J}_{2}^{4} X = S_{3}^{1}
/ S_{3}^{0}
\cong T^{*}X \otimes T^{*}X,$$ and
$${\cal J}_{2}^{4} X = S_{2}^{2}
\supset S_{2}^{1}
\supset S_{2}^{0} = {\cal J}_{1}^{4}
X,$$ with $${\cal J}_{2}^{4} X/{\cal
S}_{2}^{1} =
\odot^{2} T^{*}X,~S_{2}^{1} / S_{2}^{0} \cong
T^{*}X
\otimes (\odot^{2} T^{*}X).$$
 The exact sequences associate
to the filtration for
${\cal J}_{4}^{4} X$ are:
\begin{eqnarray*}&&0 \rightarrow
{\cal J}_{3}^{4} X \rightarrow {\cal J}_{4}^{4}
X
\rightarrow T^{*}X \rightarrow 0; \cr
&&0 \rightarrow {\cal J}_{2}^{4} X
\rightarrow {\cal J}_{3}^{4} X
\rightarrow T^{*}X \otimes T^{*}X
\rightarrow 0; \cr &&0 \rightarrow S_{1}
\rightarrow {\cal J}_{2}^{4} X
\rightarrow \odot^{2} T^{*}X \rightarrow 0; \cr
&&0 \rightarrow \odot^{4} T^{*}X \rightarrow
S_{1}
\rightarrow T^{*}X \otimes (\odot^{2} T^{*}X)
\rightarrow 0.
\end{eqnarray*}
Thus the Chern number formula:
\begin{eqnarray*}c_{1}({\cal J}_{4}^{4} X)
&=& c_{1}(T^{*} X) + c_{1}(T^{*} X \otimes
T^{*} X) + c_{1}(\odot^{2} T^{*} X) \cr&&~~~~+
c_{1}(T^{*} X
\otimes
(\odot^{2} T^{*} X)) + c_{1}(\odot^{4} T^{*}
X).\end{eqnarray*} Note that (by elementary
representation theory)
$$T^{*} X
\otimes
(\odot^{k} T^{*} X) = W_{k,1}^{*} X \oplus
(\odot^{k+1} T^{*} X)$$ where
$W_{k, 1}^{*}$ is the Weyl module associate to
the partition $\{k, 1\}$ thus:
\begin{eqnarray}c_{1}({\cal J}_{4}^{4} X) =
c_{1}(\odot^{2} T^{*} X) + \sum_{i=1}^{4}
c_{1}(\odot^{i} T^{*} X) +
\sum_{i=1}^{2} c_{1}(W_{j,1}^{*}
X).\end{eqnarray}
In particulr, if
$X$ is a curve then
$$c_{1}({\cal J}_{4}^{4} X) = (1 + 2 + 2 +
3 + 4) c_{1}( T^{*}X) = 12 c_{1}(
T^{*}X).$$Recall that $c_1(T_4^{*}X) = 10
c_1(T^{*}X)$.

\bigskip
For $m = k = 5$, we have the following
filtrations:
$${\cal J}_{5}^{5} X = S_{5}^{1}
\supset S_{5}^{0} = {\cal J}_{4}^{5}
X,~{\cal J}_{5}^{5} X/{\cal J}_{4}^{5} X =
S_{5}^{1} / S_{5}^{0}
\cong  T^{*}X,$$
$${\cal J}_{4}^{5}
X = S_{4}^{1} \supset
S_{4}^{0} = {\cal J}_{3}^{5}
X,~{\cal J}_{4}^{5}
X/{\cal J}_{3}^{5}
X = S_{4}^{1} / S_{4}^{0}
\cong
T^{*}X
\otimes T^{*}X,$$
$${\cal J}_{3}^{5}
X = S_{3}^{1} \supset
S_{3}^{0} = {\cal J}_{2}^{5}
X,~~{\cal J}_{3}^{5}
X/{\cal J}_{2}^{5}
X \cong  T^{*}X
\otimes ({\cal J}_{2}^{2}
X),$$
$${\cal J}_{2}^{5} X = S_{2}^{2} \supset
S_{2}^{1}
\supset S_{2}^{0} = {\cal J}_{1}^{5} X,~~{\cal
J}_{2}^{5} X/{\cal S}_{2}^{1} = (\odot^{2}
T^{*}X) \otimes T^{*}X,$$
$$S_{2}^{1} /
S_{2}^{0}
\cong T^{*}X
\otimes (\odot^{3} T^{*}X).$$
 The exact sequences associate
to the filtration for
${\cal J}_{5}^{5} X$ are:
\begin{eqnarray*}&&0 \rightarrow {\cal
J}_{4}^{5} X \rightarrow {\cal J}_{5}^{5} X
\rightarrow T^{*}X \rightarrow 0; \cr
&&0 \rightarrow {\cal J}_{3}^{5} X \rightarrow
{\cal J}_{4}^{5} X
\rightarrow T^{*}X \otimes T^{*}X
\rightarrow 0; \cr &&0 \rightarrow {\cal
J}_{2}^{5} X
\rightarrow {\cal J}_{3}^{5} X
\rightarrow T^{*}X
\otimes {\cal J}_{2}^{2} X \rightarrow 0; \cr
&&0 \rightarrow {\cal S}_{2}^{1} \rightarrow
{\cal J}_{2}^{5} X
\rightarrow (\odot^{2}
T^{*}X) \otimes T^{*}X
\rightarrow 0,\cr
&&0 \rightarrow \odot^{5}
T^{*}X \rightarrow
{\cal S}_{2}^{1}
\rightarrow T^{*}X \otimes (\odot^{3}
T^{*}X)
\rightarrow 0.
\end{eqnarray*}
This yields the formula:
\begin{eqnarray*}
&&c_{1}({\cal J}_{5}^{5} X) \\&&= c_{1}(T^{*}
X) + c_{1}(T^{*} X \otimes T^{*} X) +
c_{1}(T^{*} X
\otimes {\cal J}_{2}^{2} X) \\&&~~+
c_{1}((\odot^{2} T^{*} X) \otimes T^{*} X) +
c_{1}(T^{*} X
\otimes (\odot^{3} T^{*} X)) + c_{1}(\odot^{5}
T^{*} X)\\ &&= c_{1}(T^{*} X) +
c_{1}(T^{*} X \otimes T^{*} X) + c_{1}(T^{*} X
\otimes (\odot^{2} T^{*}X)) \\&&~~+ c_{1}(T^{*} X
\otimes T^{*}X) + c_{1}((\odot^{2} T^{*}
X) \otimes T^{*} X) \\&&~~+ c_{1}(T^{*} X
\otimes (\odot^{3} T^{*} X)) + c_{1}(\odot^{5}
T^{*} X)
\end{eqnarray*}
where we have used the fact that
$$c_{1}(T^{*} X
\otimes {\cal J}_{2}^{2} X) = c_{1}(T^{*} X
\otimes (\odot^{2} T^{*}X)) + c_{1}(T^{*} X
\otimes T^{*}X).$$ Recall that
$$T^{*} X \otimes T^{*} X = \odot^{2} T^{*}X
\oplus \wedge^{2} T^{*}X,$$
$$T^{*} X
\otimes (\odot^{2} T^{*}X) = \odot^{3} T^{*}X
\oplus W_{2, 1}^{*} X,$$
$$T^{*} X
\otimes (\odot^{3} T^{*}X) = \odot^{4} T^{*}X
\oplus W_{3, 1}^{*} X$$
(in general we have
$$T^{*} X
\otimes (\odot^{d} T^{*}X) = \odot^{d+1} T^{*}X
\oplus W_{d, 1}^{*} X.)$$
Thus we have:
\begin{eqnarray*}
&&c_{1}({\cal J}_{5}^{5} X) \\&&=
\sum_{j=2}^{3} c_{1}(\odot^{j} T^{*} X) +
\sum_{j=1}^{5} c_{1}(\odot^{j} T^{*} X) +
\sum_{j=1}^{2} c_{1}(W_{j,1}^{*} X) +
\sum_{j=1}^{3} c_{1}(W_{j,1}^{*} X).
\end{eqnarray*}
In particulr, if
$X$ is a curve then
$$c_{1}({\cal J}_{5}^{5} X) = (1 + 2 + 3 + 2 +
3 + 4 + 5) c_{1}( T^{*}X) = 20 c_{1}(
T^{*}X).$$

\bigskip
For $m = k = 6$, we have the following
filtrations:
$${\cal J}_{6}^{6} X = S_{6}^{1} \supset
S_{6}^{0} = {\cal J}_{5}^{6} X,~~{\cal
J}_{6}^{6} X/{\cal J}_{5}^{6} X = S_{6}^{1}
/ S_{6}^{0}
\cong  T^{*}X,$$
$${\cal J}_{5}^{6} X = S_{5}^{1} \supset
S_{5}^{0} = {\cal J}_{4}^{6} X,~~{\cal
J}_{5}^{6} X/{\cal J}_{4}^{6} X = S_{5}^{1}
/ S_{5}^{0}
\cong
T^{*}X
\otimes T^{*}X,$$
$${\cal J}_{4}^{6} X = S_{4}^{1} \supset
S_{4}^{0} = {\cal J}_{3}^{6} X,~~{\cal
J}_{4}^{6} X/{\cal J}_{3}^{6} X \cong  T^{*}X
\otimes {\cal J}_{2}^{3} X,$$
$${\cal J}_{3}^{6} X = S_{3}^{2} \supset
S_{3}^{1}
\supset S_{3}^{0} = {\cal J}_{2}^{6}
X,$$ with factors $${\cal J}_{3}^{6} X/{\cal
S}_{3}^{1} =
\odot^{2} T^{*}X,~~S_{3}^{1} /
S_{3}^{0}
\cong T^{*}X
\otimes {\cal J}_{2}^{3} X,$$
$${\cal J}_{2}^{6} X = S_{2}^{3} \supset
S_{2}^{2}
\supset S_{2}^{1} \supset S_{2}^{0} = \odot^{6}
T^{*} X,~~{\cal J}_{2}^{6} X/{\cal
S}_{2}^{2} = \odot^{3} T^{*}X,$$
$$S_{2}^{2} /
S_{2}^{1}
\cong (\odot^{2} T^{*}X)
\otimes (\odot^{2} T^{*} X), ~~S_{2}^{1} /
S_{2}^{0}
\cong T^{*}X
\otimes (\odot^{4} T^{*} X).$$
 The exact sequences associate
to the filtration for
${\cal J}_{6}^{6} X$ are:
\begin{eqnarray*}&&0 \rightarrow {\cal
J}_{5}^{6} X \rightarrow {\cal J}_{6}^{6} X
\rightarrow T^{*}X \rightarrow 0; \cr
&&0 \rightarrow {\cal J}_{4}^{6} X \rightarrow
{\cal J}_{5}^{6} X
\rightarrow T^{*}X \otimes T^{*}X
\rightarrow 0; \cr &&0 \rightarrow {\cal
J}_{3}^{6} X
\rightarrow {\cal J}_{4}^{6} X
\rightarrow T^{*}X
\otimes {\cal J}_{2}^{3} X \rightarrow 0; \cr
&&0 \rightarrow {\cal S}_{3}^{1} \rightarrow
{\cal J}_{3}^{6} X
\rightarrow \odot^{2}
T^{*}X
\rightarrow 0,\cr
&&0 \rightarrow {\cal J}_{2}^{6} X \rightarrow
{\cal S}_{3}^{1}
\rightarrow T^{*}X \otimes {\cal J}_{2}^{3} X
\rightarrow 0,\cr
&&0 \rightarrow {\cal S}_{2}^{2} \rightarrow
{\cal J}_{2}^{6} X
\rightarrow \odot^{3}
T^{*}X
\rightarrow 0,\cr
&&0 \rightarrow {\cal S}_{2}^{1} \rightarrow
{\cal S}_{2}^{2}
\rightarrow (\odot^{2}
T^{*}X) \otimes (\odot^{2}
T^{*}X)
\rightarrow 0,\cr
&&0 \rightarrow \odot^{6}
T^{*}X \rightarrow
{\cal S}_{2}^{1}
\rightarrow
T^{*}X \otimes (\odot^{4}
T^{*}X)
\rightarrow 0.
\end{eqnarray*}
This yields the formula:
\begin{eqnarray*}
c_{1}({\cal J}_{6}^{6} X) &=& c_{1}(T^{*} X) +
c_{1}(T^{*} X \otimes T^{*} X) + 2c_{1}(T^{*} X
\otimes {\cal J}_{2}^{3} X)\\&& +
c_{1}(\odot^{2} T^{*} X) + c_{1}(\odot^{3} T^{*}
X)
+ c_{1}((\odot^{2} T^{*} X)
\otimes (\odot^{2} T^{*} X)) \\ &&+
c_{1}(T^{*}X
\otimes (\odot^{4} T^{*} X)) +
c_{1}(\odot^{6} T^{*} X)\\&=& c_{1}(T^{*} X) +
3c_{1}(T^{*} X \otimes T^{*} X) + 2c_{1}(T^{*} X
\otimes (\odot^{3} T^{*}X))\\&& +
c_{1}(\odot^{2} T^{*} X) + c_{1}(\odot^{3} T^{*}
X)
+ c_{1}((\odot^{2} T^{*} X)
\otimes (\odot^{2} T^{*} X)) \\ &&+
c_{1}(T^{*}X
\otimes (\odot^{4} T^{*} X)) +
c_{1}(\odot^{6} T^{*} X).
\end{eqnarray*}
where we have used the fact that
$$c_{1}(T^{*} X
\otimes {\cal J}_{2}^{3} X)) = c_{1}(T^{*} X
\otimes (\odot^{3} T^{*}X)) + c_{1}(T^{*} X
\otimes T^{*}X).$$ In particulr, if $X$ is
a curve then
$$c_{1}({\cal J}_{6}^{6} X) = (1 + 6 + 8 + 2
+ 3 + 4 + 5 + 6) c_{1}( T^{*}X) = 35 c_{1}(
T^{*}X).$$

\bigskip
The calculation before can be carried
out in a much simpler fashion as follows. A
{\sl partition} of a natural number
$m$ is a set of {\it positve integers} $k_1, ...,
k_q$ such that $m = k_1 + ... + k_q$.
A partition can be expressed as $$m =
\sum_{j=1}^{k} j i_{j}$$ where the
integers $i_j = \#$ of $j$'s in $\{k_1, ...,
k_q\}$ are non-negative. Obviously we have $1
\le q \le k$ and $1 \le  k_i \le k$ for all
$i$. The following result is well-known in
representation theory and in combinatorics
(see [H-W]):

\bigskip
\noindent
{\bf Theorem 1.15}~~{\it The number, denoted
$p(m)$, of classes of
$S_m~($the symmetric group on $m$
elements$)$ is equal to the number of partitions of
$m$ and also to the number of $($inequivalent$)$
irreducible representations of $S_m$. The number
$p(m)$ is asymptotically approximated by the
formula of Hardy-Ramanujan
$$p(m) \sim {e^{\pi \sqrt{2m/3}} \over 4m
\sqrt{3}}.$$}

\bigskip
\noindent
{\bf Remark 1.16}~~The first few numbers are as
follows:
$$p(1) = 1, p(2) = 2, p(3) = 3, p(4) = 5, p(5) =
7, p(6) = 11, p(7) = 15,$$
$$ p(8) = 22, p(9) = 30, p(10) = 42, p(11) = 56,
p(12) = 77, p(13) = 101.$$

\bigskip
We are interested in the case
$$k = \lambda_1 +
\lambda_2 + ... + \lambda_{\rho_{\lambda}}$$
where
$\lambda_1 \ge \lambda_2 ... \ge
\lambda_{\rho_{\lambda}} \le 1$. Define $l_i =
\lambda_i +
\rho_{\lambda} - i, i = 1, ...,
\rho_{\lambda}$. Then the dimension
$d_{\lambda}$ of the representation
$V_{\lambda},
\lambda = (\lambda_1, ...,
\lambda_{\rho_{\lambda}})$ associated to the
partition
$\lambda$ is given by the formula $d_{\lambda}
= 1$ if $\rho_{\lambda} = 1$ and for
$\rho_{\lambda}
\ne 1$ (see [F-H], p. 50):
\begin{eqnarray}
d_{\lambda} = {k! \over l_1! ...
l_{\rho_{\lambda}}!}
\prod_{1 \le i < j \le \rho_{\lambda}} (l_i -
l_j)
\end{eqnarray}
The number $\rho_{\lambda}$ shall be referred
to as {\it the length of the partition
$\lambda$}.

\bigskip
We consider also the case of partitioning a
number by a partition of fixed length $k$. Denote
by
$p_{k}(m)$ the number of solutions of the
equation
\begin{eqnarray*}x_1 + ... + x_k = m\end{eqnarray*}
with the condition that
$1 \le x_k \le x_{k-1} \le ... \le x_1$. This
number is obviously equal to the number of
solutions of the equation
\begin{eqnarray*}y_1 + ... + y_k = m -
k\end{eqnarray*} with the condition that
$0 \le y_k \le y_{k-1} \le ... \le y_1$. If
there are exactly $i$ of the integers $y_i$
which are positive then these are the solutions
of $x_1 + ... + x_i = m - k$ ($x_i
\leftrightarrow y_i+1$) and so there are
$p_{i}(m-k)$ of such solutions; consequently we
have:

\bigskip
\noindent
{\bf Theorem 1.17}~~{\it With the notations above
we have
$$p_{k}(m) = \sum_{i=0}^{k} p_{i}(m-k)$$ if $1 \le
k
\le m$ and with the convention that $p_0(0)
= 1, p_0(m) = 0$ if $m > 0$ and $p_{k}(m) = 0$ if
$k > m$.}

\bigskip
The following identity is easily established by
induction:

\bigskip
\noindent
{\bf Theorem 1.18}~~{\it The number
$p_{k}(m)$ satisfies the following recursive
relation: $p_k(m) = p_{k-1}(m-1) +
p_k(m-k)$.}

\bigskip
Obviously we have $p_1(m) = p_m(m) = 1$ and
$p_2(m) = m/2$ or $(m - 1)/2$ according to $m$
being even or odd. Thus Theorem 1.19 yields
$p_3(m) = p_2(m-1) + p_3(m - 3)$, $p_4(m) =
p_3(m-1) + p_4(m - 4),$
$p_5(m) = p_4(m-1) + p_5(m - 5)$
and we get for example
$$p_1(6) = 1, p_2(6) = 3, p_6(6) = 1$$
$$p_3(6) = p_2(5) + p_3(3) = 3,$$
$$p_4(6) = p_3(5) = p_2(4) = 2,$$
$$p_5(6) = p_4(5) = p_3(4) = p_2(3) = 1$$hence as
$p(m) = \sum_k p_k(m)$ we have
$$p(6) = \sum_{k=1}^{6} p_{k}(6) = 1 + 3 + 3 + 2
+ 1 + 1= 11.$$
For $m = 7$ we have
$$p_1(7) = 1, p_2(7) = 3, p_7(7) = 1$$
$$p_3(7) = p_2(6) + p_3(4) = p_2(6) + p_2(3)
= 4,$$
$$p_4(7) = p_3(6) = 3,$$
$$p_5(7) = p_4(6) = 2,$$
$$p_6(7) = p_5(6) = 1$$
$$p(7) = \sum_{k=1}^{7} p_{k}(7) = 1 + 3 + 4 + 3
+ 2 + 1 + 1 = 15.$$

\bigskip
The {\it total length of all partitions} $L(m)$
of a positive integer
$m$ is defined to be
\begin{eqnarray*}
L(m) = \sum_{j=1}^{m} j p_j(m).
\end{eqnarray*}
For example if $m = 6$ then $L(6) = 1 + 6 + 9 +
8 + 5 + 6 = 35$ and for $m = 7, L(7) = 1 + 6 + 12 +
12 + 10 + 6 + 7 = 54.$ More generally for $k \le m$
\begin{eqnarray}L_k(m) = \sum_{j=1}^{k} j
p_j(m)\end{eqnarray}shall be referred to as the total length
of partitions of $m$ of length at most $k$.
The following Lemma is easily established from the
definitions:

\bigskip
\noindent
{\bf Lemma 1.19}~~{\it With the notations
above we have
$$\sum_{\lambda} \sum_{j=1}^{k} i_{j} =
\sum_{\lambda}
\rho_{\lambda} = \sum_{j=1}^{k} j p_{j}(m)$$
where the sum on the right is taken over all
partition $\lambda = (\lambda_{1}, ...,
\lambda_{\rho_{\lambda}})$ of
$m, 1 \le \lambda_{\rho_{\lambda}} \le ... \le
\lambda_{2} \le
\lambda_{1}, \rho_{\lambda} \le k$ and
$i_{j} = \#$ of $j$'s in $\{\lambda_{1}, ...,
\lambda_{\rho_{\lambda}}\}$.}

\bigskip
One has the following well-known
asymptotic formula:

\bigskip
\noindent
{\bf Theorem 1.20}~~{\it For $m \rightarrow
\infty$ the number $p_k(m)$ is asymptotically
given by:
$$p_k(m) \sim {m^{k-1} \over (k - 1)! k!}.$$}

\bigskip
The preceding discussions yield the following
Theorem:

\bigskip
\noindent
{\bf Theorem 1.21}~~{\it Let $X$ be a
non-singular pojective curve then the Chern
number of ${\cal J}_{k}^{m}X$ is given by
$$c_{1}({\cal J}_{k}^{m}X) = L_k(m)
c_{1}({\cal K}_{X}) = \sum_{j=1}^{k} j p_{j}(m)
c_{1}({\cal K}_{X}) =
\sum_{j=1}^{k} j p_{j}(m) c_{1}({\cal K}_{X})$$
where ${\cal K}_{X}$ is the canonical bundle of $X$. If we
fix $k$ and let $m \rightarrow \infty$ then asymptotically:
$$c_{1}({\cal J}_{k}^{m}X) \sim k p_k(m) \sim {m^{k-1}
\over (k - 1)! (k - 1)!}.$$}

\bigskip
We give as examples the explicit
calculation of the above.  For $m = k = 3$, we have
$p(3) = 3$ and the possible indices are tabulated
below:

\bigskip
\begin{center}
\begin{tabular}{|l|l|l|l|l|l|l|r|} \hline
 &$\lambda$&$\rho_{\lambda}$ &$d_{\lambda}$&
$i_1$ &
$i_2$ &
$i_3$ &$\sum_{j=1}^{k} i_j$\\ \hline
1 & (1, 1, 1)&3&1& 3 & 0 & 0&3 \\
\hline
2 &(2, 1)&2&2& 1 & 1 &0&2 \\ \hline
3 &(3) &1 &1& 0 & 0 &1&1\\ \hline
\end{tabular}
\end{center}

\bigskip
\noindent
The Chern number of a curve $X$ is obtained by
summing the last column:
$$c_{1}({\cal J}_{3}^{3}X) = (1 + 2 + 3)
c_{1}(T^{*}X) = 6c_{1}(T^{*}X).$$

\bigskip
For $m = k = 4$, we have $p(4) = 5$ and the
possible indices are listed below

\bigskip
\begin{center}
\begin{tabular}{|l|l|l|l|l|l|l|l|r|} \hline
 &$\lambda$&$\rho_{\lambda}$&$d_{\lambda}$&
$i_1$ & $i_2$ &
$i_3$ &$i_4$ &$\sum_{j=1}^{k} i_j$\\ \hline
1 &(1, 1, 1, 1)&4&1& 4 & 0 & 0 &0 &4 \\
\hline
2 &(2, 1, 1)&3&3& 2 & 1 &0 &0 &3 \\ \hline
3 &(3, 1)&2&3& 1 & 0 &1 &0 &2\\ \hline
4 &(2, 2)&2&2& 0 & 2 &0 &0 &2\\ \hline
5 &(4)&1&1& 0 & 0 &0 &1 &1\\ \hline
\end{tabular}
\end{center}

\bigskip
\noindent
The Chern numberof a curve $X$ is obtained by
summing the last column:
$$c_{1}({\cal J}_{4}^{4}X) = 12c_{1}(T^{*}X).$$

\bigskip
For $m = k = 5$, we have $p(5) = 7$ and the
possible indices are listed below

\bigskip
\begin{center}
\begin{tabular}{|l|l|l|l|l|l|l|l|l|r|} \hline
 &$\lambda$&$\rho_{\lambda}$&$d_{\lambda}$&
$i_1$ &
$i_2$ &
$i_3$ &$i_4$ & $i_5$ &$\sum_{j=1}^{k} i_j$\\
\hline 1 &(1, 1, 1, 1, 1)&5&1& 5 & 0 & 0 &0 & 0&5
\\
\hline
2 &(2, 1, 1, 1)&4&4& 3 & 1 &0 &0 & 0&4 \\ \hline
3 &(3, 1, 1)&3&6& 2 & 0 &1 &0 & 0&3\\ \hline
4 &(2, 2, 1)&3&5& 1 & 2 &0 &0 & 0&3\\ \hline
5 &(4, 1)&2&4& 1 & 0 &0 &1 & 0&2\\ \hline
6 &(3, 2)&2&15& 0 & 1 &1 &0 & 0&2\\ \hline
7 &(5)&1&1& 0 & 0 &0 &0 & 1&1\\ \hline
\end{tabular}
\end{center}

\bigskip
\noindent
The Chern numberof a curve $X$ is obtained by
summing the last column:
$$c_{1}({\cal J}_{5}^{5}X) = 20c_{1}(T^{*}X).$$

\bigskip
For $m = k = 6$, we have $p(6) = 11$ and the
possible indices are listed below

\bigskip
\begin{center}
\begin{tabular}{|l|l|l|l|l|l|l|l|l|l|r|} \hline
 &$\lambda$&$\rho_{\lambda}$&$d_{\lambda}$&
$i_1$ &
$i_2$ &
$i_3$ &$i_4$ & $i_5$& $i_6$  &$\sum_{j=1}^{k}
i_j$\\
\hline 1 &(1, 1, 1, 1, 1, 1)&6&1& 6 & 0 & 0 &0 &
0 &0 &6
\\
\hline
2 &(2, 1, 1, 1, 1)&5&5& 4 & 1 &0 &0 & 0&0&5 \\
\hline 3 &(3, 1, 1, 1)&4&10& 3 & 0 &1 &0 &
0&0&4\\ \hline 4 &(2, 2, 1, 1)&4&9& 2 & 2 &0
&0 & 0&0&4\\ \hline 5 &(4, 1, 1)&3&10& 2 &
0 &0 &1 & 0&0&3\\
\hline 6 &(3, 2, 1)&3&36& 1 & 1 &1 &0 &
0&0&3\\ \hline 7 &(2, 2, 2)&3&5& 0 &
3 &0 &0 & 0&0&3\\
\hline
8 &(5, 1)&2&30& 1 & 0 &0
&0 & 1&0&2\\ \hline
9 &(4, 2)&2&9& 0 & 1 &0 &1 &
0&0&2\\ \hline 10 &(3, 3)&2&5& 0 & 0 &2
&0 & 0&0&2\\ \hline 11 &(6)&1&1& 0 &
0 &0 &0 & 0&1&1\\
\hline
\end{tabular}
\end{center}

\bigskip
\noindent
The Chern numberof a curve $X$ is obtained by
summing the last column:
$$c_{1}({\cal J}_{6}^{6}X) = 35c_{1}(T^{*}X).$$

\bigskip
For $m = k = 7$, we have $p(7) = 15$ and the
possible indices are listed below

\bigskip
\begin{center}
\begin{tabular}{|l|l|l|l|l|l|l|l|l|l|l|r|}
\hline
 &$\lambda$&$\rho_{\lambda}$&$d_{\lambda}$&
$i_1$ &
$i_2$ &
$i_3$ &$i_4$ & $i_5$& $i_6$ & $i_7$
&$\sum_{j=1}^{k} i_j$\\
\hline 1 &(1, 1, 1, 1, 1, 1, 1)&7&1& 7 & 0 & 0
&0 & 0 &0 &0&7
\\
\hline
2 &(2, 1, 1, 1, 1, 1)&6&6& 5 & 1 &0 &0 & 0&0&0&6
\\ \hline 3 &(3, 1, 1, 1, 1)&5&15& 4 & 0 &1 &0
& 0&0&0&5\\ \hline 4 &(2, 2, 1, 1, 1)&5&14& 3
& 2 &0 &0 & 0&0&0&5\\
\hline 5 &(4, 1, 1, 1)&4&20& 3 &
0 &0 &1 & 0&0&0&4\\
\hline 6&(3, 2, 1, 1)&4&35 & 2 & 1 &1 &0 &
0&0&0&4\\ \hline 7&(2, 2, 2, 1)&4&14 & 1 & 3 &0
&0 & 0&0&0&4\\ \hline
8&(5, 1, 1)&3 &15& 2 &
0 &0 &0 & 1&0&0&3\\
\hline 9 &(4, 2, 1)&3&35& 1 &
1 &0 &1 & 0&0&0&3\\
\hline 10&(3, 3, 1)&3&21 & 1 & 0 &2 &0 &
0&0&0&3\\ \hline 11&(3, 2, 2)&3&21 & 0 & 2 &1 &0
& 0&0&0&3\\ \hline
12&(6, 1)&2&6 & 1 &
0 &0 &0 & 0&1&0&2\\
\hline 13&(5, 2)&2&14 & 0 &
1 &0 &0 & 1&0&0&2\\
\hline 14&(4, 3) &2&14& 0 & 0 &1 &1 & 0&0&0&2\\
\hline 15 &(7)&1&1& 0 & 0 &0 &0 & 0&0&1&1\\
\hline
\end{tabular}
\end{center}

\bigskip
\noindent
The Chern numberof a curve $X$ is obtained by
summing the last column:
$$c_{1}({\cal J}_{7}^{7}X) = 54c_{1}(T^{*}X).$$
We list below the next few values of $L(k)$:
$$L(8) = 86, L(9) = 128, L(10) = 192, L(11) = 275,
L(12) = 399, L(13) = 556$$ $$L(14) = 780,
L(15) = 1068, L(16) = 1463.$$

\bigskip
\bigskip
\noindent
{\bf \S~2~Computation of Chern Classes in Complex
Surfaces}

\bigskip
 We now treat the case of a complex
surface (i.e., complex dimension 2). First we
establish some basic facts:

\bigskip
\noindent
{\bf Lemma 2.1}~~{\it Let $X$ be a
nonsingular complex surface then
\begin{eqnarray*}
&&c_1(\odot^m T^{*}X) = {m(m+1) \over 2}
c_1(T^{*}X),\\ &&c_2(\odot^m T^{*}X) = a(m)  c_1^2(T^{*}X) +
b(m) c_2(T^{*}X)
\end{eqnarray*}where $a(m) = m(m^2 - 1)(3m + 2) / 24,
b(m) = m(m + 1)(m + 2) / 6$.}

\medskip
\noindent
{\it Proof.}~~the case of the first Chern class is straight
forward and the calculation is omitted (see section 4 for a
slightly more general calculation.  To
compute the Chern numbers of
$\odot^{2} E$ we proceed formally by writing
the total Chern class $c(E) = (1 +
(\lambda_{1} + \lambda_{2}) x + \lambda_{1}
\lambda_{2} x^{2})$ then the total Chern class
of  $\odot^{2} E$ is (keep in mind
that rank
$\odot^{2} E = 3$):
$$(1 + 2 \lambda_{1} x) (1 + 2 \lambda_{2} x) (1
+ (\lambda_{1} +
\lambda_{2}) x)$$
and a calculation (mod $x^{3}$) yields:
$$1 + 3(\lambda_{1} +
\lambda_{2}) x + [4 \lambda_{1} \lambda_{2}
+ 2(\lambda_{1} +
\lambda_{2})^{2}] x^{2}.$$
This shows that
\begin{eqnarray*}c_{1}(\odot^{2} E)
= 3 c_{1}(E),~~
c_{2}(\odot^{2} E) = 2
c_{1}^{2}(E) + 4 c_{2}(E)
\end{eqnarray*}and the Lemma is verified in this case.

\bigskip
Next we compute the Chern numbers of
$\odot^{3}~E$. With a similar
formalism (and keep in mind that the rank of
$\odot^{3}~E$ is 4), we have:
\begin{eqnarray*}c(\odot^{3}E)
&=& (1 + 3 \lambda_{1} x) (1 + 3\lambda_{2} x) (1 +
(2\lambda_{1} + \lambda_{2}) x) (1 +
(\lambda_{1} + 2\lambda_{2}) x) \cr &=& 1 +
6(\lambda_{1} + \lambda_{2}) x +
\{11 (\lambda_{1} +
\lambda_{2})^{2} + 10
\lambda_{1} \lambda_{2}\} x^{2}~~({\rm
mod}~x^3).\end{eqnarray*} This shows that
$$
c_{1}(\odot^{3}E) = 6 c_{1}(E),~~
c_{2}(\odot^{3}E) = 11
c_{1}^{2}(E) + 10 c_{2}(E).
$$

\bigskip
For general $m$ we observe that
\begin{eqnarray*}
c_{1}(\odot^{m}E) = \left\{\begin{array}{ll}
(mp_2(m) + {m \over 2})c_1(E), &\mbox{if $m$ is
even,}\\ (mp_2(m) + m)c_1(E), &\mbox{if $m$ is odd}
\end{array}
\right.
\end{eqnarray*}where $p_2(m)$ is the number of
solutions of $m$ with partitions of fixed length 2
defined above. The Lemma follows by recalling that
$p_2(m) = m/2$ (resp. $(m -1)/2$) if $m$ is even
(resp. odd). For the tensor product we observe that
$$
c_{1}(\otimes^{m}E) =
\sum_{i=0} {i + (m-i) \over 2} C_{i}^{m} c_1(E)
$$ which follows from the fact that a partition
$m$ of length 2 can be written simply as
$l = (l_1 = i, l_2 = m - i)$.
Previously, this was defined by requiring that
$l_1 \ge
l_2
\ge 1$ but in the preceeding formula we include
all parititions $l = (l_1,
l_2), l_i \ge 0, l_1 + l_2
= m$. This accounts for the extra term $m/2$
(resp. $m$) in the formula for the symmetric
product and also the factor $1/2$ in the formula
for tensor product.~~

\bigskip
The second Chern class is somewhat more
complicated. Given an integer $m$ the {\it
non-negative} partitions of $m$ of length 2 are
$\{(m - i, i), i = 0, ...,
m\}$ and
$$c(\odot^m E) = \prod_{i=0}^{m} (1 + ((m-i)
\lambda_1 + i \lambda_2)x)~~({\rm
mod}~x^{3}).$$The coefficients of $x^2$ is the
second Chern class and is given by the following
sums if $m$ is even:
$$s_0 = \sum_{i=0}^{{m \over 2} - 1} ((m - i)
\lambda_1 + i \lambda_2)(i \lambda_1 + (m - i)
\lambda_2)$$
$$s_1 = \sum_{i=1}^{m - 1} \{m
\lambda_1 ((m - i) \lambda_1 + i
\lambda_2) + m
\lambda_2 (i \lambda_1 + (m - i)
\lambda_2)\}$$
\begin{eqnarray*}s_2 &=& \sum_{i=2}^{m - 2} \{((m -
1)
\lambda_1 + \lambda_2) ((m - i) \lambda_1 + i
\lambda_2) +\\&& + (\lambda_1 + (m - 1)
\lambda_2) (i \lambda_1 + (m - i)
\lambda_2)\}\end{eqnarray*}
$$...$$ $$...$$
$$...$$
\begin{eqnarray*}s_j &=& \sum_{i=j}^{m - j} \{((m -
j)
\lambda_1 + j \lambda_2) ((m - i) \lambda_1 + i
\lambda_2) + \\&&+ (j\lambda_1 + (m - j)
\lambda_2) (i \lambda_1 + (m - i)
\lambda_2)\}\end{eqnarray*}
$$...$$ $$...$$
$$...$$
$$s_{(m/2)-1} = \{({m \over 2} + 1) \lambda_1 + ({m
\over 2} - 1) \lambda_2)\}\{{m \over 2}
\lambda_1 + {m
\over 2} \lambda_2)\}$$
$$c_2(\odot^m E) = s_0 + s_1 + ... +
s_{(m/2)-1}.$$By simple algebra, we have
\begin{eqnarray*}s_0 &=&
\sum_{i=0}^{{m \over 2} - 1} i(m
- i) (\lambda_1^2 + \lambda_2^2) + \sum_{i=0}^{{m
\over 2} - 1} (i^2 + (m - i)^2) \lambda_1
\lambda_2\\&=&
\sum_{i=0}^{{m \over 2} - 1} i(m
- i) (\lambda_1 + \lambda_2)^2 + \sum_{i=0}^{{m
\over 2} - 1} (m - 2i)^2 \lambda_1
\lambda_2\\&=&
\sum_{i=0}^{{m \over 2} - 1} i(m
- i) c_1(E)^2 + \sum_{i=0}^{{m
\over 2} - 1} (m - 2i)^2 c_2(E).
\end{eqnarray*}For $s_j, j \ge 1$ the main
observation is that each of these can
be expressed as $(\lambda_1 +
\lambda_2)^2$ and so invoves only $c_1^2$, indeed
we have, for $1 \le j \le (m/2) - 1$,
\begin{eqnarray*}s_j &=&
\sum_{i=j}^{m - j} (m^2 - m(i + j) + 2ij)
(\lambda_1 +
\lambda_2)^2\\&=&
\sum_{i=j}^{m - j} (m^2 - m(i + j) +
2ij) c_1(E)^2.\end{eqnarray*}

\bigskip
If $m$ is odd:
$$s_0 = \sum_{i=0}^{{m - 1 \over 2}} ((m - i)
\lambda_1 + i \lambda_2)(i \lambda_1 + (m - i)
\lambda_2)$$
and $s_j, j = 1, ..., {m-1 \over 2}$ are defined
as before with
$c_2(\odot^m E) = s_0 + s_1 + ... +
s_{{m-1 \over 2}}.$ By simple algebra, we have
\begin{eqnarray*}s_0 =
\sum_{i=0}^{{m-1 \over 2}} i(m
- i) c_1(E)^2 + \sum_{i=0}^{{m-1
\over 2}} (m - 2i)^2 c_2(E).
\end{eqnarray*}Thus
\begin{eqnarray*}c_{2}(\odot^{m} E) &=&
\sum_{i=0}^{{m
\over 2} - 1} (m - 2i)^2 c_2(E) + \sum_{i=0}^{{m \over 2} - 1} i(m
- i) c_1(E)^2 +\\~~&&+ \sum_{j =
1}^{{m
\over 2} - 1}
\sum_{i=j}^{m - j} (m^2 - m(i + j) + 2ij)
c_1(E)^2
\end{eqnarray*}if $m$ is even and
\begin{eqnarray*}c_{2}(\odot^{m} E) &=&
\sum_{i=0}^{{m-1
\over 2}} (m - 2i)^2 c_2(E) + \sum_{i=0}^{{m-1
\over 2}} i(m - i) c_1(E)^2 +\\~~&&+ \sum_{j =
1}^{{m - 1
\over 2}}
\sum_{i=j}^{m - j} (m^2 - m(i + j) + 2ij)
c_1(E)^2
\end{eqnarray*}if $m$ is odd. The Lemma follows by
simplifying the preceding formulas.~QED

\bigskip
\noindent
{\bf Lemma 2.2}~~{\it Let $E_i, i = 1, ..., k$ be
holomorphic vector bundles, of rank $r_i$
respectively, over a non-singular complex surface
$X$ then
\begin{eqnarray*}
&&(i)~c_{1}(\otimes_{i=1}^{k} E_i) = \sum_{i=1}^k
 (r_1 ... r_{i-1} r_{i+1} ...
r_k) c_1(E_i),\\ &&(ii)~c_{2}(\otimes_{i=1}^{k} E_i) =
\prod_{i=1}^k r_i
\sum_{i=1}^{k} ({c_2(E_i) \over r_i} - {c_1^2(E_i) \over
2r_i}) + {\prod_{i=1}^k r_i^2 \over 2} \sum_{i=1}^{k}
({c_1(E_i) \over r_i})^2.
\end{eqnarray*}}

\medskip
\noindent
{\it Proof.}~~Consider first the case $k =
2$ then by expressing formally $E_1 = L_1 \oplus ... \oplus
L_{r_1}, E_2 = F_1 \oplus ... \oplus
F_{r_2}$ as direct sums of line bundles we get
$$E_1 \otimes E_2 = \sum_{i=1}^{r_1} L_i \otimes (F_1 \oplus
...
\oplus F_{r_2})$$ hence the first Chern class is given by
$$c_1(E_1 \otimes E_2) = \sum_{i=1}^{r_1} (r_2 c_1(L_i) +
c_1(E_2)) = \sum_{i=1}^{r-1}  r_2 c_1(E_1) + r_1
c_2(E_2).$$ The case of general $k$ is similar. For
$c_2(E_1 \otimes E_2)$ we have
\begin{eqnarray*}c_2(E_1 \otimes E_2) &=&
c_2(\sum_{i=1}^{r_1} L_i
\otimes E_2) \\&=&
\sum_{i < j} c_1(L_i \otimes E_2) c_1(L_j \otimes E_2) +
\sum_{i} c_2(L_i \otimes E_2).\end{eqnarray*} The formula
of the Lemma follows from the above and the following
formulas,
$$c_l(L_i \otimes E_2) = \sum_{p=0}^{l} C_{l-p}^{r_2 - p}
c_1^{l-p} (L_p) c_i(E_2).$$ The calculation of the general
case is achieved via induction.~~QED

\bigskip
The next formulas are consequences of the preceding lemmas:

\bigskip
\noindent
{\bf Corollary 2.3}~~{\it Let $X$ be a
non-singular complex surface
$X$ then
$$c_{1}(\odot^{i_1} T^{*}X \otimes ... \otimes \odot^{i_k}
T^{*}X) = {i_1 + ... + i_k \over 2} (i_1 + 1) ... (i_k + 1),
$$
\begin{eqnarray*}\lefteqn{c_2(\odot^{i_1} T^{*}X \otimes ...
\otimes
\odot^{i_k} T^{*}X)}\\ &=& \sum_{j=1}^{k} {\prod_{l=1}^{k}
(i_l + 1) \over i_j + 1} c_2(\odot^{i_j} T^{*}X) \\&&+
\sum_{j=1}^{k} {\prod_{l \ne j} (i_l + 1) \{(\prod_{l \ne
j} (i_l + 1)) - 1\}
\over 2} c_1^2(\odot^{i_j} T^{*}X)
\\&&+
\prod_{1 \le j_1 < j_2 \le k} {\prod_{l} (i_l + 1)^2
\over (i_{j_1} + 1)(i_{j_2} + 1)} c_1(\odot^{i_{j_1}}
T^{*}X) c_1(\odot^{i_{j_2}} T^{*}X).
\end{eqnarray*}}

\bigskip
Let $m$ be a positive integer and for each fixed
positive integer $k$ denote by $q_{k}(m)$ be
the number of solutions of
 the equation $$i_1 + 2i_2 + ... + k i_k = m.$$ A
solution of the preceding equation shall be
referred to as a weighted partition of $m$ of
length $k$. It is easy to see that

\bigskip
\noindent
{\bf Lemma 2.4}~{\it With the notations above we have
$q_k(m) = p_1(m) + ... + p_k(m)$.}

\bigskip
With this we have the following asymptotic
estimate:

\bigskip
\noindent
{\bf Theorem 2.5}~{\it For $m \rightarrow
\infty$ the number $q_k(m)$ is asymptotically
given by:
$$q_k(m) \sim {m^{k-1} i^2) m^{j-1}
\over (k - 1)! (k - 1)!}.$$}

\medskip
\noindent
{\it Proof.}~~By Theorem 1.20, we have
$$p_j(m) \sim {m^{j-1} \over (j - 1)! j!}.$$By Lemma 2.4,
$$q_k(m) = \sum_{j=1}^{k} p_j(m) \sim \sum_{j = 1}^{k}
{m^{j-1} \over (j - 1)! j!} \sim {m^{k-1} i^2) m^{j-1}
\over (k - 1)! (k - 1)!}.$$QED

\bigskip
With the preceding results
the computation of the Chern numbers for
${\cal J}_k^m X$ can now be carried out by using
the Theorem of Green on Griffiths. First we
compute the Chern classes for the sheaves of each
of the weighted partitions. Then the Chern numbers
of ${\cal J}_k^m X$ is computed from these by the
following Lemma. To state the Lemma we denote by
$${\cal I}_{km} = \{I = (i_1, ..., i_k)~|~i_j \in {\bf N},
i_1 + 2i_2 + ... + ki_k = m\}.$$Moreover fixing an
ordering of the set ${\cal I}_{km}$ then

\bigskip
\noindent
{\bf Lemma 2.6}~{\it Let $X$ be a non-singular surface then
$$c_1({\cal J}_k^m X) = \sum_{I \in {\cal I}_{km}}
c_1({\cal
S}_{I}),$$
$$c_2({\cal J}_k^m X) = \sum_{I} c_2({\cal
S}_{I}) + \sum_{I < J, I, J \in {\cal I}_{km}} c_1({\cal
S}_{I}) c_1({\cal S}_{J})$$where
${\cal S}_{I} = \odot^{i_1} T^{*}X
\otimes ... \otimes \odot^{i_k} T^{*}X$.}

\medskip
\noindent
{\it Proof.}~~This is a consequence of Theorem 1.13:
$${\cal J}_{k-1}^{m} X = {\cal F}_{k}^{0}
\subset {\cal F}_{k}^{1} \subset ... \subset
{\cal F}_{k}^{[m/k]} = {\cal J}_{k}^{m} X$$
$($where $[m/k]$ is the greatest integer smaller
than or equal to $m/k)$ such that
$${\cal F}_{k}^{i}/{\cal F}_{k}^{i-1} \cong
{\cal J}_{k-1}^{m-ki} X  \otimes (\odot^{i} T^{*}
X).$$From the exact sequence
$$0 \rightarrow {\cal F}_{k}^{[m/k] - 1} \rightarrow {\cal
J}_{k}^{m} X \rightarrow {\cal J}_{k-1}^{m-k[m/k]} X \otimes
(\odot^{[m/k]} T^{*} X) \rightarrow 0$$we see that
$$c_1({\cal
J}_{k}^{m} X) = c_1({\cal F}_{k}^{[m/k] - 1}) + c_1({\cal
J}_{k-1}^{m-k[m/k]} X \otimes (\odot^{[m/k]} T^{*} X)),$$
\begin{eqnarray*}c_2({\cal
J}_{k}^{m} X) &=& c_1({\cal F}_{k}^{[m/k] - 1}) c_1({\cal
J}_{k-1}^{m-k[m/k]} X \otimes (\odot^{[m/k]} T^{*} X)) \\
&&+
c_2({\cal F}_{k}^{[m/k] - 1}) + c_2({\cal
J}_{k-1}^{m-k[m/k]} X \otimes (\odot^{[m/k]} T^{*}
X)).\end{eqnarray*}We then use filtrations of
${\cal F}_{k}^{[m/k] - 1}$ and of ${\cal
J}_{k-1}^{m-k[m/k]} X$ to compute the Chern classes.
Eventually the Chern classes are expressed by the Chern
classos of the bundles
${\cal S}_{I} = \odot^{i_1} T^{*}X
\otimes ... \otimes \odot^{i_k} T^{*}X$ for each $I \in
{\cal I}_{km}$.~~QED

\bigskip
 We shall compute the explicit numbers for the
following cases (I)
$k = 2, 1
\le m
\le 6$, (II) $k = 3, m = 6$ which will be needed
later. We shall aso compute (III) $k = m \le 5$ for
comparison with the result of section 1. We shall
write, for simplicity:
$$c_1 = c_1(T^{*}X), c_2 = c_2(T^{*}X).$$

\bigskip
\noindent
$(I_{22}) ~k = 2, m = 2$

\medskip There are two weighted partitions $P_1 =
(i_1 = 2, i_2 = 0)$ and $P_2 = (i_1 = 0, i_2 = 1)$
corresponding to the two solutions of $i_1 + 2
i_2 = 2$. The corresponding sheaves are ${\cal S}_{1} =
\odot^2 T^{*}X, {\cal S}_{2} = T^{*}X$. Denote by
$\Delta({\cal
S}_i) = c_1({\cal S}_{i}) - c_2({\cal S}_{i})$ and
$\mu({\cal S}_i) = c_1(\mu({\cal S}_i))/{\rm
rank}~\mu({\cal S}_i)$.

\medskip
\begin{eqnarray*}
\begin{tabular}{|l|l|l|l|l|l|r|} \hline
$P_i$ & ${\cal S}_i$ & rank &
$c_1({\cal S}_i)$ &$c_2({\cal S}_i)$ & $\Delta({\cal
S}_i)$ &$\mu({\cal S}_i)$\\
\hline (2, 0) & $\odot^2 T^{*}X$ & 3 & $3c_1$ &$2c_1^2 +
4c_2$ & $7c_1^2 -
4c_2$&1
\\
\hline
(0, 1) & $T^{*}X$ & 2 &$c_1$ &$c_2$ & $c_1^2 -
4c_2$&1/2 \\ \hline
\end{tabular}
\end{eqnarray*}Thus $c_1({\cal
J}_2^2 X) = 4c_1(T^*X), c_2({\cal
J}_2^2 X) = 5c_1^2(T^*X) + 5 c_2(T^*X)$, hence
$$\Delta({\cal
J}_2^2 X) = c_1^2({\cal
J}_2^2 X) - c_2({\cal J}_2^2 X) = 11 c_1^2(T^*X) -
5 c_2(T^*X),~\mu({\cal
J}_2^2 X) = 4/5.$$
We remark that the formula given in [G-G] is
$c_1^2({\cal J}_2^2X) - c_2({\cal J}_2^2X) = 7 c_1^2(T^*X) -
5 c_2(T^*X).$

\bigskip
\noindent
$(I_{23}) ~k = 2, m = 3$

\medskip
There are two weighted partitions $P_1 =
(i_1 = 3, i_2 = 0)$ and $P_2 = (i_1 = 1, i_2 = 1)$
corresponding to the two solutions of $i_1 + 2
i_2 = 3$.

\medskip
\begin{eqnarray*}
\begin{tabular}{|l|l|l|l|l|l|r|} \hline
$P_i$ & ${\cal S}_i$ & rank &
$c_1({\cal S}_i)$ &$c_2({\cal S}_i)$ & $\Delta({\cal
S}_i)$ &$\mu({\cal S}_i)$\\
\hline (3, 0) & $\odot^3 T^{*}X$ & 4 & $6c_1$ &$11c_1^2 +
10c_2$ & $25c_1^2 -
10c_2$&3/2
\\
\hline
(1, 1) & $T^{*}X \otimes T^{*}X$ & 4
&$4c_1$ &$6c_1^2 + 4c_2$ & $10c_1^2 -
4c_2$&1 \\
\hline
\end{tabular}
\end{eqnarray*}Thus $c_1({\cal
J}_2^3 X) = 10c_1(T^*X), c_2({\cal
J}_2^3 X) = 41c_1^2(T^*X) + 14 c_2(T^*X)$, hence
$$\Delta({\cal
J}_2^3 X) = 59 c_1^2(T^*X) -
14 c_2(T^*X),~\mu({\cal
J}_2^3X) = 5/4.$$

\bigskip
\noindent
$(I_{24}) ~k = 2, m = 4$

\medskip
There are 3 weighted partitions $P_1 =
(i_1 = 4, i_2 = 0), P_1 =
(i_1 = 2, i_2 = 1)$ and $P_3 = (i_1 = 0, i_2 =
2)$ corresponding to the 3 solutions of $i_1 + 2
i_2 = 4$.

\medskip
\begin{eqnarray*}
\begin{tabular}{|l|l|l|l|l|l|r|} \hline
$P_i$ & ${\cal S}_i$ & rank &
$c_1({\cal S}_i)$ &$c_2({\cal S}_i)$ & $\Delta({\cal
S}_i)$ &$\mu({\cal S}_i)$\\
\hline (4, 0) & $\odot^4 T^{*}X$ & 5 & $10c_1$ &$35c_1^2 +
20c_2$ & $65c_1^2 -
20c_2$&2
\\
\hline
(2, 1) & $\odot^2 T^{*}X \otimes T^{*}X$ & 6
&$9c_1$ &$34c_1^2 +11c_2$ & $57c_1^2 -
11c_2$&3/2 \\
\hline
(0, 2) & $\odot^2 T^{*}X$ & 3 &$3c_1$ &$2c_1^2 + 4c_2$ &
$7c_1^2 -
4c_2$&1 \\
\hline
\end{tabular}
\end{eqnarray*}Thus $c_1({\cal
J}_2^4X) = 22c_1(T^*X), c_2({\cal
J}_2^4X) = 203c_1^2(T^*X) + 35 c_2(T^*X)$, hence
$$\Delta({\cal
J}_2^4X)  = 281 c_1^2(T^*X) - 35
c_2(T^*X),~\mu({\cal
J}_2^4X) = 11/7.$$

\bigskip
\noindent
$(I_{25}) ~k = 2, m = 5$

\medskip
There are 3 weighted partitions $P_1 =
(i_1 = 5, i_2 = 0), P_1 =
(i_1 = 3, i_2 = 1)$ and $P_3 = (i_1 = 1, i_2 =
2)$ corresponding to the 3 solutions of $i_1 + 2
i_2 = 5$.

\medskip
\begin{eqnarray*}
\begin{tabular}{|l|l|l|l|l|l|r|} \hline
$P_i$ & ${\cal S}_i$ & rank &
$c_1({\cal S}_i)$ &$c_2({\cal S}_i)$ & $\Delta({\cal
S}_i)$ &$\mu({\cal S}_i)$\\
\hline (5, 0) & $\odot^5 T^{*}X$ & 6 & $15c_1$ &$85c_1^2 +
35c_2$ & $120c_1^2 -
35c_2$&5/2
\\
\hline
(3, 1) & $\odot^3 T^{*}X \otimes T^{*}X$ & 8
&$16c_1$ &$112c_1^2 +24c_2$ & $144c_1^2 -
24c_2$&2 \\
\hline
(1, 2) & $T^{*}X
\otimes \odot^2 T^{*}X$ & 6 &$9c_1$ &$34c_1^2 + 11c_2$ &
$47c_1^2 -
11c_2$&3/2 \\
\hline
\end{tabular}
\end{eqnarray*}Thus $c_1({\cal
J}_2^5X) = 40c_1(T^*X), c_2({\cal
J}_2^5X) = 750c_1^2(T^*X) + 70 c_2(T^*X)$, hence
$$\Delta({\cal
J}_2^5X) = c_1^2({\cal
J}_2^5X) - c_2({\cal J}_2^5X) = 850 c_1^2(T^*X) - 70
c_2(T^*X),~\mu({\cal
J}_2^5X) = 2.$$

\bigskip
\noindent
$(I_{26}) ~k = 2, m = 6$

\medskip
There are 4 weighted partitions $P_1 =
(i_1 = 6, i_2 = 0), P_2 =
(i_1 = 4, i_2 = 1), P_3 =
(i_1 = 2, i_2 = 1)$ and $P_4 = (i_1 = 0, i_2 =
3)$ corresponding to the 3 solutions of $i_1 + 2
i_2 = 6$.

\medskip
\begin{eqnarray*}
\begin{tabular}{|l|l|l|l|l|l|r|} \hline
$P_i$ & ${\cal S}_i$ & rank &
$c_1({\cal S}_i)$ &$c_2({\cal S}_i)$ & $\Delta({\cal
S}_i)$ &$\mu({\cal S}_i)$\\
\hline (6, 0) & $\odot^6 T^{*}X$ & 7 & $21c_1$ &$175c_1^2 +
56c_2$ & $226c_1^2 -
56c_2$&3
\\
\hline
(4, 1) & $\odot^4 T^{*}X \otimes T^{*}X$ & 10 &$25c_1$
&$330c_1^2 +45c_2$ & $295c_1^2 -
45c_2$&5/2 \\
\hline
(2, 2) & $\odot^2 T^{*}X \odot^2 T^{*}X$ & 9 &$18c_1$
&$147c_1^2 +
24c_2$ & $177c_1^2 -
24c_2$&2 \\
\hline
(0, 3) & $\odot^3 T^{*}X$ & 4 &$6c_1$ &$11c_1^2 + 10c_2$ &
$25c_1^2 -
10c_2$&3/2 \\
\hline
\end{tabular}
\end{eqnarray*}Thus $c_1({\cal
J}_2^6 X) = 70c_1(T^*X), c_2({\cal
J}_2^6 X) = 662c_1^2(T^*X) + 135 c_2(T^*X)$, hence
$$\Delta({\cal
J}_2^6 X) = 4238 c_1^2(T^*X) - 135
c_2(T^*X),~\mu({\cal
J}_2^6 X) = 7/3.$$

\bigskip
\noindent
$(II_{36}) ~k = 3, m = 6$

\medskip
There are 7 weighted partitions $P_1 =
(i_1 = 6, i_2 = 0, i_3 = 0), P_2 =
(i_1 = 4, i_2 = 1, i_3 = 0), P_3 =
(i_1 = 3, i_2 = 0, i_3 = 1), P_4 =
(i_1 = 2, i_2 = 2, i_3 = 0), P_5 =
(i_1 = 1, i_2 = 1, i_3 = 1), P_6 =
(i_1 = 0, i_2 = 3, i_3 = 0)$ and $P_7 = (i_1 = 0,
i_2 = 0, i_3 = 2)$ corresponding to the 7 solutions
of
$i_1 + 2 i_2 + 3 i_3 = 6$.

\medskip
\begin{eqnarray*}
\begin{tabular}{|l|l|l|l|l|l|r|} \hline
$P_i$ & ${\cal S}_i$ & rank &
$c_1({\cal S}_i)$ &$c_2({\cal S}_i)$ & $\Delta({\cal
S}_i)$ &$\mu({\cal S}_i)$\\
\hline (6, 0, 0) & $\odot^6 T^{*}X$ & 7 & $21c_1$ &$175c_1^2
+ 56c_2$ & $260c_1^2
- 56c_2$&3
\\
\hline
(4, 1, 0) & $\odot^4 T^{*}X \otimes T^{*}X$ & 10 &$25c_1$
&$330c_1^2 + 45c_2$ & $295c_1^2
- 45c_2$&5/2
\\
\hline
(3, 0, 1) & $\odot^3 T^{*}X \otimes T^{*}X$ & 8 &$16c_1$
&$112c_1^2 + 24c_2$ & $144c_1^2
- 24c_2$&2 \\
\hline
(2, 2, 0) & $\odot^2 T^{*}X \otimes \odot^2 T^{*}X$ & 9
&$18c_1$ &$147c_1^2 + 24c_2$ & $177c_1^2
- 24c_2$&2 \\
\hline
(1, 1, 1) & $T^{*}X \otimes T^{*}X \otimes T^{*}X$ & 8
&$12c_1$ &$66c_1^2 + 12c_2$ & $78c_1^2
- 12c_2$&3/2
\\
\hline
(0, 3, 0) & $\odot^3 T^{*}X$ & 4 &$6c_1$ &$11c_1^2 + 10c_2$
& $25c_1^2
- 10c_2$&3/2 \\
\hline
(0, 0, 2) & $\odot^2 T^{*}X$ & 3 &$3c_1$ &$2c_1^2 + 4c_2$ &
$7c_1^2
- 4c_2$&1 \\
\hline
\end{tabular}
\end{eqnarray*}us $c_1({\cal
J}_3^6 X) = 101c_1(T^*X), c_2({\cal
J}_3^6 X) = 5026c_1^2(T^*X) + 175 c_2(T^*X)$ and
$$\Delta({\cal
J}_3^6 X) = 5175 c_1^2(T^*X)
- 175 c_2(T^*X),~\mu({\cal
J}_3^6 X) = 101/49.$$

\bigskip
\noindent
$(III_{33})~k = m = 3$

\bigskip
In this case there are 3 weighted partitions: $P_1
= (3, 0, 0), P_2 = (1, 1, 0)$  and $P_3 = (0, 0,
3)$. The tabulation is given by

\medskip
\begin{eqnarray*}
\begin{tabular}{|l|l|l|l|l|l|r|} \hline
$P_i$ & ${\cal S}_i$ & rank &
$c_1({\cal S}_i)$ &$c_2({\cal S}_i)$ & $\Delta({\cal
S}_i)$ &$\mu({\cal S}_i)$\\
\hline (3, 0, 0) & $\odot^3 T^{*}X$ & 4 & $6c_1$ &$11c_1^2
+ 10c_2$ &$25c_1^2
- 10c_2$ &3/2
\\
\hline
(1, 1, 0) & $T^{*}X \otimes T^{*}X$ & 4 &$4c_1$
&$6c_1^2 + 4c_2$ & $10c_1^2
- 4c_2$&1 \\
\hline
(0, 0, 1) & $T^{*}X$ & 2 &$c_1$ &$c_2$ & $c_1^2
- c_2$&1 \\
\hline
\end{tabular}
\end{eqnarray*}Thus we have
$$c_1({\cal J}_3^3 X) = 11 c_1^2(T^*X), c_2({\cal
J}_3^3 X) = 51 c_1^2(T^*X) + 15 c_2(T^*X)$$and so
$\mu({\cal J}_3^3 X) = 11/10$ and
$$c_1^2({\cal J}_3^3 X) - c_2({\cal J}_3^3 X) = 70
c_1^2(T^*X) - 15 c_2(T^*X).$$
The formula given in [G-G]
is $c_1^2({\cal J}_3^3 X) - c_2({\cal J}_3^3 X) =
85 c_1^2(T^*X) - 49 c_2(T^*X).$

\bigskip
\noindent
$(III_{44})~k = m = 4$

\bigskip
In this case there are 5 weighted partitions: $P_1
= (4, 0, 0, 0), P_2 = (2, 1, 0, 0), P_3 = (1, 0,
1, 0), P_4 = (0, 2, 0, 0)$ and $P_5 = (0, 0, 0,
1)$. The tabulation is given by

\medskip
\begin{eqnarray*}
\begin{tabular}{|l|l|l|l|l|l|r|} \hline
$P_i$ & ${\cal S}_i$ & rank &
$c_1({\cal S}_i)$ &$c_2({\cal S}_i)$ & $\Delta({\cal
S}_i)$ &$\mu({\cal S}_i)$\\
\hline (4, 0, 0, 0) & $\odot^4 T^{*}X$ & 5 & $10c_1$ &$35c_1^2 +
20c_2$ & $65c_1^2 -
20c_2$&2
\\
\hline
(2, 1, 0, 0) & $\odot^2 T^{*}X \otimes T^{*}X$ & 6 &$9c_1$
&$34c_1^2 + 11c_2$ & $47c_1^2 -
11c_2$&3/2 \\
\hline
(1, 0, 1, 0) & $T^{*}X \otimes T^{*}X$ & 4 &$4c_1$
&$6c_1^2 + 4c_2$ & $10c_1^2 -
4c_2$&1 \\
\hline
(0, 2, 0, 0) & $\odot^2 T^{*}X$ & 3 &$3c_1$ &$2c_1^2 + 4c_2$ &
$7c_1^2 - 4c_2$&1 \\
\hline
(0, 0, 0, 1) & $T^{*}X$ & 2 &$c_1$ &$c_2$ & $c_1^2 - c_2$&1/2 \\
\hline
\end{tabular}
\end{eqnarray*}
Thus we have
$$c_1({\cal J}_4^4 X) = 27 c_1^2(T^*X), c_2({\cal
J}_4^4 X) = 338 c_1^2(T^*X) + 40 c_2(T^*X)$$and
so $\mu({\cal J}_4^4 X) = 27/20$ and
$$c_1^2({\cal J}_4^4 X) - c_2({\cal J}_4^4 X) =
391 c_1^2(T^*X) - 40 c_2(T^*X).$$

\bigskip
\noindent
$(III_{55})~k = m = 5$

\bigskip
In this case there are 3 weighted partitions: $P_1
= (5, 0, 0, 0, 9), P_2 = (3, 1, 0, 0, 0), P_3 =
(2, 0, 1, 0, 0), P_4 = (1, 2, 0, 0, 0), P_5 =
(1, 0, 0, 1, 0) P_6 = (1, 2, 0, 0, 0)$ and $P_7 =
(1, 2, 0, 0, 0)$. The tabulation is given by

\medskip
\begin{eqnarray*}
\begin{tabular}{|l|l|l|l|l|l|l|r|} \hline
$P_i$ & ${\cal S}_i$ & rank &
$c_1({\cal S}_i)$ &$c_2({\cal S}_i)$ & $\Delta({\cal
S}_i)$ &$\mu({\cal S}_i)$\\
\hline (5, 0, 0, 0, 0) & $\odot^5 T^{*}X$ & 6 &
$15c_1$ &$85c_1^2 + 35c_2$ & $140c_1^2 -
35c_2$&5/2
\\
\hline
(3, 1, 0, 0, 0) & $\odot^3 T^{*}X \otimes T^{*}X$
& 8 &$16c_1$ &$112c_1^2 + 24c_2$ & $144c_1^2 -
24c_2$&2 \\
\hline
(2, 0, 1, 0, 0) & $\odot^2 T^{*}X \otimes T^{*}X$
& 6 &$9c_1$ &$34c_1^2 + 11c_2$ & $47c_1^2 -
11c_2$&3/2
\\
\hline
(1, 2, 0, 0, 0) & $T^{*}X \otimes \odot^2 T^{*}X$
& 6 &$9c_1$ &$34c_1^2 + 11c_2$ & $47c_1^2 -
11c_2$&3/2
\\
\hline
(1, 0, 0, 1, 0) & $T^{*}X \otimes T^{*}X$ & 4
&$4c_1$ &$6c_1^2 + 4c_2$ & $10c_1^2 - 4c_2$&1 \\
\hline
(0, 1, 1, 0, 0) & $T^{*}X \otimes T^{*}X$ & 4
&$4c_1$ &$6c_1^2 + 4c_2$ & $10c_1^2 - 4c_2$&1
\\
\hline
(0, 0, 0, 0, 1) & $T^{*}X$ & 2
&$c_1$ &$c_2$ & $c_1^2 - c_2$&1\\ \hline
\end{tabular}
\end{eqnarray*}
Thus we have
$$c_1({\cal J}_5^5 X) = 58 c_1^2(T^*X), c_2({\cal
J}_5^5 X) = 1622 c_1^2(T^*X) + 90 c_2(T^*X)$$and
so $\mu({\cal J}_5^5 X) = 29/18$ and
$$c_1^2({\cal J}_5^5 X) - c_2({\cal J}_4^4 X) =
1742 c_1^2(T^*X) - 90 c_2(T^*X).$$
We remark that the inequality (1.21) in [G-G]
is incorrect (for example set $k = 2$ or $k = 3$ and
compare these to the formulas obtained above; indeed
for $k = 2$ (see (1.21) in [G-G]) reduces to
$c_1^2(X) - c_2(X) > 0$).

\bigskip
\bigskip
\noindent
{\bf \S~3~~Weighted Projective Spaces and
Projectivized Jet Bundles}

\bigskip
For a vector bundle, e.g., the $k$-jet bundle
$T^kX$, a standard approach of studying the
bundle is to projectivized it and then study the
line bundles over the projectivization. We are
going to do the same for the ${\bf
C}^{*}$-bundle $J^{k}X$ using the well-known results in the
former case as a guide. The fiber of
the projectivized bundle are certain types of
weighted projective space. Thus we shall first
recall some basic facts about weighted
projective spaces. For more detailed discussions and
further references the readers are referred to the articles
[B-R], [Do] and the monograph [Di].

\bigskip
Let $Q = (q_{0}, q_{1}, ..., q_{r})$ ($r
\ge 1$) be an $(r + 1)$-tuple of positive
integers. The tuple $Q$ is said to be {\sl reduced} if
the greatest common divisor (gcd) of $(q_{0}, q_{1},
.... q_{r})$ is 1. In general if the gcd is $d$ the tuple
$$Q_{{\rm red}} = Q/d = (q_0/d, ..., q_r/d)$$ is called the
reduction of
$Q$. Let
$d_0 = gcd(q_{1}, ..., q_{r}), d_r = gcd(q_{0}, ...,
q_{r-1})$ and
$$d_i = {\rm gcd}(q_{0}, q_{1}, ..., q_{i-1}, q_{i+1},
..., q_{r}),~~1 \le i \le r - 1.$$ Let $a_0 = lcm (d_{1},
..., d_{r}), a_r = lcm(d_{0}, ..., d_{r-1})$ and
$$a_i = {\rm lcm}(d_{0}, d_{1}, ..., d_{i-1},d_{i+1},
..., d_{r}),~~1 \le i \le r - 1$$where $lcm$ means least
common multiple. Define the normalization of $Q$ by
$$Q_{{\rm norm}} = (q_0/a_0, ..., q_r/a_r).$$ A tuple $Q$ is
said to be {\sl normalized} if $Q = Q_{{\rm norm}}$.

\bigskip
Let
$({\bf C}^{r+1}, Q)$ be the
$(r+1)$-dimensional complex vector space such
that the variable $z_{i}$ is assigned the weight
(or degree) $q_{i}$. A ${\bf C}^{*}$-action is
defined on
$({\bf C}^{r+1}, Q)$ by:
\begin{eqnarray}
\lambda.(z_{0}, ..., z_{r})
= (\lambda^{q_{0}}z_{0}, ...,
\lambda^{q_{r}}z_{r}),~\lambda \in {\bf C}^{*}.
\end{eqnarray}
The quotient space, ${\bf P}(Q) = ({\bf C}^{r+1},
Q)/{\bf C}^{*}$, is called the weighted
projective space of type $Q$. The equivalence
class of an element $(z_{0}, ..., z_{r})$ is
denoted by $[z_{0}, ..., z_{r}]_{Q}$. For $Q =
(1, ..., 1) = {\bf 1}, {\bf P}(Q) = {\bf P}^{r}$
is the usual complex projective space of
dimension $r$ and an element of
${\bf P}^{r}$ is denoted simply by $[z_{0},
..., z_{r}]$. Indeed for the special case $r = 1$ it
can be shown that, for any tuple $(q_0, q_1)$, ${\bf P}(q_0,
q_1)
\cong {\bf P}^1$. This is not so if $r \ge 2$, however,
we do have:

\bigskip \noindent
{\bf Theorem 3.1}~{\it Let $Q = (q_0, ..., q_r)$ be an
$(r+1)$-tuple of positive integers then $${\bf P}(Q) \cong
{\bf P}(Q_{{\rm red}}) \cong {\bf P}(Q_{{\rm norm}}).$$}

\bigskip
\noindent
{\bf Example 3.2}~~It is clear that a normalized tuple is
reduced. The converse is not true in general. Let
$Q = (4, 6, 12)$ then
$Q_{{\rm red}} = (2, 3, 6)$ is reduced but is not
normalized. In fact $Q_{{\rm norm}} = (Q_{{\rm
red}})_{{\rm norm}} = (1, 1, 6)$. The tuple $(6, 10, 15)$
is reduced but is not normalized, in fact its normalization
is $(1, 1, 1)$ hence ${\bf P}(6, 10, 15) \cong {\bf P}^2$.

\bigskip
Define a map $\rho_{Q} : ({\bf C}^{r+1}, {\bf 1})
\rightarrow ({\bf C}^{r+1}, Q)$ by
\begin{eqnarray}
\rho_{Q}(z_{0}, ..., z_{r}) = (z_{0}^{q_{0}}, ...,
z_{r}^{q_{r}}).
\end{eqnarray}
It is easily seen that $\mu_{Q}$ is compatible
with the respective ${\bf C}^{*}$-actions and
hence descends to a well-defined morphism:
\begin{eqnarray}
\bar{\rho}_{Q} : {\bf P}^{r} \rightarrow {\bf P}(Q),~
\bar{\rho}_{Q}([z_{0}, ..., z_{r}]) = [z_{0}^{q_{0}}, ...,
z_{r}^{q_{r}}]_{Q}.
\end{eqnarray} The weighted projective
space can aso be described as follows. Denote by
$\Theta_{q_{i}}$ the group of
$q_{i}$-th roots of unity. Then the group
$\Theta_{Q} = \oplus_{i=0}^{r} \Theta_{q_{i}}$
acts on ${\bf P}^{r}$ by coordinate wise
multiplication:
$$(\theta_{0}, ..., \theta_{r}).[z_{0},
..., z_{r}] = [\theta_{0} z_{0}, ..., \theta_{r}
z_{r}],~ \theta_{i} \in \Theta_{q_{i}}$$ and it is easily
verified that
${\bf P}(Q) = {\bf
P}^{r}/\Theta_{Q}.$

\bigskip
\noindent
{\bf Theorem 3.3}~~{\it The weighted
projective space
${\bf P}(Q)$ is isomorphic to the quotient ${\bf
P}^{r}/\Theta_{Q}$. In particular, ${\bf P}(Q)$
is irreducible and normal {\rm (}the
singularities are cyclic quotients and hence
rational{\rm )}.}

\medskip
Denote by $S_{Q}(m)$ the space of
homogeneous polynomials of degree $m > 0$ in the
variables
$z_{i}$ (assigned with the degree $q_{i}$).
In other words, a polynomial $P$ is in $S(Q)(m)$
if
$$P(\lambda.(z_{0}, ..., z_{r})) = \lambda^{m}
P(z_{0}, ..., z_{r}).$$ We may express such a
polynomial explicitly:
\begin{eqnarray}
P = \sum_{(i_{0}, ..., i_{r}) \in {\cal I}_{Q,m}} a_{i_{0} ... i_{r}}
z_{0}^{i_{0}} ...
z_{r}^{i_{r}}
\end{eqnarray}
where  the index set ${\cal I}_{Q,m}$ is
defined by:
$${\cal I}_{Q,m} = \{(i_{0}, ...,
i_{r})~|~\sum_{j=0}^{r} q_{j}i_{j} = m\}.$$ The
sheaf ${\cal O}_{{\bf P}(Q)}(m)$ is the sheaf
over ${\bf P}(Q)$ whose global regular sections
are precisely the elements of $S_{Q}(m)$:
\begin{eqnarray}
H^{0}({\bf P}(Q), {\cal O}_{{\bf P}(Q)}(m))
= S_{Q}(m).
\end{eqnarray}
 For negative integer $- m, m > 0$ the
sheaf ${\cal O}_{{\bf P}(Q)}(- m)$ is defined
to be the dual of ${\cal O}_{{\bf P}(Q)}(m).$

\bigskip
\noindent
{\bf Theorem 3.4}~~{\it {\rm (i)} For any $m
\in {\bf Z}, {\cal O}_{{\bf P}(Q)}(m)$ is a reflexive
coherent sheaf. {\rm (ii)} The sheaf ${\cal O}_{{\bf
P}(Q)}(m)$ is locally free if $m$ is divisible
by each
$q_{i}~($hence by the least common multiple$
)$. {\rm (iii)} Let $m_{Q}$ be the least common
multiple of $\{q_{0}, ..., q_{r}\}$ then
${\cal O}_{{\bf P}(Q)}(m_{0})$ is ample. {\rm
(iv)} There exists an interger $n_{0}$ depending
only on $Q$ such that ${\cal O}_{{\bf
P}(Q)}(nm_{Q})$ is very ample for all
$n \ge n_{0}$. {\rm
(v)} For any $\alpha, \beta \in {\bf Z}$ we have ${\cal
O}_{{\bf P}(Q)}(\alpha m_{Q}) \otimes {\cal
O}_{{\bf P}(Q)}(\beta) \cong {\cal
O}_{{\bf P}(Q)}(\alpha m_{Q} + \beta)$.}

\bigskip
For any subset $J \subset \{0, 1,..., r\}$
denote by $m_{J}$ the least common multiple of
$\{q_{j}, j \in J\}$ and define
$$m(Q) = - |Q| + {1 \over r} \sum_{\nu =
2}^{r+1} {\sum_{\#J = \nu} m_{J}
\over C_{\nu - 2}^{r-1}}$$ where $C_{a}^{b}$ is
the usual binomial coefficient and $|Q| = q_{0} +
... + q_{r}$. It is known that assertion (iv) holds if
$n > m(Q)$. In general the line sheaf ${\cal O}_{{\bf
P}(Q)}(m)$ is not invertible if $m$ is not an integer
multiple of $m_{Q}$. It can be shown that for $Q = (1, 1,
2)$ the sheaf  ${\cal O}_{{\bf
P}(Q)}(1)$ is not invertible and hence, neither is ${\cal
O}_{{\bf P}(Q)}(1) \otimes {\cal O}_{{\bf P}(Q)}(1)$. This
also shows that
${\cal O}_{{\bf P}(Q)}(1) \otimes {\cal
O}_{{\bf P}(Q)}(1) \not\cong {\cal
O}_{{\bf P}(Q)}(2)$ as ${\cal
O}_{{\bf P}(Q)}(2)$ is invertible by part (ii) of the
preceding Theorem.

\bigskip
\noindent
{\bf Theorem 3.5}~~{\it Let $Q$ be a $(r + 1)$-tuple of
positive integers then
\begin{eqnarray*}&&(i)~~H^{i}({\bf P}(Q),
{\cal O}_{\bf P}(Q)(p)) = \{0\}, ~p \in {\bf Z}~{\rm
if}~i \ne 0, r;\\
&&(ii)~H^{0}({\bf P}(Q),
{\cal O}_{\bf P}(Q)(p)) = S_{Q}(p) ~p \in {\bf
Z};\\&&(iii)~H^{r}({\bf P}(Q), {\cal O}_{\bf P}(Q)(p)) \cong
S(Q)(- p - |Q|), ~p \in {\bf
Z}\end{eqnarray*}where $|Q| = q_0 + .... + q_r$.}

\bigskip
Denote by $Pic~({\bf P}(Q))$ and $Cl({\bf P}(Q))$ the
Picard group and respectively the divisor class group.

\bigskip
\noindent
{\bf Theorem 3.6}~~{\it Let $Q = Q_{{\rm norm}}$ be a
normalized $(r + 1)$-tuple of positive integers then
$(i)~Pic~({\bf P}(Q)) \cong {\bf Z}$
is generated by
$[{\cal O}_{\bf P}(Q)(m_{Q})]; ~(ii)~Cl({\bf P}(Q)) \cong
{\bf Z}$
is generated by
$[{\cal O}_{\bf P}(Q)(1)].$}

\bigskip
Let $Q$ be a
$(r + 1)$-tuple of positive integers define for $k = 1,
..., r$:
\begin{eqnarray*}l_{Q,k} = lcm ~\{{q_{i_0} ... q_{i_k}
\over gcd~(q_0, ..., q_{i_k})}~|~0 \le i_0
< .... < i_k \le r\}.\end{eqnarray*}

\bigskip
\noindent
{\bf Theorem 3.7}~~{\it Let $Q$ be a
$(r + 1)$-tuple of positive integers then
\begin{eqnarray*}H^{i}({\bf P}(Q); {\bf Z}) \cong
\left\{\begin{array}{ll} {\bf
Z}, &\mbox{if $i$ is even,}\\
0, &\mbox{if $i$ is odd.}
\end{array}
\right.\end{eqnarray*}
Moreover, let $\bar{\rho}_{Q} : {\bf P}^r \rightarrow {\bf
P}(Q)$ be the quotient map as defined by $(28)$ then the
following diagram commutes,
\begin{eqnarray*}
&&H^{2k}({\bf P}(Q); {\bf Z})
~~\stackrel{\bar{\rho}_{Q}^{*}}{\longrightarrow}
~~H^{2k}({\bf P}^r; {\bf Z})\\
&&~~~~\cong \downarrow ~~~~~~~~~~~~~~~~~~~\cong \downarrow\\
&&~~~~~~~{\bf Z}
~~~~~~~~~\stackrel{l_{Qk}}{\longrightarrow}~~~~~~~
{\bf Z}
\end{eqnarray*}where the lower map is the multiplication by
the number $l_{Qk}$.}

\bigskip
Note that the number $l_{Qr}$ is precisely the number of
preimages of a point in ${\bf P}(Q)$ under the quotient
map $\bar{\rho}_{Q}$. The proof of the preceding Theorem
for $k = r$ is quite easy. for the general case we refer
the readers to [Ka]. We shall only be concerned with the
case where
$n, k \ge 1$ are positive integers and
\begin{eqnarray*}
Q = ((\underbrace{1, ...,
1}_{n}), (\underbrace{2, ..., 2}_{n}), ...,
(\underbrace{k, ..., k}_{n})).
\end{eqnarray*}
In this case we shall write ${\bf P}_{n,k}$
for ${\bf P}(Q)$. Note that
$r =$ dim
${\bf P}_{n,k} = nk - 1$ In this case the least common
multiple of $Q$ is $m_{Q} = k!$ and $l_{Qr} = (k!)^n$.

\bigskip
Let $\pi : ({\cal E}, h) \rightarrow X$ be a
holomorphic hermitian vector bundle over a
compact K\"{a}hler manifold $X$. Denote by ${\cal
L}({\cal E})$ be the "hyperplane bundle" defined
over the projectivized bundle ${\bf P}({\cal
E})$. It is defined as follows:
\begin{eqnarray*}
&&\pi^{*}{\cal E} ~~~~\longrightarrow~~~~ {\cal E}\cr
&&~\downarrow pr ~~~~~~~~~~~\downarrow p
\cr &&{\bf
P}({\cal E})~~
\stackrel{\pi}{\longrightarrow}~~~~
X
\end{eqnarray*}
the tautological
sub-sheaf is defined by:
$$\{((x, [\xi]), \eta) \in \pi^{*}{\cal E}~|~(x,
[\xi])
\in {\bf P}({\cal E}),~ p([\xi]) = x,~ [\eta] =
[\xi]\}$$
and
${\cal L}_{k}$ is defined to be the dual of the
tautological line bundle. In other words, since
the fiber ${\bf P}({\cal E})$ over a point $x
\in X$ is a projective space, the
restriction of ${\cal L}_{k}X$ to ${\bf
P}({\cal E}_{x})$ is the hyperplane line bundle
${\cal O}_{{\bf P}^{r-1}}(1)$ (here $r =$ rank
${\cal E}$).
We shall often use the notation ${\cal
O}_{{\bf P}({\cal E})}(1)$ for ${\cal L}_{k}$ and
the tensor product ${\cal L}_{k}^{m}$ by ${\cal
O}_{{\bf P}({\cal E})}(m)$ for any integer $m \in
{\bf Z}$.
The following is a classical Theorem of
Grothendieck:

\bigskip
\noindent
{\bf Theorem 3.8}~~{\it Let ${\cal E}$ be a
holomorphic vector bundle over a
complex manifold $X$ then for any $m, j \ge 0$, the
j-th direct image sheaf of the $m$-fold tensor product of
${\cal L}(m)$ is isomorphic to the $m$ fold symmetric
product of
$E$, i.e.,
$R^j_{*}{\cal L}^m({\cal E})
\cong
\odot^m {\cal E}$ and
$$H^{j}(X, \odot^m {\cal E} \otimes {\cal S}) \cong
H^{q}(X, {\cal L}^m({\cal E})
\otimes p^{*} {\cal S})$$ where ${\cal S}$ is any
sheaf on $X$.}

\bigskip
Let $\pi : J^{k}X \rightarrow X$ be the
(restricted) $k$-jet bundle of a complex manifold
$X$. Denote by ${\cal L}_{k}$ the
"hyperplane sheaf" defined over the
projectivized $k$-jet bundle ${\bf
P}(J^{k}X)$. It is defined as follows.
Consider the commutative diagram:
\begin{eqnarray*}
&&~\pi^{*}J^{k}X ~~~~\longrightarrow~~
J^{k}X\cr &&~~~\downarrow pr
~~~~~~~~~~~~~~~\downarrow p
\cr &&{\bf
P}(J^{k}X)~~~
\stackrel{\pi}{\longrightarrow}~~~~
~X
\end{eqnarray*}
the tautological
sub-sheaf is defined by:
$$\{((x, [\xi]), \eta) \in \pi^{*}J^{k}X~|~(x,
[\xi])
\in {\bf P}(J^{k}X),~ p([\xi]) = x,~ [\eta] =
[\xi]\}$$
and
${\cal L}_{k}$ is defined to be the dual the
tautological line sheaf. In other words, since
the fiber ${\bf P}(J_{x}^{k}X)$ over a point $x
\in X$ is a  weighted projective space of type $Q
= ((1, ..., 1); ...; (k, ..., k))$ the
restriction of ${\cal L}_{k}X$ to ${\bf
P}(J_{x}^{k}X)$ is the line sheaf ${\cal
O}_{{\bf P}(Q)}(1)$ as defined in the
preceding section. We shall use the
notation ${\cal
O}_{{\bf P}(J^{k}X)}(1)$ for ${\cal L}_{k}$. More generally
for
any integer $m, {\cal O}_{{\bf P}(J^{k}X)}(m)$ is the sheaf
on ${\bf P}(J^{k}X)$ which restricts to the bundle ${\cal
O}_{{\bf P}(Q)}(m)$ along each fiber of the projection map
$p : {\bf P}(J^k X) \rightarrow X$. The proof of
the preceding Theorem relies on the classical
Vanishing Theorem of cohomologies on projective
spaces. The analogoue of this for weighted
projective spaces is provided by Theorem 3.3 and
hence we have (see [G-G] and [K-O]):

\bigskip
\noindent
{\bf Theorem 3.9}~~{\it Let $X$ be a complex manifold and
${\cal S}$ be a sheaf over $X$ then for any $m, j \ge 0$ we
have $R^j_{*}{\cal O}_{{\bf
P}(J^{k}X)}(m) \cong {\cal J}_{k}^{m}X$ and
$$H^{j}(X, {\cal J}_{k}^{m}X \otimes {\cal
S}) \cong H^{j}({\bf P}(J^{k}X), {\cal O}_{{\bf
P}(J^{k}X)}(m)
\otimes p^{*}{\cal S}).$$}

\bigskip
The following is a consequence of Theorem 3.4:

\bigskip
\noindent
{\bf Theorem 3.10}~~{\it Let $X$ be a complex manifold
then

\smallskip
{\rm (i)} for any $m
\in {\bf Z}, {\cal O}_{{\bf P}(J^kX)}(m)$ is a reflexive
coherent sheaf$;$

\smallskip{\rm (ii)} the sheaf ${\cal O}_{{\bf
P}(J^kX)}(m)$ is locally free if $m$ is divisible
by each
$q_{i}~($hence by the least common multiple $
k!);$

\smallskip
{\rm
(iii)} for any $\alpha, \beta \in {\bf Z}, {\cal
O}_{{\bf P}(J^kX)}(k! \alpha) \otimes {\cal
O}_{{\bf P}(J^kX)}(\beta) \cong {\cal
O}_{{\bf P}(J^kX)}(k! \alpha + \beta)$.}

\bigskip
Due to the fact that ${\cal
O}_{{\bf P}(J^kX)}(1)$ is not locally free and that, in
general, ${\cal
O}_{{\bf P}(J^kX)}(a) \otimes {\cal
O}_{{\bf P}(J^kX)}(a) \cong {\cal
O}_{{\bf P}(J^kX)}(a + b)$ some of the proof of the results
that are valid on projectivized vector bundle are not valid
even though modifications of the results can be obtained
via alternative proofs. We establish some of the results
(the counterparts in the case of projectivized vector
bundle are well-known) that will be essential in the next
section.

\bigskip
\noindent
{\bf Lemma 3.11}~~{\it Let $X$ be a complex manifold of
dimension $n$ and let $p : {\bf P}(J^kX)
\rightarrow X$ be the projection map. Then the natural
morphism:
$$\phi : p^{*}p_{*} {\cal O}_{{\bf P}(J^kX)}(k!)
\rightarrow {\cal O}_{{\bf P}(J^kX)}(k!)$$ is surjective and
$$\sum_{i=0}^n (-1)^{i} c_1^{r-i}({\cal O}_{{\bf
P}(J^kX)}(k!)) ~.~ p^{*} c_i({\cal J}_k^{k!} X) = 0$$
 where ${\cal F}$ is the kernel of $\phi$
and $r = nk - 1$ is the fiber dimension of $p$.}

\medskip
\noindent
{\it Proof.}~~For simplicity we write ${\cal O}(k!)$ for
${\cal O}_{{\bf P}(J^kX)}(k!)$. By definition the
restriction of ${\cal O}(k!)$ to a fiber of the projection
map is ${\cal O}_{{\bf P}(Q)}(k!)$ where $$Q = ((1, ...,
1); ...; (k, ..., k)).$$ Thus the least common multiple of
the indices is
$k!$ and ${\cal O}_{{\bf P}(Q)}(k!)$ is ample by Theorem
3.4. This implies that the map $\phi$ is surjective (see
for example [B-S]). By Theorem 3.9
$p^{*}p_{*} {\cal O}(k!) = p^{*}{\cal J}_k^{k!} X$ and so
the sequence:
$$0 \rightarrow {\cal F} \rightarrow p^{*}{\cal J}_k^{k!}
X \rightarrow {\cal O}(k!) \rightarrow 0$$ is eact. By
Whitney's formula
$$\sum_{i=0}^{r} p^{*} c_i({\cal J}_k^{k!} X) = (1
+ c_1({\cal O}(k!)) . \sum_{i=0}^{r-1}
c_i({\cal F})$$ and hence
$$p^{*} c_i({\cal J}_k^{k!} X) = c_1({\cal O}(k!)) .
c_{i-1}({\cal F}) + c_i({\cal F})$$for $0 \le i \le r$ with
$c_{-1}({\cal F}) = c_r({\cal F}) = 0$ (as rank ${\cal F} = r
- 1$). We can eliminate the Chern classes of ${\cal F}$ by
first multiplying the preceding identity by
$c_1^{r-i}({\cal O}(k!))$ and then take alternating
sum; namely:
\begin{eqnarray*}
c_1^{r-1}({\cal O}(k!)) ~.~ p^{*} c_1({\cal J}_k^{k!} X)
&=& c_1^{r}({\cal O}(k!)) + c_1^{r-1}({\cal O}(k!)) ~.~
c_1({\cal F})\\
c_1^{r-2}({\cal O}(k!)) ~.~ p^{*}c_2({\cal J}_k^{k!} X) &=&
c_1^{r-1}({\cal O}(k!)) ~.~ c_1({\cal F}) + c_1^{r-2}({\cal
O}(k!)) . c_2({\cal F})\\
c_1^{r-3}({\cal O}(k!)) ~.~p^{*}c_3({\cal J}_k^{k!} X) &=&
c_1^{r-1}({\cal O}(k!)) ~.~ c_2({\cal F}) + c_1^{r-3}({\cal
O}(k!)) ~.~ c_3({\cal F})\\ &...&\\
c_1({\cal O}(k!)) ~.~ p^{*}c_{r-1}({\cal
J}_k^{k!} X) &=& c_1^2({\cal O}(k!)) ~.~ c_{r-2}({\cal F})
+ c_1({\cal O}(k!)) . c_{r-1}({\cal F})\\
p^{*}c_r({\cal J}_k^{k!} X) &=&
c_1({\cal O}(k!)) ~.~ c_{r-1}({\cal
F})
\end{eqnarray*} and multiply the $i$-identity above by
$(-1)^{r}$ and then taking the sum from $i = 1$ to $i = r$
yields
$$\sum_{i=1}^r (-1)^{i} c_1^{r-i}({\cal O}(k!)) ~.~
p^{*} c_i({\cal J}_k^{k!} X) = - c_1^{r}({\cal O}(k!)).$$
Moving the RHS to the LHS yields the identity of the
Lemma. ~~QED

\bigskip
Note that $c_i({\cal J}_k^{k!} X) = 0$ if $i \ge n =$ dim
$X$.

\bigskip
\noindent
{\bf Lemma 3.12}~~{\it \it Let $X$ be a compact complex
manifold of complex dimension $n$ then for any $x \in X$,
$$\int_{{\bf P}(J^{k}X)_x} c_1^{nk-1}({\cal O}_{{\bf
P}(J^kX)}(k!)|_{{\bf P}(J^{k}X)_x}) = (k!)^n$$ where
${\bf P}(J^{k}X)_x$ is the fiber over $x$.}

\medskip
\noindent
{\it Proof.}~~By definition the fiber, ${\bf
P}(J^{k}X)_x$, over any point $x \in
X$ of the projection map  of
$p : {\bf P}(J^{k}X) \rightarrow X$ is the weighted
projective space $${\bf P}((\underbrace{1, ...,
1}_{n}), (\underbrace{2, ..., 2}_{n}), ...,
(\underbrace{k, ..., k}_{n}))$$
of dimension $nk - 1$. By Theorem 3.7 the quotient map
$\bar{\rho}_{Q} : {\bf P}^{nk-1} \rightarrow {\bf P}(Q)$ is
a finite morphism with sheet number $l_{Q,nk-1} =
(k!)^{n}$. The generator of
$H^{2(nk-1)}({\bf P}^{nk-1}; {\bf Z})$ is represented by the
$(nk-1)$-th power,
$\omega_{FS}^{nk-1}$, of the the Fubini-Study
metric $\omega_{FS} = c_1({\cal O}_{{\bf
P}^{nk-1}}(1)$. The Lemma follows readily as we have:
$$\int_{{\bf P}^{nk-1}} c_1^{nk-1}({\cal O}_{{\bf
P}^{nk-1}}(1)) = \int_{{\bf P}^{nk-1}} \omega_{FS}^{nk-1} =
1.$$ QED

\bigskip
\noindent
{\bf Theorem 3.13}~~{\it Let $X$ be a compact complex
manifold of complex dimension $n$ then the following
intersection formulas hold$:$
$$c_1^{nk + j - 1}({\cal O}_{{\bf
P}(J^kX)}(k!)) ~.~p^{*}D_1 ~.~\cdots ~.~p^{*}D_{n-j} =
(k!)^n \Delta_j ~.~D_1 ~.~\cdots~.~D_{n-j} $$ for
divisors $D_1, ..., D_{n-j}, j = 0, 1, ..., n$ on $X$. The
numbers $\Delta_j$ is defined by setting $\Delta_0 = 1,
\Delta_1 = c_1({\cal
J}_k^{k!}X)$ and by the recursive relation$:$
$$\Delta_j = \sum_{i=1}^j (-1)^{i+1} \Delta_{j-i} ~.~
c_i({\cal J}_k^{k!}X),~~j \ge 2.$$}

\medskip
\noindent
{\it Proof.}~~Note that dim ${\bf P}(J^kX) = n(k+1) - 1$
and the fiber dimension, dim ${\bf P}(J^kX)_x = nk - 1$.
Thus, by fiber integration (Lemma 3.12),
$$c_1^{nk - 1}({\cal O}_{{\bf
P}(J^kX)}(k!)) ~.~p^{*}D_1 ~.~\cdots~.~p^{*}D_{n} =
(k!)^n \Delta_0 D_1 ~.~\cdots~.~D_{n}$$which is the case
$j = 0$. By Lemma 3.11 with
$r = nk - 1$,
\begin{eqnarray}\sum_{i=0}^r (-1)^{i} c_1^{r-i}({\cal
O}_{{\bf P}(J^kX)}(k!)) ~.~ p^{*} c_i({\cal J}_k^{k!} X) =
0\end{eqnarray}
and, multiplying by
$p^{*}D_1~.~
\cdots ~.~p^{*}D_{n-1}$, we get
\begin{eqnarray*}&&c_1^{nk-1}({\cal O}_{{\bf
P}(J^kX)}(k!))~.~p^{*}D_1~.~
\cdots ~.~p^{*}D_{n-1} \\&&= c_1^{nk-2}({\cal O}_{{\bf
P}(J^kX)}(k!))~.~p^{*}c_1({\cal J}_k^{k!} X)~.~p^{*}D_1~.~
\cdots ~.~p^{*}D_{n-1}\end{eqnarray*} as the rest of
the terms vanish for dimension reason. Multiplying the
above by $c_1({\cal O}_{{\bf
P}(J^kX)}(k!))$ yields,
\begin{eqnarray*}&&c_1^{nk}({\cal O}_{{\bf
P}(J^kX)}(k!))~.~p^{*}D_1~.~
\cdots ~.~p^{*}D_{n-1} \\&&= c_1^{nk-1}({\cal O}_{{\bf
P}(J^kX)}(k!))~.~p^{*}c_1({\cal J}_k^{k!} X)~.~p^{*}D_1~.~
\cdots ~.~p^{*}D_{n-1}.
\end{eqnarray*}Fiber integration shows that the term on
the right above equals $$(k!)^n c_1({\cal J}_k^{k!}
X)~.~D_1~.~
\cdots ~.~D_{n-1} = (k!)^n \Delta_1~.~D_1~.~
\cdots ~.~D_{n-1}.$$ This establish the Theorem for
the case $j = 1$.

\bigskip
If we multiply (31) by
$p^{*}D_1~.~
\cdots ~.~p^{*}D_{n-2}$ we are left with 3 terms (again
for dimension reason):
\begin{eqnarray*}&&c_1^{nk-1}({\cal O}_{{\bf
P}(J^kX)}(k!))~.~p^{*}D_1. \cdots
.p^{*}D_{n-2} \\&&= c_1^{nk-2}({\cal O}_{{\bf
P}(J^kX)}(k!))~.~p^{*}c_1({\cal J}_k^{k!}
X) ~.~p^{*}D_1~.~
\cdots ~.~p^{*}D_{n-2} \\&&~~-
c_1^{nk-3}({\cal O}_{{\bf
P}(J^kX)}(k!))~.~p^{*}c_2({\cal J}_k^{k!}
X)~.~p^{*}D_1~.~
\cdots ~.~p^{*}D_{n-2}. \end{eqnarray*}
Now multiply the above by $c_1^{2}({\cal O}_{{\bf
P}(J^kX)}(k!))$ then the LHS above is given by
$$c_1^{nk+1}({\cal O}_{{\bf
P}(J^kX)}(k!))~.~p^{*}D_1. \cdots
~.~p^{*}D_{n-2} $$while the first term on the
RHS is given by (the case $j = 1$):
\begin{eqnarray*}\lefteqn{c_1^{nk}({\cal O}_{{\bf
P}(J^kX)}(k!))~.~p^{*}c_1({\cal J}_k^{k!}
X) ~.~p^{*}D_1~.~
\cdots ~.~p^{*}D_{n-2}} \\&=&
(k!)^n \Delta_1~.~c_1({\cal J}_k^{k!} X)~.~D_1~.~
\cdots ~.~D_{n-2}\end{eqnarray*}and the second term on the
right is given by (the case $j = 0$):
\begin{eqnarray*}\lefteqn{- c_1^{nk-1}({\cal O}_{{\bf
P}(J^kX)}(k!))~.~p^{*}c_2({\cal J}_k^{k!}
X)~.~p^{*}D_1~.~
\cdots ~.~p^{*}D_{n-2}}\\&=& - (k!)^n \Delta_0
c_2({\cal J}_k^{k!} X)~.~D_1~.~
\cdots ~.~D_{n-2}\end{eqnarray*} Combining the above
yields the case $j = 2$:
\begin{eqnarray*}\lefteqn{c_1^{nk+1}({\cal O}_{{\bf
P}(J^kX)}(k!))~.~p^{*}D_1. \cdots
~.~p^{*}D_{n-2}}\\&=& (k!)^n \{\Delta_1~.~c_1({\cal
J}_k^{k!} X) - \Delta_0~
c_2({\cal J}_k^{k!} X)\}~.~D_1~.~
\cdots ~.~D_{n-2}\\&=& (k!)^n \Delta_2~.~D_1~.~
\cdots ~.~D_{n-2}
\end{eqnarray*}as, by definition, $$\Delta_2 =
\Delta_1~.~c_1({\cal J}_k^{k!} X) - \Delta_0~
c_2({\cal J}_k^{k!} X) = c_1^2({\cal J}_k^{k!} X) -
c_2({\cal J}_k^{k!} X).$$
Thus the case $j = 2$ is also established.
Inductively, the prcedure above yields:
\begin{eqnarray*}
\lefteqn{
c_1^{n-k+j-2}({\cal O}_{{\bf
P}(J^kX)}(k!))~.~p^{*}D_1~.~ \cdots ~.~p^{*}D_{n-3}} \\
 & = & \sum_{i=1}^{j} (-1)^{i-1} \Delta_{j-i}~.~c_i({\cal
J}_k^{k!} X)~.~D_1~.~ \cdots
~.~D_{n-3}\\
    & = & (k!)^n \Delta_{j}~.~D_1~.~ \cdots ~.~D_{n-3}.
\end{eqnarray*}
QED

\bigskip

\bigskip
\noindent
{\bf Theorem 3.14}~~{\it Let $X$ be a non-singular
projective surface and assume that $(i)~c_1^2({\cal J}_k^{k!}
X) - c_2({\cal J}_k^{k!} X) > 0$ and $(ii)~h^2({\cal
J}_k^{k!m}) = O(m^{(n + 1)k - 2)!})$. Then
${\cal J}_k^{k!} X$ is big.}

\medskip
\noindent
{\it Proof.}~~Let ${\bf P}(J^kX)$ be the projectivized
$k$-jet bundle. Then
$\dim {\bf P}(J^kX) = (n + 1)k - 1$. Riemann-Roch
applied to the line bundle ${\cal O}_{{\bf P}(J^kX)}(k!)$
yields
$$\chi ({\cal O}_{{\bf P}(J^kX)}(k!m)) = {c_1^{(n + 1)k -
1}({\cal O}_{{\bf P}(J^kX)}(k!)) \over ((n + 1)k - 1)!}
m^{(n + 1)k - 1} + O(m^{(n + 1)k - 2}). $$ Theorem 3.13 and
assumption (i) imply that there exists positive constant
$c > 0$ and positive integer $m_0^{'}$ such that
\begin{eqnarray*}\chi ({\cal O}_{{\bf P}(J^kX)}(k!m)) &=&
{c_1^2({\cal J}_k^{k!} X) - c_2({\cal J}_k^{k!} X) \over ((n
+ 1)k - 1)!} m^{(n + 1)k - 1} + O(m^{(n + 1)k -
2})\\&\ge& cm^{(n + 1)k -
1}\end{eqnarray*}for all
$m \ge m_0'$. Theorem 3.8 implies that the same is true
for
${\cal J}_k^{k!m}$ i.e. $\chi ({\cal J}_k^{k!m}) \ge
cm^{r+1}$ and, a priori:
$$h^0({\cal J}_k^{k!m}) + h^2({\cal
J}_k^{k!m}) > cm^{(n + 1)k -
1}$$ for all $m \ge m_0'$. The Theorem follows now from
assumption (ii).

\bigskip
\bigskip
\noindent
{\bf \S~4~~Surfaces of General Type}

\bigskip
We recall first some well-known results on manifolds of
general type.

\bigskip
\noindent
{\bf Theorem 4.1}~{\it Let $X$ be a minimal surface of
general type then $c_1^2(T^{*}X) > 0, c_2(T^{*}X) > 0$
and $c_1^2(T^{*}X) \ge 3c_2(T^{*}X)$. Moreover, we have
$$5
c_1^2(T^{*}X) - c_2(T^{*}X) + 36
\ge 0,~~{\rm if}~m~{\rm is~even,}$$
$$5
c_1^2(T^{*}X) - c_2(T^{*}X) + 30
\ge 0,~~{\rm if}~m~{\rm is~odd.}$$}

\bigskip
Let $L_0$ be a nef line bundle on a
non-singular surface $X$. A coherent sheaf $E$
over $X$ is said to be semi-stable (resp. stable) with
respect to $L_0$ if $c_1(E)~.~c_1(L_0) \ge 0$ and if, for
any coherent subsheaf $0 \ne {\cal S}$ of $E$, we have:
\begin{eqnarray}\mu_{{\cal S}, L_0}  \stackrel{\rm def}{=}
{c_1({\cal S})~.~c_1(L_0) \over {\rm rank}~{\cal S}} \le
\mu_{E, L_0}  \stackrel{\rm def}{=} {c_1(E)~.~c_1(L_0)
\over {\rm rank}~E}
\end{eqnarray} (resp.
$\mu_{{\cal S}, L_0} <
\mu_{E, L_0}).$

\bigskip
If $X$ is of general type then (see Maruyama [Ma])

\bigskip
\noindent
{\bf Theorem 4.2}~~{\it Let $X$ be a surface of general
type then $\otimes^m T^{*}X, \odot^m T^{*}X$ are semi-stable
with respect to the canonical bundle ${\cal K}_X = \det
T^{*}X$.}

\bigskip
Indeed we have:

\bigskip
\noindent
{\bf Theorem 4.3}~~{\it Let $X$ be a minimal surface of
general type. If
$D$ is a divisor in
$X$ such that
$H^0(X, E_k
\otimes [- D])
\ne 0$ where $E_k = (\odot^{i_1} T^{*}X \otimes ... \otimes
\odot^{i_k} T^{*}X), i_1, ..., i_k$ being positive
integers then
$$c_1(E_k).c_1(D) \le \mu_{E_k} \le {m
\over 2} c_1^2(T^{*}X)$$ with $m = i_1 + 2i_2 + ... + ki_k$}

\medskip
\noindent
{\it Proof.}~~This follows from the calculation of the
Chern number $c_1(E_k)$ in section 2. The computation there
shows that $$\mu_{E_k} \le {m
\over 2} c_1^2(T^{*}X)$$ with equality if and only if
$k = 1$, i.e.,
$$\mu_{\odot^m T^{*}X} = {c_1(\odot^m T^{*}X) \over {\rm
rank}~\odot^m T^{*}X}~.~c_1(T^{*}X) = {{m(m+1)\over 2}
\over m+1} c_1^2(T^{*}X) = {m
\over 2} c_1^2(T^{*}X).$$QED

\bigskip
Note
that in general, if $E$ is a vector bundle of rank
$r$ then
\begin{eqnarray}{\rm rank}~\odot^m E = {(m + r - 1)! \over
(r - 1)! m!}.\end{eqnarray} The Chern number $c_1(\odot^m
E)$ is given by (compare section 2)
\begin{eqnarray}c_1(\odot^m E) = {1 \over r!} {(m + r - 1)!
\over (m - 1)!} c_1(E).\end{eqnarray}
This is done by
induction on the rank of
$E$. If rank $E = 1$ then clearly we have $c_1(\odot^m E) =
m c_1(E)$. If rank $E = 2$ we may formally split the bundle
$E$ as direct sum of line bundles, i.e., we have an exact
sequence:
$$0 \rightarrow L_1 \rightarrow E \rightarrow L_2
\rightarrow 0$$ so that there is a filtration
$$\odot^m E = F_0 \supset F_1 \supset ... \supset F_{m+1} =
0$$ with $F_i/F_{i+1} \cong L_1^{i} \otimes L_2^{m-i}$ and
so, by Whitney's formula:
\begin{eqnarray*}
c_1(\odot^m E) &=& \sum_{i=0}^{m} c_1(F_i/F_{i+1}) \\ &=&
\sum_{i=0}^{m} c_1(L_1^{i} \otimes L_2^{m-i})\\&=&
\sum_{i=0}^{m}  i c_1(L_1) + \sum_{i=0}^{m} (m - i)
c_1(L_2)\\ &=& {m(m+1) \over 2} (c_1(L_1) + c_1(L_2)) \\&=&
{m(m+1) \over 2} c_1(E).
\end{eqnarray*} If rank $E = 3$ then we split the bundle
into a rank 2 bundle $A$ and a line bundle $L$, i.e., we
have an exact sequence:
$$0 \rightarrow F \rightarrow E \rightarrow L
\rightarrow 0$$
so that there is a filtration
$$\odot^m E = F_0 \supset F_1 \supset ... \supset F_{m+1} =
0$$ with $F_i/F_{i+1} \cong \odot^i F \otimes  L^{m-i}$ and
the Chern number is given by:
\begin{eqnarray*}
c_1(\odot^m E) &=& \sum_{i=0}^{m} c_1(F_i/F_{i+1}) \\ &=&
\sum_{i=0}^{m} \{c_1(\odot^i F) + ({\rm rank}~\odot^i F) (m
- i) c_1(L)\}\end{eqnarray*} and by induction the RHS above
is $$
\sum_{i=0}^{m}  {i(i+1) \over 2} c_1(F) + \sum_{i=0}^{m} (i
+ 1) (m - i) c_1(L)$$ hence
\begin{eqnarray*}
c_1(\odot^m E) &=& \sum_{i=0}^{m} {i(i+1) \over 2} c_1(F) +
\sum_{i=0}^{m} (i + 1) (m - i) c_1(L)\\ &=& {1 \over
6} {(m + 2)! \over (m-1)!} c_1(F) +
m \sum_{i=0}^{m} (i + 1) c_1(L) - \sum_{i=0}^{m}
i(i + 1) c_1(L)\\&=& {1 \over
6} {(m + 2)! \over (m-1)!} c_1(F) +
{m(m + 1) \over 2} c_1(L) - {1 \over 3} {(m + 2)! \over (m -
1)!} c_1(L)\\ &=& {1 \over
6} {(m + 2)! \over (m-1)!} c_1(E).
\end{eqnarray*}Ne that we have used the formula:
$$\sum_{i=0}^{m} i(i + 1) = {1 \over 3}{(m+2)! \over
(m-1)!}.$$The general case is proved by induction using the
formula:
$$\sum_{i=0}^{m} i(i + 1)(i + 1) ... (i + k) = {1 \over
k + 2}{(m+k+1)!
\over (m-1)!}.$$QED

\bigskip
Examples of surfaces of general type are provided by
complete intersections in ${\bf P}^n$. A Smooth complete intersection of
type $(d_1, ..., d_{n-r}), 1
\leq r \leq n - 1$, in
${\bf P}^n$ is a smooth variety $X$ of dimension r which is
the transversal intersections of $(n - r)$
hypersurfaces of degree $d_1, ..., d_{n-r}$ respectively.
By the adjunction formula, the canonical bundle of a
complete intersection $X$ of type
$(d_1, ..., d_{n-r})$ is given by the formula:
$${\cal K}_X = {\cal O}_{P^n}(d_1 + ... + d_{n-r} - (n +
1))|_X = {\cal O}_X(d_1 + ... + d_{n-r} - (n + 1)).$$
The normal bundle ${\cal N}_{X|Y}$ of a smooth
hypersurface $X$ in a smooth variety $Y$ is given by
$${\cal N}_{X|Y} = {\cal O}_Y(X)|_X = {\cal O}_X(X).$$
Thus for a hypersurface $X_1$ of degree $d_1$ in ${\bf
P}^n$, the normal bundle
$${\cal N}{X_1|P^n} = {\cal O}_{P1}(d_1)|_{X_1} = {\cal
O}_{X_1}(d_1).$$
Inductively, for a smooth complete intersection $X$ of type
$(d_1, ..., d_{n-r})$ we get
$${\cal N}_{X|P^n} = \oplus _{1\leq i \leq n-r} {\cal
O}_X(d_i).$$
To compute Chern classes of $X$ we apply the Whitney
formula to the exact sequence:
$$0 \longrightarrow TX \longrightarrow TP^n|_X
\longrightarrow {\cal N}_{X|P^n} = \oplus _{1 \leq i \leq
n-r} {\cal O}_X(d_i) \longrightarrow 0.$$ which yields the
following formula for the total Chern classes:
$$c(TX)~.~c({\cal N}_{X|P^n}) = c(T{\bf P}^n|_X).$$
Operating symbolically, we get:
$$1 + c_1(TX) + ... + c_r(TX) = (1 + \theta)^{n+1}/
\prod_{1\leq i
\leq n-r} (1 + d_i\theta)$$
where
$$
\theta^r = \prod_{1 \leq i \leq n-r}  d_i.$$
Expanding formally the RHS above yields:
$$(1 + \theta)^{n+1} = 1 + C_1^{n+1} \theta
+ C_2^{n+1} \theta^2 + ... + C_r^{n+1} \theta^r$$
$$(1 + d_i\theta)^{-1} = 1 - d_i\theta + (d_i\theta)^2 - ...
+ (-1)^r(d_i\theta)^r, 1 \leq i \leq n - r.$$
Define
polynomials $p_q (0 \leq q \leq n - r)$ in $d_1, ...,
d_{n-r}$ by $p_0(d_1, ..., d_{n-r}) = 1$, $$p_q(d_1, ...,
d_{n-r}) =
\sum_{
 1\leq i_1 \leq ...\leq i_q \leq n-r}   d_{i_1} ....
d_{i_q}~1 \leq q \leq n - r.$$
Then for $0 \leq q \leq n - r$, the Chern classes of $X$ are
given by:
\begin{eqnarray*}
c_q(TX) =  \sum_{i=0}^q (-1)^i C_{q-i}^{n+1} ~p_i(d_1, ...,
d_{n-r})~
\theta ^q.
\end{eqnarray*}
For hypersurface ($r = n - 1$) the formulas above reduce to
$$c_q(TX) =  \sum_{i=0}^q (-1)^i C_{q-i}^{n+1} ~ d_i \theta
^q, ~~0 \leq q \leq n - 1. $$ For surfaces of complete
intersections ($r = 2$) and the formulas reduce to:
$$c_1(TX) = ((n + 1) - \sum_{i=1}^{n-2} d_i) \theta, $$
$$c_2(TX) = \{{n(n+1) \over 2} - (n+1) \sum_{i=1}^{n-2} d_i
+
 \sum_{1 \leq i \leq j \leq n-2}     d_id_j\} \theta^2.$$
In particular, if $d_1 = ... = d_{n-2} = d$ then
$$c_1(X) = \{(n + 1) - (n - 2)d\} \theta,$$
$$c_2(X) = \{{n(n+1)\over 2}  - (n+1)(n-2) d +
{(n-1)(n-2) \over 2} d^2\} \theta^2.$$For $n = 3$ then
$$c_1(TX) = (4 - d) \theta,~c_2(TX) = (6  -
4 d + d^2) \theta^2;$$equivalently, for the cotangent
bundle, we have:
$$c_1(T^{*}X) = (d - 4) \theta,~c_2(T^{*}X) = (6  -
4 d + d^2) \theta^2.$$

\bigskip
We shall need a vanishing Theorem (see [G-G]) which is a
consequence of a result of Bogomolov ([B1], [B2]):

\bigskip
\noindent
{\bf Theorem 4.4}~~{\it Let $X$ be a minimal surface of
general type and if the geometric genus $p_g(X) > 0$ then
$$H^2(X, {\cal J}_k^m X) = 0$$ if
$k
\ge 1$ and $m > 2k$.}

\bigskip
Actually, it was asserted in [G-G] that the preceding
Theorem holds without the assumption that $p_g(X) > 0$. At
the momoent I can only get through the proof with this
additional assumption.

\bigskip
Thus the condition of Theorem 3.14 is satisfied for a
minimal surface of general type and we obtained,

\bigskip
\noindent
{\bf Corollary 4.5}~~{\it Let $X$ be a smooth minimal surface
of general type with $p_g(X) > 0$ and if $c_1^2({\cal J}_k^m
X) - c_2({\cal J}_k^m X) > 0$ then there exists
$c > 0$ and
$m_0 > 0$ such that $$\dim H^0(X, {\cal J}_{k}^{k! m} X)
\ge c m^{n(k+1) - 1},
$$if $k \ge 1$ and $m > 1,$ i. e., ${\cal J}_{k}^{k!} X$ is
big.}

\bigskip
\noindent
{\bf Corollary 4.6}~~{\it Let $X$ be a smooth minimal surface
of general type with $p_g(X) > 0$ then ${\cal J}_{k}^{k!} X$
is big for $k \ge 3$.}

\medskip
\noindent
{\it Proof.}~~By the calculation in section
2, $$\Delta({\cal
J}_3^6 X) = c_1^2({\cal
J}_3^6 X) - c_2({\cal J}_3^6 X) = 5175 c_1^2(T^*X) - 175
c_2(T^*X)$$which is clearly $> 0$ in view of Theorem 4.1. The
Corollary now follows from Corollary 4.5.~~QED

\bigskip
If $X$ is a smooth hypersurface the preceding Corollary can
be expressed in terms of the degree:

\bigskip
\noindent
{\bf Corollary 4.6}~~{\it Let $X$ be a non-singular
hyper surface of degree $d$ in ${\bf P}^3$. Then ${\cal
J}_2^2 X$ is big if $d \ge 9$ and ${\cal J}_3^6
X$ is big if $d \ge 5$.}

\medskip
\noindent
{\it Proof.}~~By the calculation in sectin 2, we have:
$$c_{1}^{2}({\cal J}_{2}^{2}X) -
c_{2}({\cal J}_{2}^{2}X) = 11
c_{1}^{2}(T^{*}X) - 5
c_{2}(T^{*}X)$$
and as noted before, for a
 smooth hypersurface $X$ in ${\bf P}^{3}$ of degree $d$, the
Chern numbers are given by
$$c_{1}(T^{*}X) = d - 4,~~c_{2}(T^{*}X) = d^{2} - 4d +
6,$$ we conclude that:
$$c_{1}^{2}({\cal J}_{2}^{2}X) - c_{2}({\cal
J}_{2}^{2}X) = 11d^2(d - 4)^{2} - 5d^2(d^{2}
- 4d + 6) > 0$$ if $d \ge 9$. Computing
similarly we conclude that
$c_{1}^{2}({\cal J}_{3}^{6}X) -
c_{2}({\cal J}_{3}^{6}X) > 0$ if $d \ge
5$. Moreover, by Noether;s Theorem (i. e.,
Riemann-Roch):
$$1 - q(X) + p_g(X) = {1 \over 12} (c_1^2(T^{*} X) +
c_2(T^{*} X))$$ implies that $p_g(X) > 0$ because the
irregularity $q(X) = 0$.~~QED

\bigskip
We need one last observation to deal with the fact that ${\cal
J}_k^m X$ is not semi-stable as can be seen from the
calculation in section 2. In fact each of the factors
$\odot^{i_1} T^{*} X \otimes ... \otimes \odot^{i_k}
T^{*}X, i_1 + 2i_2 + ... + ki_k = m$ which is a subsheaf of
${\cal J}_k^m X$ (note that not all of them are) is a
destabling subsheaf). However, we also observe that each of
these sheaves is semi-stable (by Theorem 4.3). Moreover the
ratio:
$c_1(\odot^{i_1}
T^{*}X \otimes ... \otimes \odot^{i_k} T^{*}X) / {\rm
rank}~(\odot^{i_1} T^{*}X \otimes ... \otimes \odot^{i_k}
T^{*}X) \le m/2$ thus we have:

\bigskip
\noindent
{\bf Theorem 4.7}~~{\it Let $X$ be a complex surface such
that $Pic X \cong {\bf Z}$. If
$H^0(X, {\cal J}_k^m X \otimes [-D])
\ne
\{0\}$ where $D$ is a divisor in
$X$ then $$c_1([D]) ~.~c_1(T^{*}X) \le {m \over 2}
c_1^2(T^{*}X).$$ }

\medskip
\noindent
{\it Proof.}~~This follows from the filtration:
$${\cal J}_{k-1}^{m} X = {\cal F}_{k}^{0}
\subset {\cal F}_{k}^{1} \subset ... \subset
{\cal F}_{k}^{[m/k]} = {\cal J}_{k}^{m} X$$
$($where $[m/k]$ is the greatest integer smaller
than or equal to $m/k)$ such that
$${\cal F}_{k}^{i}/{\cal F}_{k}^{i-1} \cong
{\cal J}_{k-1}^{m-ki} X  \otimes (\odot^{i} T^{*}
X).$$From the exact sequence
$$0 \rightarrow {\cal F}_{k}^{[m/k] - 1} \rightarrow {\cal
J}_{k}^{m} X \rightarrow {\cal J}_{k-1}^{m-k[m/k]} X \otimes
(\odot^{[m/k]} T^{*} X) \rightarrow 0$$we obtain an exact
sequence:
\begin{eqnarray*}\lefteqn{0 \rightarrow H^0(X, {\cal
F}_{k}^{[m/k] - 1}
\otimes [-D])
\rightarrow H^0(X, {\cal J}_{k}^{m} X \otimes [-D])
\rightarrow}\\
&\rightarrow &H^0(X, {\cal J}_{k-1}^{m-k[m/k]} X \otimes
(\odot^{[m/k]} T^{*} X) \otimes [-D])
\end{eqnarray*}which shows that if $H^0(X, {\cal J}_{k}^{m} X
\otimes [-D]) \ne \{0\}$ then either $H^0(X, {\cal
F}_{k}^{[m/k] - 1}
\otimes [-D]) \ne \{0\}$ or $$H^0(X, {\cal
J}_{k-1}^{m-k[m/k]} X
\otimes (\odot^{[m/k]} T^{*} X) \otimes [-D]) \ne \{0\}$$and
eventually this means that
$${\rm either}~H^0(X, \odot^{i_1} T^{*}X \otimes ...
\otimes \odot^{i_k} T^{*}X
\otimes [-D]) \ne \{0\}$$
for at least one of the factors $\odot^{i_1} T^{*}X \otimes
... \otimes \odot^{i_k} T^{*}X, i_1 + 2 i_2 + ... + ki_k = m$.
With this the Theorem follows from Theorem 4.3.~~QED

\bigskip
From the computation in section 2 we see that
$$\mu({\cal J}_k^{k!} X) = {c_1({\cal J}_k^{k!} X) \over
{\rm rank}~ {\cal J}_k^{k!} X} < {m \over 2} c_1(T^{*}
X).$$Thus the estimate is weaker than one would get if it
were stable, however this is the best that one can do and
this weaker estimate is sufficient for our purpose.

\bigskip
We assume from now on that $(i) ~X$ is a minimal
surface of general type, (ii) $p_g(X) > 0$ and $(iii) ~Pic(X)
\cong {\bf Z}$. Then ${\cal J}_k^{k!} X$ is big for $k \ge
3$. This implies that ${\cal J}_k^{k!} X \otimes [- D]$ is
big for any effective ample divisor $D$ in $X$. Schwarz
Lemma implies that the lifting of any holomorphic curves in
${\bf P}(J^kX)$ is contained in the zero set of a
non-trivial section of ${\cal O}(k!m) \otimes p^{*}[-D])$. We
proceed to consider the subvarieties of
${\bf P}(J^kX), k \ge 2$.  Let
$Y_1$ be an irreducible effective horizontal (i. e., not of
the form $p^{*} D$ for some effective divisor $D$ in $X$)
divisor in
${\bf P}(J^kX)$ then:
$$[Y_1] = {\cal O}_{{\bf
P}(J^k X)}(m_1) \otimes p^{*}[- D_1]$$ where $D_1$ is a
divisor in $X$ and $m_1 \in {\bf N}$, we may
assume that $m_1$ is divisible by $k!$ by replacing $Y_1$
with $k! Y_1$ (so it is non-reduced but set theoretically it
has only one irreducible component). Thus we may write
$m_1 = k! \alpha_1$. For simplicity of notations we shall
write
${\cal O}(j)$ instead of ${\cal O}_{{\bf P}(J^k X)}(j)$.
Since dim
${\bf P}(J^k X) = 2(k + 1) - 1 = 2k + 1$ we get from Theorem
4.7 and Theorem 3.13:
\begin{eqnarray*}
\lefteqn {c_1^{2k}({\cal O}(k!)|_{Y_{1}})} \\&=&
c_1^{2k+1}({\cal O}(k!)).(c_1({\cal O}(k! \alpha_1) -
p^{*} c_1([D_1])\\&=& (k!)^2 \{\alpha_1 c_1^{2k+
1}({\cal O}(k!)) - c_1^{2k}({\cal O}(k!))~.~p^{*}
c_1([D_1])\}\\&=& (k!)^2 \{\alpha_1 (c_1^2({\cal J}_k^{k!}
X) - c_2({\cal J}_k^{k!} X)) - c_1({\cal
J}_k^{k!} X)~.~c_1([D_1])\}
\\&=& (k!)^2 \{\alpha_1 (c_1^2({\cal
J}_k^{k!} X) - c_2({\cal
J}_k^{k!} X)) - a(k, k!) c_1(T^{*} X)~.~c_1([D_1])\}\\&\ge&
(k!)^2 \{\alpha_1 (a(k, k!)^2 c_1^2(T^{*} X) - c_2({\cal
J}_k^{k!} X)) - {\alpha_1 k! a(k, k!) \over 2}
c_1^2(T^{*}X)\}\\&=&  (k!)^2 \alpha_1 \{a(k, k!) (a(k, k!)
- {k!
 \over 2}) c_1^2(T^{*} X) - c_2({\cal
J}_k^{k!} X)\}.
\end{eqnarray*}This means that ${\cal O}(k!)|_{Y_1}$ is
again big if $$a(k, k!) (a(k, k!) - {k!
 \over 2}) c_1^2(T^{*} X) - c_2({\cal
J}_k^{k!} X) > 0.$$ For example if $k = 3, a(k, k!) = 101$ by
the calculation in section 2; the preceding inequality
yields:
\begin{eqnarray*}
c_1^{2k}({\cal O}(k!)|_{Y_{1}}) &\ge&
  (k!)^2 \alpha_1 \{101 (101 - 3) c_1^2(T^{*} X) - c_2({\cal
J}_k^{k!} X)\}\\&=& (k!)^2 \alpha_1 \{(9292 - 5026)
c_1^2(T^{*} X) - c_2(T^{*} X)\}\\&=& (k!)^2 \alpha_1 (4266
c_1^2(T^{*} X) - 175 c_2(T^{*} X)) \\&>& 0.
\end{eqnarray*}This means that ${\cal O}(k!)|_{Y_1}$ is
again big. The Schwarz Lemma in the appendix again implies
that the lifting of any holomorphic curves in
${\bf P}(J^kX)$ is contained in the zero set of a
non-trivial section of ${\cal O}(k!m)|_{Y_1} \otimes
p|_{Y_1}^{*}[-D])$.

\bigskip
Next we consider divisor $Y_2$ in $Y_1$ which is  of the
form:
$$[Y_2] = ({\cal O}(m_2) \otimes p^{*}[- D_2])|_{Y_1}$$
where $D_2$ is a divisor in $X$ and $m_2 \in {\bf N}$ which
we may assume to be divisible by $k!$, i.e., $m_2 = \alpha_2
k!$. We remark that for the investigation of degeneration of
liftings of a holomorphic curve
$f : {\bf C} \rightarrow X$ these are the only type of
subvarieties that we have to deal with. We have:
\begin{eqnarray*}
\lefteqn {c_1^{2k-1}({\cal O}(k!)|_{Y_{2}})} \\&=&
c_1^{2k-1}({\cal O}(k!)).(c_1({\cal O}(k! \alpha_1) -
p^{*} c_1([D_1]).(c_1({\cal O}(k! \alpha_2) -
p^{*} c_1([D_2])\\&=& \alpha_1 \alpha_2 c_1^{2k+1}({\cal
O}(k!)) - \alpha_1 c_1^{2k}({\cal O}(k!))~.~p^{*}
c_1([D_2]) \\&&- \alpha_2 c_1^{2k}({\cal
O}(k!))~.~p^{*} c_1([D_1]) +
p^{*}c_1([D_1])~.~p^{*}c_1([D_2])\\&\ge& (k!)^2 \{\alpha_1
\alpha_2 (c_1^2({\cal
J}_k^{k!} X) - c_2({\cal
J}_k^{k!} X)) - \alpha_1
c_1({\cal J}_k^{k!}X)~.~c_1([D_1]) \\&&-
\alpha_2 \alpha_1
c_1({\cal J}_k^{k!}X)~.~ c_1([D_1]) +
c_1([D_1])~.~c_1([D_2])\}
\\&=& (k!)^2 \{\alpha_1
\alpha_2 (c_1^2({\cal
J}_k^{k!} X) - c_2({\cal
J}_k^{k!} X)) - \alpha_1 a(k, k!) c_1(T^{*}
X)~.~c_1([D_2])\\&&- \alpha_2 a(k, k!) c_1(T^{*}
X)~.~c_1([D_1]) +
c_1([D_1])~.~c_1([D_2])\}\\&\ge&
(k!)^2 \{\alpha_1
\alpha_2 (c_1^2({\cal
J}_k^{k!} X) - c_2({\cal
J}_k^{k!} X)) - 2{\alpha_1  \alpha_2 k! a(k, k!) \over 2}
c_1^2(T^{*} X)\}\\&=&  (k!)^2 \alpha_1  \alpha_2 \{a(k, k!)
(a(k, k!) - k!) c_1^2(T^{*} X) - c_2({\cal
J}_k^{k!} X)\}.
\end{eqnarray*} Proceeding inductively, we get a sequence of
subvarieties $Y_1 \supset Y_2 \supset \cdots \supset
Y_{2k}$ where each $Y_i$ is of codimension $i$ and of the
form $$[Y_{i+1}] = ({\cal O}(m_{i+1}) \otimes p^{*}[-
D_{i+1}])|_{Y_{i}}.$$A similar calculation shows that:
\begin{eqnarray*}
\lefteqn{c_1^{2k-i+1}({\cal O}(k!)|_{Y_{i}})} \\&\ge&
(k!)^2 \alpha_1  \cdots \alpha_i \{a(k, k!)
(a(k, k!) - k! {i\over 2} ) c_1^2(T^{*} X) - c_2({\cal
J}_k^{k!} X)\},
\end{eqnarray*}$i = 1, ..., 2k$. For $k = 3$ we have:
\begin{eqnarray*}
\lefteqn{a(3, 3!)
(a(3, 3!) - 3! {i\over 2} ) c_1^2(T^{*} X) - c_2({\cal
J}_k^{k!} X)}\\&\ge& 101 (101 - (3!)3) c_1^2(T^{*} X) -
c_2({\cal J}_k^{k!} X)\\&=& 3357 c_1^2(T^{*} X) - 175
c_2(T^{*}X) \\&>& 0
\end{eqnarray*}for all $i = 1, ..., 2k = 6$. Thus we arrive
at the following Theorem:

\bigskip
\noindent
{\bf Theorem 4.8}~~{\it Let $X$ be a smooth minimal surface
of general type such that $(i) ~Pic X \cong {\bf Z}$ and
$(ii)~p_g(X) > 0$. Then $X$ is hyperbolic. Consequently, a
generic smooth hypersurface in ${\bf P}^3$ of degree $d \ge 5$ is
hyperbolic.}

\bigskip
In [D-E] certain types of 2-jet differentials ${\cal A}$
were used and the authors establised that $c_1^2({\cal A}) -
c_2({\cal A}) = 13 c_1^2(T^{*}X) - 9 c_2(T^{*}X)$ on any
hypersurface $X$ of degree $\ge 42$. This is weaker than
what we have, namely $c_1^2({\cal J}_2^2X) - c_2({\cal
J}_2^2X) = 11 c_1^2(T^{*}X) - 5 c_2(T^{*}X)$. Actually I
have some trouble using this stronger estimate to get
hyperbolicity due to non-semistability (recall that neither
${\cal J}_k^m$ nor $T^{*}_k X$ is semi-stable) so
that $c_1^2({\cal J}_2^2X) - c_2({\cal J}_2^2X) = 11
c_1^2(T^{*}X) - 5 c_2(T^{*}X)$ is not big enough to reach
hyperbolicity. It appears that
${\cal A}$ is not semi-stable either.

\bigskip
\bigskip
\noindent
{\bf Appendix A:~The Lemma of Logarithmic Derivatives}

\bigskip
One of the main tool in Nevanlinna Theory is the
classical Lemma of Logarithmic Derivatives
(abbrev. LLD). For example, LLD implies that,
even though there is no pointwise estimate
between (the absolute value of) a holomorphic
function and (the absolute value of) its
derivatives, such estimates do exist in the sense
of integral averages (i.e, their characteristic
functions bound each other). The purpose of this appendix
is to extend the classical LLD to all jet differentials of
logarithmic type and in
particular all regular jet differentials. The
proof is based on the very simple observation
that (the absolute value of) any jet differential
of logarithmic type is bounded by (the absolute
value of) those of the classical type (hence the
classical LLD applies).

\bigskip
\noindent
{\bf Theorem A1}~(Lemma of Logarithmic
Derivatives)~~{\it Let
$X$ be a projective variety and let $(i)~D$ be an
effective divisor with simple normal crossings or
$(ii)~D$ is the trivial divisor in
$X~($i.e. the support of $D$ is empty or
equivalently, the line bundle associate to
$D$ is trivial$)$. Let $f : {\bf C}
\rightarrow X$
be a holomorphic map
and $\omega \in H^0(X, {\cal J}_k^mX(\log
D))~($resp. $H^0(X, {\cal J}_k^mX)$ in case
$(ii))$ a jet differential such that $\omega \circ j^kf$ is
not identically zero, then
$$T_{\omega \circ j^{k}f}(r) = \int_{0}^{2\pi} \log^+
|\omega(j^kf(re^{\sqrt{-1}\theta}))| \frac{d\theta}{2\pi}
 ~.\leq.~   O(\log T_f(\omega_X; r))  + O(\log r).$$
Here $\omega_{X}$ can be taken to be $c_{1}({\cal L})$
of any ample line bundle
${\cal L}$ on $X$.}

\medskip
\noindent
{\it Proof}.~We claim that there exists a finite
number of rational functions
${t_1, ..., t_q}$ on $X$ such that:
\begin{quote}
$ (\dag)	$  {\it the logarithmic jet differentials
$\{(d^{(j)}t_i/t_i)^{m/j} ~|~ 1 \leq i
\leq q, 1 \leq j \leq k\}$ span the fibers of
	${\cal J}_k^mX(\log D)$ over every point of $X$.}
\end{quote}
Without loss of generality we may assume that $D$
is ample; otherwise we may replace $D$ be $D +
D^{'}$ so that $D + D^{'}$ is ample. Observe that
if $s$ is a function holomorphic on a
neighborhood $U$ such that $[s = 0] = D \cap U$
then $[s^{\tau} = 0] = \tau D \cap U$ where
$\tau$
is a rational number. Thus
$d^{(j)}(\log s^{\tau}) = \tau d^{(j)}(\log s)$ is still a
jet differential with logarithmic
singularity along $D \cap U$ so the mutiplicity causes
no problem. This means that we may assume without loss of
generality that $D$ is very ample by replacing
$D$ with $\tau D$ for some $\tau$ so that $\tau D$ is very
ample.

\medskip
	Let $u \in H^0(X, [D])$ be a section such that $D = [u =
0]$. At a point
$x \in D$ choose a section $v_1 \in H^0(X, [D])$ so
that $E_1 = [v_1 = 0]$
 is smooth, $D + E_1$ is of simple normal crossings
and $v_1$ is non-vanishing at $x$ (this is possible because
the line bundle $[D]$ is very ample).  The rational
function $t_1 = u_1/v_1$ is regular on the affine open
neighborhood $X \setminus E_1$ of $x$ and $(X \setminus
E_1) \cap [t_1 = 0] = (X \setminus E_1) \cap D$.
 Choose rational functions $t_2 = u_2/v_2, ..., t_n  =
u_n/v_n$ where $u_i$ and $v_i$
 are sections of a very ample bundle ${\cal L}$ so
that $t_2, ..., t_n$ are regular at
$x$, the divisors $D_i = [u_i = 0], E_i = [v_i = 0]$
are smooth and $D + D_2 + ... + D_n + E_1 + ... + E_n$ is of
simple normal crossings. Moreover, since the bundles
involved are very ample the sections can be chosen so that
$dt_1 \wedge ... \wedge dt_n$ is non-vanishing at $x$; the
complete system of
 sections provides an embedding, hence at each point there
are $n + 1$ sections with the property that $n$ of the the
quotients of these $n + 1$ sections form a local coordinate
system on some open neighborhood $U_x$ of $x$. This implies
that $ (\dag)	$ is satisfied over $U_x$.
Since $D$ is compact it is covered by a
finite number of such open neighborhoods, say $U_1,
..., U_p$ and a finite number of rational functions
(constructed as above for each $U_i$) on X so that $ (\dag)
$ is satisfied on
$\cup_{1 \leq i \leq p} U_i$. Moreover,
there exists relatively compact open subsets $U_i'$ of $U_i$
($1 \leq i \leq p$) such that $\cup_{1 \leq i \leq p} U_i'$
still covers $D$.

\medskip
	Next we consider a point $x$ in the compact set
$X \setminus \cap_{1 \leq  i \leq p} U_i'$. Repeating the
procedure as above we can find rational functions $s_1 =
a_1/b_1,
 ..., s_n = a_n/b_n$ where $a_i$ and $b_i$ are sections of
some very ample line ${\cal L}$ bundle so that $s_1, ...,
s_n$ from a holomorphic local coordinate on some open
 neighborhood $V_x$ of $x$. Thus $ (\dag)	$ is satisfied on
$V_x$ by the rational functions $s_1, ..., s_n$. Note that
we must also choose these sections  so that the divisor
 $H = [s_1~ ...~ s_n = 0]$ together with those
divisors (finite in number), which had been already
constructed above, is still a divisor with simple normal
crossings (this is possible by the very ampleness of the
line bundle ${\cal L}$.) Since $X \setminus \cap_{1 \leq i
\leq p} U_i'$ is compact, it is cover by a finite of such
coordinate neighborhoods. The coordinates are rational
functions and finite in number and by construction it is
clear that the condition $(\dag)$ is satisfied on
$X \setminus \cap_{1 \leq i \leq p} U_i'$. Since $\cup_{1
\leq i \leq p} U_i$ together with
$X \setminus \cap_{1 \leq i \leq p} U_i'$ covers $X$, the
condition $(\dag)$ is satisfied on
$X$.

\medskip
	If $D$ is the trivial divisor, then it is enough to
use only the second part of the construction above and again
$(\dag)$ is verified with ${\cal J}_k^mX(\log D) = {\cal
J}_k^mX$

\medskip
	To obtain the estimate of the Theorem observe that
the function,
$$ \rho : J^kX(- \log D) \rightarrow [0, \infty] $$
defined by
\begin{eqnarray*}
\rho(\xi) = \sum_{i=1}^{q}
\sum_{j=1}^{k} |(d^{(j)}t_i/t_i)^{m/j}(\xi)|^2,~
\xi \in J^kX(- \log D)
\end{eqnarray*}
($\{t_i\}$ is the family of rational functions satisfying
the condition $(\dag)$) is continuous in the extended
sense; it is continuous, in the ususal sense,
 outside the fibers over the divisor $E$ (the sum of the
divisors associated to the rational functions $\{t_i\}$;
note that $E$ contains $D$). Over the fiber of each point
$x \in X - E$, $|(d^{(j)}t_i/t_i)^{m/j}(\xi)|^2$ is
finite for
 $\xi \in J^kX(- \log D)_x$, thus $\rho$ is not
identically infinite. Moreover, since
$$\{(d^{(j)}t_i/t_i)^{m/j} ~|~ 1 \leq i \leq q, 1 \leq
j \leq k\} $$ span the fiber of ${\cal J}_k^mX(\log D)$ over
every point of $X$, $\rho$ is strictly positive ($+ \infty$
allowed) outside the zero-section of $J_kX(-\log D)$. The
quotient
$$|\omega|^2/\rho : J^kX(-\log D) \rightarrow [0, \infty]$$
does not take on the extended value $\infty$ when restricted
to $J^kX(-\log D) \setminus \{zero-section\}$ becuase, as we
have just observed, $\rho$
 is non-vanishing (even though it blows up along the fibers
over $E$ so that the reciprocal $1/\rho$ is zero there) and
the singularity of $|\omega|$ is no worst then that of
$\rho$ becuase the singularity of $\omega$ ocuurs only
along $D$ (which is contained in $E$) and is of log type.
Thus the restriction to $J_kX(-\log D) \setminus
\{zero-section\}$,
$$|\omega|^2/\rho :
J^kX(-\log D) \setminus \{zero-section\} \rightarrow [0,
\infty)$$ is a continuous non-negative funtion.
Moreover, since
$|\omega|$ and
$\rho$ have the same homogenity:
$$ |\omega(\lambda.\xi)|^2 = |\lambda|^{2m}
|\omega(\lambda.\xi)|^2~~ and~~
\rho(\lambda.\xi) = |\lambda|^{2m}\rho(\xi)$$
for all $\lambda \in {\bf C}^*$ and $\xi \in J^kX(-\log
D)$ we see that
$|\omega|^2/\rho$ descends to a well-defined function on
${\bf P}(E_{k,D}) = (J^kX(-\log D)
\setminus \{zero-section\}) / {\bf C}^*$, i.e.,
$$ |\omega|^2/\rho : {\bf P}(E_{k,D}) \rightarrow
[0, \infty)$$ is a well-defined continuous function and so,
by compactness, there exists a constant $c$ with the
property that
\begin{eqnarray*}
 |\omega|^2 \leq c\rho.
\end{eqnarray*}
This implies that
\begin{eqnarray*}T_{\omega \circ j^kf}(r) &=&
\int_{0}^{2\pi}
\log^+
|\omega(j^kf(re^{\sqrt{-1}\theta}))| \frac{d\theta}{2\pi}
\\&\leq&
 \int_{0}^{2\pi} \log^+
|\rho(j^kf(re^{\sqrt{-1}\theta}))| \frac{d\theta}{2\pi}
  + O(1).\end{eqnarray*}
Since $t_i$ is a rational function on $X$ the function
$$ (d^{(j)}t_i/t_i)^{m/j}(j^kf) = ((t_i \circ
f)^{(j)}/t_i \circ f)^{m/j} $$ ($m$ is divisible by $k!$) is
meromorphic on ${\bf C}$ and so, by the definition
 of $\rho$,
$$ \log^+ |\rho(j^kf)| \leq O(\max_{1 \leq i \leq q,1 \leq j
\leq k}
\log^+ | (t_i \circ f)^{(j)}/t_i \circ f |) + O(1). $$
Now by the classical lemma of logarithmic derivatives
for meromorphic functions,
$$ \int_{0}^{2\pi} \log^+ | (t_i \circ f)^{(j)}/t_i \circ f
|) \frac{d\theta}{2\pi}
  ~.\leq.~  O(\log rT_{t_i \circ f}(r)). $$
Since $t_i$ is a rational function,
$$ \log T_{t_i \circ f}(r) \leq O(\log T_f(\omega_X; r)) +
O(1) $$ and we arrive at the estimate
\begin{eqnarray*}
&&\int_{0}^{2\pi} \log^+ |\rho(j^kf(re^{\sqrt{-1}\theta})|)
\frac{\theta}{2\pi}\\
&&~\leq O(\int_{0}^{2\pi} \log^+ |(t_i \circ
f)^{(j)}/t_i \circ f|
\frac{d\theta}{2\pi}) + O(1)\\
&&.\leq.~O(\log T_f(r)) + O(\log
r).
\end{eqnarray*}
This implies that
$$ T_{\omega \circ j^kf}(r) ~.\leq.~ O(\log T_f(r)) + O(\log
r) $$ as claimed.~~ QED

\bigskip
We obtain, as imediate consequence, the following
Schwarz's type Lemma for logarithmic jet
differentials.

\bigskip
\noindent
{\bf Corollary A2}~~{\it Let $X$ be a projective
variety and $D$ be an effective divisor
$($possibly the trivial divisor$)$ with simple
normal crossings. Let
$f : {\bf C}
\rightarrow X
\setminus D$ be a holomorphic map. Then
$$\omega(j^{k}f) \equiv 0~~{\rm for~all}~~\omega
\in H^{0}(X, {\cal J}_{k}^{m}X(\log D) \otimes
[-H])$$ where $H$ is a generic hyperplane section
$($and hence any hyperplane section$)$.}

\medskip
\noindent
{\it Proof.}~If $f$ is constant then the
Corollary holds trivially. So we may assume that
$f$ is non-constant and suppose that $\omega \circ
j^{k}f \not\equiv 0$. Moreover, since $F$ is
non-constant, we may assume without loss of
generality that $\log r = o(T_{f}(H; r))$ by
replacing $f$ with $f \circ \phi$ where $\phi$ is
a transcendental function on ${\bf C}$. By
Theorem 4.1, we have
$$\int_{0}^{2 \pi} \log^{+} |\omega \circ
j^{k}f|~{d\theta \over 2 \pi} = T_{\omega \circ
j^{k}f}(r) ~.\leq.~O(\log rT_{f}(H; r)).$$
On the other hand, since $\omega$ vanishes on $H$
and $H$ is generic (see (1) or (2) in section 1),
we obtain via Jensen's Formula:
\begin{eqnarray*}T_{f}(H; r) &\leq& N_{f}(H; r) +
O(\log rT_{f}(H; r)) \\&=&
\int_{0}^{2
\pi}
\log |\omega
\circ j^{k}f|~{d\theta \over 2 \pi} + O(\log
rT_{f}(H; r))\end{eqnarray*} which, together with
the preceding estimate, imply that:
$$T_{f}(H; r)~.\leq.~O(\log rT_{f}(H; r)).$$
This is impossible hence we must have $\omega
\circ j^{k}f \equiv 0$. If $H_{1} = [s_{1} = 0]$
is any hyperplane section then it is linearly
equivalent to a generic hyperplane section $H =
[s = 0]$. If $\omega$ vanishes along $H^{'}$ then
$(s/s_{1}) \omega$ vanishes along $H$. The
preceding discussion implies that $(s/s_{1})
\omega(j^{k}f) \equiv 0$. This implies that
actually $\omega(j^{k}f) \equiv 0$ as we may
choose a generic section $H$ so that the image of
$f$ is not entirely contained in $H$.~~QED

\bigskip
Actually the proof of Theorem A.1 gives a little
more. In fact the same proof yields:

\bigskip
\noindent
{\bf Theorem A3}~~{\it Let $\rho_{k}$ be a
pseudo singular jet metric on $J^{k}X(-\log D)$
with the property that there exists a constant $c
> 0$ such that $\rho_{k} \leq \rho$ where $\rho$
is the singular jet metric on $J^{k}X(-\log D)$
defined by the family of rational functions
$(\dagger)~($see
$(28))$. Then $$T_{j^{k}f}(\rho_{k};
r) =
\int_{0}^{2\pi}
\log^+ |\rho_{k}(j^kf(re^{\sqrt{-1}\theta}))|
\frac{d\theta}{2\pi}
 ~.\leq.~   O(\log rT_f(\omega_X;
r))$$ for any K\"{a}hler mertic
$\omega_{X}$ on
$X$. In particular, if $\rho_{k}$ is a
non-singular pseudo metric on $J^{k}X$ then the
preceding estimate holds.}

\bigskip
The Schwarz Lemma can be further extended as follows.

\bigskip
\noindent
{\bf Theorem A4}~~{\it Let $Y \subset {\bf P}(J^kX)$ be a
subvariety and suppose that there exists a non-trivial section
$\sigma \in H^{0}(Y, {\cal O}_{{\bf P}(J^kX)})(m)|_{Y} \otimes
p|_{Y}^{*}[-D]$ where $D$ is a generic ample divisor in $X$ and
$p : {\bf P}(J^kX) \rightarrow X$ is the projection map. If
the image of the lifting
$[j^k f] : {\bf C} \rightarrow {\bf P}(J^kX)$ of a holomorphic
curve $f : {\bf C} \rightarrow X$ is contained in $Y$ then
$\sigma([j^kf]) \equiv 0$.}

\bigskip
\bigskip
\noindent
{\bf Appendix B}

\bigskip
Let $S_{n}$ be the symmetric group on $n$
elements then the order of $S_{n}$ is $n!$.

\bigskip
\noindent
{\bf Definition B1}~~A maximal set of mutually
conjugate elements of $S_n$ is said to be a
class of $S_n$.

\bigskip
\noindent
{\bf Definition B2}~~A partition of a natural
number
$n$ is a set of positve integers $k_1, ...,
k_q$ such that $n = k_1 + ... + k_q$.

\bigskip
A
partition can be expressed as $$n =
\sum_{i=1}^{n} i a_{i}$$ where the
integers $a_i$ are non-negative.

\bigskip
\noindent
{\bf Theorem B3}~~{\it The number, denoted
$p(n)$, of classes of
$S_n$ is equal to the number of partitions of
$n$ and also to the number of inequivalent
irreducible representations of $S_n$. The number
$p(n)$ is asymptotically approximated by the
formula of Hardy-Ramanujan
$$p(n) \sim {e^{\pi \sqrt{2n/3}} \over 4n
\sqrt{3}}.$$}

\bigskip
The alternating subgroup $A_n$ (i.e. the even
permutations) is the commutator subgroup of $S_n$
and is obviously of index 2. Thus there are two
1-dimensional representations of $S_n$: the
trivial representation and the representation
$\Gamma_{\sigma}$ defined by $P \mapsto
\sigma(P)$ where $\sigma$ is the signature of a
permutation $P$ (i.e. $P \mapsto \pm 1$
depending on whether $P$ is even or odd (i.e.,
can be expressed as an even or odd number of
transpositions: interchanging two of the $n$
elements).

\bigskip
\noindent
{\bf Lemma B4}~~{\it Let $X = X(n)$ be a set of
$n$ elements and let $Y_1, ..., Y_k$ be $k$ not
necessarily distinct subsets of $X$. For any
subset $J$ of the index set $\{1, ..., k\}$,
denote by
$$n(J) = \# \cap_{j \in J} Y_j;$$ and for $0 \le
i \le k$, denote by
$$n_0 = n,~n_i = \sum_{\# J = i} n(J),~1 \le i
\le k.$$ Then the number of elements not in any
of the subsets $Y_i, i = 1, ..., k$ is given by
the formula
$$\# (X \setminus (\cup_{i=1}^{k} Y_i)) = n -
n_1 + n_2 - ... + (-1)^{k} n_k = \sum_{i=0}^{k}
(-1)^{i} n_i.$$}

\medskip
\noindent
{\it Proof.}~~If $k = 2$ the formula can be
expressed as usual:
$$\# (X \setminus (Y_1 \cup Y_2)) = \# X - \#
Y_1 - \# Y_2 + \# (Y_1 \cap Y_2).$$ One way to
prove the Lemma is by induction on $k$.
Alternatively one can also argue as follows:~~QED

\bigskip
An element $P$ of ${\cal S}_n$ is said to be a
derangement if $P(i) \ne i$ for $i = 1, ...,
n$. The number of derangements is denoted by
$d_n$. Then

\bigskip
\noindent
{\bf Corollary B5}~~{\it The number of
derangements is given by the formula
$$d_n = n! \sum_{i=0}^{n} {(-1)^{i} \over
i!}.$$ In particular we see that asymptotically
$d_n \sim e^{-1}$.}

\medskip
\noindent
{\it Proof.}~~Apply Lemma 4 with $X = {\cal
S}_n$ and $Y_i = \{P \in {\cal S}_{n}~|~P(i) =
i\}, i = 1, ..., n.$

\medskip
Alternatively, the formula can be obtained
by considering the power series:
$$e^{x} \sum_{i=0}^{\infty} d_i {x^i \over
i!} = \sum_{i=0}^{\infty} (\sum_{j=0}^{i} {i!
\over j! (i - j)!}~ d_{i-j}) {x^i \over i!} =
\sum_{i=0}^{\infty} x^{i} = {1 \over 1 -
x}.$$Thus we have
$$\sum_{i=0}^{\infty} d_i {x^i \over
i!} = e^{-x} (1 - x)^{-1}$$which yields the
formual of the corollary.~~QED

\bigskip
The formula of the Corollary can also be
obtained via the recursive formula:
$$d_n - n d_{n-1} + (-1)^{n}.$$

\bigskip
\noindent
{\bf Corollary B6}~~{\it The number of
surjections from a set $A$ of $n$ elements to a
set $B$ of $k$ elements is given by the formula
$$\sum_{i=0}^{k} (-1)^{i} {k! \over i! (k -
i)!} (k - i)^{n}.$$}

\medskip
\noindent
{\it Proof.}~~Apply Lemma 4 to the set $X$ of
all maps from $A$ to $B$ and $Y_i, = 1, ..., k$
be the subset consisting of those maps such that
$i$ is not in the image.~~QED

\bigskip
Note that Corollary 6 implies trivially that
\begin{eqnarray*}
\sum_{i=0}^{k} (-1)^{i} {k! \over i! (k -
i)!} (k - i)^{n} =
\left\{\begin{array}{ll} n! &\mbox{{\rm if}
$k = n$},\\
0 &\mbox{{\rm if} $k > n$.}
\end{array}
\right.
\end{eqnarray*}
There is a more general formula which can be
proved in a similar fashion:
\begin{eqnarray*}
\sum_{i=0}^{n} (-1)^{i} {n! \over i! (n - i)!}
{(m + n - i)! \over (m + n - k)! (k - i)!} =
\left\{\begin{array}{ll} m! / k!
(m - k)! &\mbox{{\rm if} $m \ge k$},\\
0 &\mbox{{\rm if} $m < k$.}
\end{array}
\right.
\end{eqnarray*}

\bigskip
\noindent
{\bf Theorem B7}~~{\it The number of non-negative
integer solutions of the equation
$$x_1 + ... + x_k = n$$ is $(n + k - 1)!/(k -
1)! n!.$ On the other hand the number of
positive integer solutions is $(n - 1)!/(k - 1)!
(n - k)!$.}

\medskip
\noindent
{\it Proof.}~~So we have to find the number of
ways to put $n$ {\it black} (otherwise
identical) balls in
$k$ slots. If we insert white balls in between
the slots we end up with a total of $n + k - 1$
balls $k - 1$ of them white. This is the same as
choosing $k - 1$ balls from a total of $n + k -
1$ balls and the first assertion follows.

\medskip
The second assertion follows from the first
by making the substitution $y_i = x_i - 1$.
resulting in the equation
$$y_1 + ... + y_k = n - k.$$QED

\bigskip
The number $(n + k - 1)!/(k -
1)! n!$ is the coefficient of $x^{n}$ in the
expansion of the function
\begin{eqnarray}{1 \over (1 - x)^{k}} =
\sum_{i=0}^{\infty} c_{n} x^{n}.\end{eqnarray}

\bigskip
Let $c_{n,k}$ be the number of elements of
${\cal S}_n$ consisting of exactly $k$ cycles.

\bigskip
\noindent
{\bf Theorem B8}~~{\it With the notations above
we have
$$c_{n,k} = (n - 1) c_{n-1,k} + c_{n-1, k-1}$$
and these numbers are the coefficients of the
expansion of the function $x(x + 1) ... (x + n
-1)$
$$x(x + 1) ... (x + n
-1) = \sum_{k=0}^{n} c_{n,k} x^{k}$$ and also
$${x! \over (x - n)!} = \sum_{k=0}^{n}
(-1)^{n-k} c_{n,k} x^{k}.$$ Moreover these
numbers are the coefficients of the expansion of
the function
$$\log (1 + x)^{k} = k! \sum_{n=k}^{\infty}
c_{n,k} {x^{n} \over n!}.$$}

\medskip
\noindent
{\it Proof.} The recursive relation follows from
the observation that there are exactly $n - 1$
different ways to get a permutation on
$n$ elements consisting of exactly $k$ cycles
from a permutation on
$n - 1$ elements consisting of exactly $k$
cycles. These account for the first term on the
right of the recursive formula. Next we observe
that there is exactly one way to get
a permutation on
$n$ elements consisting of exactly $k$ cycles
from a permutation on
$n - 1$ elements consisting of exactly $k - 1$
cycles and these account for the second term in
the formula. The rest of the Theorem follows
from the observation that if we write
$$x(x + 1) ... (x + n
-1) = \sum_{k=0}^{n} a_{n,k} x^{k}$$ then the
coefficients satisfy the same recursive formula
as $c_{n,k}$:
$$a_{n,k} = (n - 1) a_{n-1,k} + a_{n-1,k-1}.$$
The last assetion follows by observing that
$$(1 + x)^t = e^{t \log (1 + x)} =
\sum_{k=0}^{\infty} {1 \over k!} t^{k} (\log
(1 + x))^{k}.$$ On the other hand, we have
\begin{eqnarray*}
(1 + x)^t &=& \sum_{n=0}^{\infty} {t! \over n! (t
- n)!} x^{n} \\&=& \sum_{n=0}^{\infty} {x^{n}
\over n!} \sum_{j=0}^{n} c_{n,j} t^{j}\\&=&
\sum_{j=0}^{\infty} t^{j} \sum_{n=j}^{\infty}
c_{n,j} {x^{n}
\over n!}.
\end{eqnarray*}QED

\bigskip
Denote by $P(n, k)$ the set of all partitions of
a set of $n$ elements into $k$ non-empty subsets
and let
$$p_{n,k} = \# P(n, k).$$

\bigskip
\noindent
{\bf Theorem B9}~~{\it With the notations above
we have
$$p_{n,k} = k p_{n-1,k} + p_{n-1,k-1}$$and these
numbers are the coefficients of the expansion of
the function
$$x^{n} = \sum_{k=0}^{n} p_{n,k} {x! \over (x -
k)!}$$ and also as the coefficients of the power
series expansion
$$(e^{x} - 1)^{k} = k! \sum_{n=k}^{\infty}
p_{n,k} {x^{n} \over n!}.$$}

\medskip
\noindent
{\it Proof.}~~The recursive formula follows from
the observation that a partition of $n$ elements
into $k$ subsets can be obtained from a
partition of $n - 1$ elements into $k$ subsets
by insertingthe element $n$ into any one of the
$k$ subsets. Alternatively one can also get a
partition of $n$ elements into $k$ subsets from a
partition of $n - 1$ elements into $k - 1$
subsets by simply adding one more subset
consisting of just the element $n$.

\medskip
For a positive integer $x$ there are exactly
$x^{n}$ maps from the set $\{1, ..., n\}$ of
$n$ elements to the set $\{1, ..., x\}$. On
the other hand, by definition of the number
$p_{n,k}$ we have the relation:

\medskip
$k! p_{n,k} = \#$ of surjections from a set of
$n$ element onto a set of $k$ elements.

\medskip
\noindent
Hence for any subset $Y$ of $k$
elements of
$\{1, ..., x\}$ there are $k! p_{n,k}$
surjections from
$\{1, ..., n\}$ onto the set $Y$. Since the
number of subsets of $k$ elements of $\{1, ...,
x\}$ is $x!/k! (x-k)!$ we get
$$x^{n} = \sum_{k=0}^{n} {x! \over k! (x-k)!} k!
p_{n,k} = \sum_{k=0}^{n} {x! \over (x-k)!}
p_{n,k}.$$
By Corollary 4 we have:
\begin{eqnarray}k! p_{n,k} = \sum_{i=0}^{k}
(-1)^{i} {k!
\over i! (k - i)!} (k -
i)^{n} = \sum_{i=1}^{k}
(-1)^{k-i} {k!
\over i! (k - i)!} i^{n}.\end{eqnarray}

\bigskip
If $k = 1$ then $p_{n,1} = 1$ as there is only
one such partition. The usual expansion of the
exponential function yields
$$e^{x} - 1 = \sum_{n=1}^{\infty} {x^{n} \over
n!}.$$ The case of general $k$ can be verified
by induction by differentiating the
power series
$$F_k(x) = \sum_{n=k}^{\infty}
s_{n,k} {x^{n} \over n!}$$ resulting in
\begin{eqnarray*}F_k^{'}(x) &=&
\sum_{n=k}^{\infty} p_{n,k} {x^{n-1} \over (n -
1)!} \\&=&
\sum_{n=k}^{\infty} (k p_{n-1,k} + p_{n-1,k-1})
{x^{n-1} \over (n-1)!}\\&=&
kF_{k}(x) + F_{k-1}(x).\end{eqnarray*}Bu
induction hypothesis we have:
$$F_{k-1}(x) = {1 \over (k-1)!} (e^{x} -
1)^{k-1}$$ hence the function $F_k$ satisfies
the differential equation
$$F_k^{'}(x) = kF_{k}(x) + {1 \over
(k-1)!} (e^{x} - 1)^{k-1}.$$It is clear that
$$F_k(x) = {1 \over k!} (e^{x} -
1)^{k}$$is a solution and is indeed the uniqe
solution satisfying $p_{k,k} = 1$.~~QED

\bigskip
\noindent
{\bf Theorem B10}~~{\it
The number of partitions of $n$
$$p(n) = \sum_{k=1}^{n} {k^{n}
\over k!}.$$}

\bigskip
Denote by $p_{k}(n)$ the number of solutions of
the equation
\begin{eqnarray}x_1 + ... + x_k = n\end{eqnarray}
with the condition that
$1 \le x_k \le x_{k-1} \le ... \le x_1$. This
number is obviously equal to the number of
solutions of the equation
\begin{eqnarray}y_1 + ... + y_k = n -
k\end{eqnarray} with the condition that
$0 \le y_k \le y_{k-1} \le ... \le y_1$. If
there are exactly $i$ of the integers $y_i$
which are positive then these are the solutions
of $x_1 + ... + x_i = n - k$ and so there are
$p_{i}(n-k)$ of such solutions; consequently we
have:

\bigskip
\noindent
{\bf Theorem B11}~~{\it With the notations above
we have
$$p_{k}(n) = \sum_{i=1}^{k} p_{i}(n-k).$$}

\bigskip
Consider the case $k = 3$ then the number of
solutions of $$x_1 + x_2 + x_3 = n$$ such that
$0 \le x_3 \le x_2 \le x_1$ is the same as
$p_{3}(n+3)$. Let $y_1 = x_1 - x_2 \ge 0, y_2 =
x_2 - x_3 \ge 0, y_3 = x_3 \ge 0$ then this is
also the number of solutions of the equation
$$y_1 + 2y_2 + 3y_3 = n$$ with the condition
that $y_i \ge 0$. Thus the number $p_{3}(n+3)$
is the coefficient of $x^{n}$ in the expansion of
the function (compare ())
$$(1 - x)^{-1} (1 - x^{2})^{-1} (1 - x^{3})^{-1}
= \sum_{n=0}^{\infty} p_{3}(n+3)~ x^{n}.$$We
have the factorization
$$(1 - x^{3}) = (1 - x) (1 - \theta x) (1 -
\theta^{2} x)$$ where $\theta$ is a $3$-rd root
of unity, hence
\begin{eqnarray*}&&(1 - x)^{-1} (1 -
x^{2})^{-1} (1 - x^{3})^{-1} \\&&= (1 - x)^{-3}
(1 + x)^{-1} (1 - \theta x)^{-1} (1 -
\theta^{2} x)^{-1}\\&&= {1 \over 6(1
- x)^{3}} + {1 \over 4(1
- x)^{2}} + {17 \over 72 (1 - x)}
+ {1 \over 8 (1 + x)} + \\&&~~~~~~~~~~~~+ {1
\over 9 (1 -
\theta x)} + {1 \over 9 (1 - \theta^{2}
x)}\end{eqnarray*}and we get from the expansion
of each of the term of the partial fraction
decomposition that
$$p_{3}(n+3) = {(n+3)^{2} \over 12} - {7 \over
72} + {(-1)^{n} \over 8} + {\theta^{n} +
\theta^{2n} \over 9}.$$We infer that
$$|p_{3}(n+3) - {(n+3)^{2} \over 12}| < {1
\over 2}$$ or equivalently that
$$|p_{3}(n) - {n^{2} \over 12}| < {1 \over
2}.$$

\bigskip
The following identity is easily established by
induction:

\bigskip
\noindent
{\bf Theorem B12}~~{\it The number
$p_{k}(n)$ satisfies the following recursive
relation: $p_k(n) = p_{k-1}(n-1) +
p_k(n-k)$.}

\bigskip
Obviously we have $p_1(n) = n$ and $p_2(n) =
n/2$ or $(n - 1)/2$ according to $n$ being even
or odd. Thus Theorem 10 yields $p_3(n) = p_2(n-1)
+ p_3(n - 3)$, $p_4(n) = p_3(n-1) + p_4(n - 4),$
$p_5(n) = p_4(n-1) + p_5(n - 5)$
and we get for example
$$p_1(6) = 1, p_2(6) = 3, p_6(6) = 1$$
$$p_3(6) = p_2(5) + p_3(3) = 3,$$
$$p_4(6) = p_3(5) = p_2(4) = 2,$$
$$p_5(6) = p_4(5) = p_3(4) = p_2(3) = 1$$hence
$$p(6) = \sum_{k=1}^{6} p_{k}(6) = 1 + 3 + 3 + 2
+ 1 + 1= 11.$$

\bigskip
The total partition length $L(n)$ of a positive
integer
$n$ is defined to be
\begin{eqnarray}
L(n) = \sum_{k=1}^{n} k p_k(n).
\end{eqnarray}
For example if $n = 6$ then $L(6) = 1 + 6 + 9 +
8 + 5 + 6 = 35.$

\bigskip
For $n = 7$ we have
$$p_1(7) = 1, p_2(7) = 3, p_7(7) = 1$$
$$p_3(7) = p_2(6) + p_3(4) = p_2(6) + p_2(3)
= 4,$$
$$p_4(7) = p_3(6) = 3,$$
$$p_5(7) = p_4(6) = 2,$$
$$p_6(7) = p_5(6) = 1$$hence
$$p(7) = \sum_{k=1}^{7} p_{k}(7) = 1 + 3 + 4 + 3
+ 2 + 1 + 1 = 15$$and the total partition length
$$L(7) = 1 + 6 + 12 + 12 + 10 + 6 + 7 = 54.$$

\bigskip
For general $k$ one has the following
asymptotic formula:

\bigskip
\noindent
{\bf Theorem B13}~~{\it For $n \rightarrow
\infty$ the number $p_k(n)$ is asymptotically
given by:
$$p_k(n) \sim {n^{k-1} \over (k - 1)! k!}.$$}

\medskip
\noindent
{\it Proof.}~~The number $p_k(n)$ is defined to
be the number of solutions of
$x_1 + ... + x_k = n$ with the condition that
$1 \le x_k \le x_{k-1} \le ... \le x_1$. If we
drop this last condition then the $k!$
permutations of a solution is also a solution of
$x_1 + ... + x_k = n$. However since $x_i$ may
equal $x_j$ for $i \ne j$ hence we have the
inequality:
$$C_{k-1}^{n-1} = {(n-1)! \over (k-1)! (n -
k)!} \le k!~ p_{k}(n).$$

\medskip
On the other hand, if we set $y_i = x_i + (k -
i)$ and if $x_1, ..., x_k$ is a solution with
$1 \le x_k \le x_{k-1} \le ... \le x_1$ then the
$y_i$'s are distinct and is a solution of the
equation:
$$y_1 + ... + y_k = n + {k(k-1) \over 2}.$$ From
this we obtain a reverse inequality:
$$k!~ p_{k}(n) \le C_{k-1}^{n+(k(k-1)/2) - 1} =
{\{n + (k(k-1)/2) - 1\} ! \over (k-1)! \{n +
k(k-1)/2) - k\}!}.$$The Theorem follows
immediately from these two estimates.~~QED

\bigskip
\bigskip
\noindent
{\bf References}

\bigskip
\noindent
[A-S]  Azukawa, K. and Suzuki, M.: {\sl Some examples of
algebraic geometry and hyperbolic
manifolds.} Rocky Mountain J. Math. {\bf 10},
653-659  (1980).

\bigskip
\noindent
[B] Bloch, A.: {\it Sur les systemes de
donctions uniformes satisfaisant a l'eauation
d'une vatiete algebraique dont l'irregularite
depasse la dimension.}  J. de Math. {\bf 5}
(1926), 19-66

\medskip
\noindent
[B1]  Bogomolov, F.: {\sl Families of curves on a surface of
general type.}  Sov. Math. Dokl. {\bf 18}, 1294-1297
(1977)

\medskip
\noindent
[B2]  Bogomolov, F.: {\sl Holomorphic tensors and vector
bundles on projective varieties.}
Math. USSR, Izv. {\bf 13}, 499-555 (1979)

\medskip
\noindent
[B-K] Brieskorn, E., Knorrer, H.: {\sl Ebene algebraische
Kurven.} Birkhauser, Basel-Boston-
Stuttgart (1981)

\medskip
\noindent
[B-P-V]  Barth, W., Peters, C., Van de Ven, A.: {\sl Compact
complex surfaces.} Springer-Verlag
Berlin-Heidelberg-New York-Tokyo (1984)

\medskip
\noindent
[B-R] Beltrametti, M. and Robbiano, L.: {\it
Introduction to the theory of weighted projective
spaces.} Expositiones Mat.

\medskip
\noindent
[B-S] Banica, C. and O. Stanasila :
{\it Algebraic methods in the global theory of
complex spaces.} English translation (1976, John
Wiley \& Sons Publisher) of the original Romanian
book of the same title published in 1974

\medskip
\noindent
[C-G] Carlson, J.A. and Green, M.: {\sl A Picard Theorem for
Holomorphic Curves in the Plane.}
Duke Math. J. , 1-9 (1976)

\medskip
\noindent
[D] Demailly, J.P.: {\it Algebraic
criteria for Kobayashi hyperbolic projective
varieties and jet differentials.} Preprint 1999

\medskip
\noindent
[D-El] Demailly, J.P. and J. El Goul: {\it
Hyperbolicity of generic surfaces of high degree
in projective 3-space.} Preprint 1998

\medskip
\noindent
[El] El Goul, J.: {\it Surfaces
alg\'{e}briques hyperboliques, Propri\'{e}t\'{e}s
de n\'{e}gativit\'{e} de la courbure.} Th'{e}se
de doctorat de l'ubiversit\'{e} Fourier, pp 78 (1997)

\medskip
\noindent
[D-S-W1] Dethloff, G., Schumacher, G., Wong, P. M.:
{\sl Hyperbolicity of the complements of
plane algebraic curves.} Amer. J. Math. {\bf 117},
(1995) 573-599

\medskip
\noindent
[D-S-W2] Dethloff, G., Schumacher, G., Wong, P. M.:
{\sl Hyperbolicity of the complements of
plane algebraic curves: The Three Components Case.} Duke
Math. J. {\bf 78}, (1995)

\medskip
\noindent
[Di] Dimca, A.: {\sl Singularities and
topology of hypersurfaces.} Universitext,
Springer-Verlag 1992

\medskip
\noindent
[Do 1989] Dolgachev, I.: {\sl Weighted projective
varieties.} In {\sl Group actions and vector
fields, Proceedings 1981.} Lecture Notes in Math.
956, Springer-Verlag (1982), 34-71

\medskip
\noindent
[F 1991] Faltings, G.:
{\it Diophantine approximation on abelian
varieties.}  Ann. Math., {\bf 133} (1991), 549-576

\medskip
\noindent
[F-H 1991]
Fulton, W.
and J. Harris,  {\sl Representation Theory.}
Grad. Text in Math. {\bf 129}, Springer-Verlag
(1991)

\medskip
\noindent
[G1] Grant, C.G.:  {\sl Entire holomorphic curves in
surfaces.} Duke Math. J., {\bf 53},  345-358 (1986)

\medskip
\noindent
[G2] Grant, C.G.: {\sl Hyperbolicity of surfaces modular
rational and elliptic Curves.}
{\sl Pacific J. Math.} {\bf 139}, 241-249 (1989)

\medskip
\noindent
[G-G]  Green, M. and Griffiths, P.: {\sl
Two applications of algebraic geometry to entire
holomorphic mappings.}  The Chern Symposium 1979
(Proc. Interat. Symp., Berkeley, Calif., 1979)
Springer-Verlag, New York (1980), 41-74

\medskip
\noindent
[G-P] Grauert, H. and Peternell, U.: {\sl Hyperbolicity of
the Complement of Plane Curves.}
Manus. Math. {\bf 50}, 429-441 (1985)

\medskip
\noindent
[Gr]  Grauert, H.: {\sl Jetmetriken und hyperbolische
Geometrie.}  Math. Z. {\bf 200}, 149-168 (1989)

\medskip
\noindent
[Gre]  Green, M.: {\sl Some Picard theorems for holomorphic
maps to algebraic varieties.}
Am. J. Math. {\bf 97}, 43-75 (1975 )

\medskip
\noindent
[Gri]  Griffiths, P.: {\sl Hermitian differential geometry,
Chern classes, and positive vector
bundles.} In Global Analysis (papers in honor of K. Kodaira)
edited by D. C. Spencer and
S. Iyanaga, University of Tokyo Press and Princeton University Press
(1969), 185-251

\medskip
\noindent
[Gro] Grothendieck, A.: {\sl La theorie des classes de
Chern.} Bull. Soc. math. France {\bf 86}, 137-
154 (1958)

\medskip
\noindent
[Ha] Hartshorne, R.: Algebraic Geometry.
Graduate Text in Math. {\bf 52}, Springer-Verlag
(1977)

\medskip
\noindent
[Hi] Hirzebruch, F.: {\sl Topological Methods in Algebraic
Geometry.}  Grund. der Math. Wiss. {\bf 131}
Springer-Verlag (1966)

\medskip
\noindent
[H-W 1954] Hardy, G.H. and Wright, E.M.: {\sl
An introduction to the theory of numbers.}
3rd. edition (1970) Clarendon Press

\medskip
\noindent
[I1]  Iitaka, S.: {\sl Logarithmic forms of algebraic
varieties.} J. Fac. Sci. Univ. Tokyo, Sect.
IA Math. 23, 525-544 (1976)

\medskip
\noindent
[I2] Iitaka, S.:{\sl  Birational geometry for open
varieties.} Sem. Math. Sup., Les Presses de
L'Universite de Montreal

\medskip
\noindent
[J] Jung, E. K.: {\sl Holomorphic Curves in Projective
Varieties.} Notre Dame Thesis, 1995

\medskip
\noindent
[Ka] Kawasaki, T.: {\sl Cohomology of twisted projective
spaces and lens complexes.}  Math. Ann. {\bf 206}, 243-248
(1973)

\medskip
\noindent
[Ko] Kobayashi, S.: {\sl Hyperbolic manifolds and
holomorphic mappings.}  Marcel Dekker,
New York (1970)

\medskip
\noindent
[K-O] Kobayashi, S. and Ochiai, T.:  {\sl On complex
manifolds with positive tangent bundles.}
J. Math. Soc. Japan {\bf 22}, 499-525 (1970)

\medskip
\noindent
[L] Lu, S.: {\sl On meromorphic maps into varieties of
log-general type.} Proc. Symp. Amer.
Math. Soc. {\bf 52}, 305-333 (1991)

\medskip
\noindent
[L-Y]  Lu, S., Yau, S. T.: {\sl Holomorphic curves in
surfaces of general type.} Proc. Natl.
Acad. Sci. USA {\bf 87}, 80-82 (1990)

\medskip
\noindent
[Ma]  Maruyama:  {\sl The Theorem of
Grauert-Mulich-Splindler.}  Math. Ann. {\bf 255}, 317-333
(1981)

\medskip
\noindent
[Mc]  McQuillan, M.: {\sl A new proof of the Bloch
Conjecture.} preprint June 1993

\medskip
\noindent
[Mi]  Miyaoka, Y.: {\sl On the Chern numbers of surfaces of
general type.} Invent. Math. {\bf 42},
225-237 (1977)

\medskip
\noindent
[M-N] Masuda, K. and Noguchi, J.: {\sl  A Construction of
Hyperbolic Hypersurface of ${\bf P}^n({\bf C}).$}
Preprint 1994

\medskip
\noindent
[Nad] Nadel, A.: {\sl Hyperbolic surfaces in ${\bf P}^3$.}
Duke Math. J. {\bf  58}, 749-771 (1989)

\medskip
\noindent
[Nam]  Namba, M., {\sl Geometry of Projective Algebraic
Curves.} Monographs and Textbooks
in Pure and Appl. Math. {\bf 88} Marcel Dekker, New York
(1984)

\medskip
\noindent
[O-S-S 1980] Okonek, C., Schneider and Spindler,
H.: {\sl Vector bundles on complex projective
spaces.} Progress in Math. {\bf 3},
Birkh\"{a}user 1980

\medskip
\noindent
[Sa]  Sakai, F.: {\sl Semi-stable curves on algebraic
surfaces.} Math. Ann. {\bf 254}, 89-120 (1980)

\medskip
\noindent
[St]  Stoll, W.: {\it Value Distribution Theory in Several
Complex Variables.} Shandong Science and Technology Press
(1995)

\medskip
\noindent
[S-W]  Stoll, W. and P. M. Wong:  {\it Second Main
Theorem of Nevanlinna Theory for
Non-equidimen-sional Meromorphic Maps.}  Amer. J.
Math., vol. 116 (1994), 1031-1071

\medskip
\noindent
[S-Y1]  Siu, Y.T. and Yeung, S.K.: {\it Defect for ample
divisors of abelian varieties, Schwarz Lemma, and hyperbolic
hypersurfaces of low degrees.} (1996)

\medskip
\noindent
[S-Y2]  Siu, Y.T. and Yeung, S.K.:  {\sl Hyperbolicity of the
Complement of A Generic Smooth
Curve of High Degree in the Complex Projective Plane.}
Preprint 1994, June.

\medskip
\noindent
[T1]  Tsuji, H.: {\sl An inequality of Chern numbers for open
algebraic varieties.} Math. Ann.
{\bf 277}, 483-487 (1987)

\medskip
\noindent
[T2] Tsuji, H.: {\it Stability of
tangent bundles of minimal algebraic varieties.}
Topology {\bf 22} (1988), 429-441

\bigskip
\noindent
[W1 1989]
 Wong, P.
M.,  On the Second Main Theorem of Nevanlinna
Theory.  Amer. J. of Math., vol. 111 (1989),
549-583.

\medskip
\noindent
[W2] Wong, P. M., Holomorphic curves in spaces of constant
curvature. In "Complex
Geometry" edited by G. Komatsu and Y. Sakane,
201-223. Lecture Notes in Pure and
Applied Mathematics, vol 143, Marcel Dekker (1993)

\medskip
\noindent
[W3] Wong, P. M.: {\it Recent results in
Hyperbolic Geometry and Diophantine Geometry.}
Proc. Inter-national Symposium on Holomorphic
mappings, Diophantine Geometry and Related Topics
(1992), edited by J. Noguchi, R.I.M.S., Kyoto
University (1993), 120-135

\medskip
\noindent
[W4] Wong, P. M.: {\it Holomorphic curves in varieties with nef jet
differentials.} preprint 1999 (revision of an earlier version,
1995)

\medskip
\noindent
[W5] Wong, P. M.: {\it Topics on Higher Dimensional
Nevanlinna Theory and Intersection Theory.}  Six lectures at
the Center of Theoretical Science, National Tsing-Hua
University, 1999

\medskip
\noindent
[W6] Wong, P. M.: {\it On Complex hyperbolic manifolds.} preprint
2000

\medskip
\noindent
[W7] Wong, P. M. and Cherng-Yih Yu: {\it Negatively Curved
Jet Metrics On Quasi-Projective Varieties.} preprint 2000

\medskip
\noindent
[X]  Xu, G.: {\sl Subvarieties of general hypersurfaces in
projective space.} J. Diff.
Geom. (1992)

\medskip
\noindent
[Z] Zaidenberg, M.G.: {\sl Stability of hyperbolic
embeddedness and construction of examples} Math. USSR Sbornik
{\bf 63}, 351-361 (1989)

\bigskip
\bigskip
\bigskip
\noindent
e-mail address: wong.2@nd.edu or ppmpwong@juno.com

\end{document}